\documentclass[12pt]{article}
\usepackage{latexsym}
\usepackage{amsmath,amssymb}
\usepackage{bm}
\usepackage[dvipdfmx]{graphicx,color}
\usepackage{amssymb}
\usepackage{cite}
\usepackage{amsfonts, amscd, amsthm}
\usepackage[all]{xy}
\usepackage{here}
\usepackage{ulem}
\usepackage{mathabx}

\newcommand{\DDelta}{\widehat{\Delta}}
\newcommand{\vvarphi}{\widehat{\varphi}}
\renewcommand{\theequation}{\thesection.\arabic{equation}}

\newcommand{\dd}{{\rm d}}

\makeatletter
\@addtoreset{equation}{section}
\renewcommand{\theequation}{\thesection.\@arabic\c@equation}
\makeatother

\makeatletter
\renewcommand\appendix{\par%\newpage
  \setcounter{section}{0}%
  \setcounter{subsection}{0}%
  \gdef\thesection{Appendix \@Alph\c@section }
  \renewcommand{\theequation}
  {\Alph{section}.\arabic{equation}}
}
\makeatother

\makeatletter
\newcounter{subeqncnt}
\def\thesubeqncnt{\alph{subeqncnt}}
\def\subequations{\begingroup%
\stepcounter{equation}\edef\@tempa{\theequation}%
\let\c@equation\c@subeqncnt\c@subeqncnt\z@
\edef\theequation{\@tempa\noexpand\thesubeqncnt}}

\makeatother
\allowdisplaybreaks

\begin{document}

\titlepage

%%%%%
\title{Differential Equations \\of General Genus 
Hyperelliptic $\wp$ Functions\\
} 
\author{
Kazuyasu Shigemoto\thanks{shigemot@tezukayama-u.ac.jp} \\
Tezukayama University, Nara 631-8501, Japan\\
}
\date{\empty}

%%%%%%%%%%%%%%%%%%%%%%%%%%%%%%%%%%%%%%%%

\maketitle
\abstract{%
We review the Baker's method to obtain differential equations of the 
general genus hyperelliptic $\wp$ functions. Further, we demonstrate to 
obtain differential equations of genus four hyperelliptic 
differential equations, which agree with our previous result.
}

%%%%%%%%%%%%%%%%%%%%%%%%%%%%%%%%%%%%%%%%%%%%%%%%%
\section{Introduction} 
\setcounter{equation}{0}
Hyperelliptic $\wp$ differential equations are considered as the 
higher dimensional generalization of the KdV equation, that is, 
integrable differential equations. Then, there exists the potential, the sigma
function, in the form 
$\displaystyle{\wp_{ij}(u_1,u_2, \cdots, u_g)=-\frac{\partial }{\partial u_i}
\frac{\partial}{\partial u_i} \log \sigma(u_1,u_2, \cdots, u_g)}$. 
The sigma function corresponds to the $\tau$ function in the 
KdV equation.

There are three kinds of approach to derive differential equations
for the hyperelliptic $\wp_{ij}$ functions. The first approach is 
to use the generalized fundamental relation, the Bolza's formula (II) by Baker~\cite{Baker4}.
This approach is the systematic approach to derive differential equations, 
and we can understand the structure of differential equations.
The second approach is to use one series of differential equations 
$\wp_{gggi}=6 \wp_{gg} \wp_{gi} +6 \wp_{g, i-1}-2 \wp_{g-1,i}
+\lambda_{2g} \wp_{g i}+\delta_{g,i} \lambda_{2g-1}/2$ for general genus $g$ case 
and the Kummer surface type relations by 
Buchstaber-Enolskii-Leykin~\cite{Buchstaber1,Buchstaber2}. The Kummer surface type
relations are necessary to obtain the full differential equations 
because one series of differential equations contain only
$\lambda_{2g}$ and $\lambda_{2g-1}$ but does not contain $\lambda_i\ (i=0,\cdots, 2g-2)$.
In genus two case, we use the Kummer surface identities itself.
This approach is quite elegant, but it is unclear how to systematically obtain the 
full differential equations.
The third approach is to use the defining relation of the $\zeta_i$ functions
by ours~\cite{Hayashi}.  This approach is quite straightforward but 
it is difficult to understand the structure of differential equations.

Thus, we consider that the Baker's paper is quite important, but it is difficult to follow.
Hence, we review the Baker's method and give the detailed derivation 
of various formulae. As Baker demonstrated his method in genus three case,  we 
demonstrate his method in genus four case and we obtain the same result as our 
previous paper~\cite{Hayashi}.
%\newpage

%%%%%%%%%%%%%%%%%%%%%%%%%%%%%%%%%%%%%%%%%%%%%%%%%%%%%%%%%%%
\section{Review of the Baker's formulation}
\setcounter{equation}{0}
We first review the basic tools, Jacobi's inversion problem, various formulae of $\wp_{ij}$,
and various formulae of $\zeta_i$ functions. Next, we review the fundamental relation, 
which is the starting point of the Baker's paper~\cite{Baker4}. We also review to rewrite the 
fundamental relation into the Balza's formula (I) and Bolza's 
formula (II)\cite{Baker1,Bolza1,Bolza2}.
Actually, Baker starts from the fundamental relation in the form of the Bolza's formula (II).
%%%%%%%%%%%%%%%%%%%%%%%%%%%%%%%%%%%%%%%%%%%%%%%%%%
\subsection{Basic tools}
We introduce $F(x)$ by
\begin{align}
F(x)=F(x; x_1, \cdots, x_g)=\prod_{i=1}^g (x-x_i)=x^g+\sum_{j=1}^g (-1)^j h_j(x_1,x_2,\cdots, x_g) x^{g-j}.
\label{2e1}
\end{align}
where $h_j(x_1, \cdots, x_g)$ is the $j$-th fundamental symmetric polynomial basis of 
$\{x_1, \cdots, x_g\}$.
By substituting $x=x_j$, we have the relation
\begin{equation}
x_i^g    -h_1(x_1,x_2,\cdots,x_g)x_i^{g-1}
         +h_2(x_1,x_2,\cdots,x_g)x_i^{g-2}
+\cdots
+(-1)^{g} h_g(x_1,x_2,\cdots,x_g)=0 .
\label{2e2}
\end{equation}
By the partial fraction theorem, we can express  $x^k/F(x)\ $$(0\le k \le g-1)$ in the
form  
\begin{align}
\frac{x^k}{F(x)}=\sum^g_{i=1} \frac{a_i}{(x-x_i)},\ (0\le k \le g-1), 
\label{2e3}
\end{align}
where $a_i$ are constants. In order to determine these constant values, we multiply
$(x-x_j)$ in Eq.(\ref{2e3}) and take the limit $x \rightarrow x_j$. Then we have

\begin{align}
\frac{x_j^k}{F'(x_j)}=\sum^g_{i=1} a_i \delta_{i,j}=a_j, 
\label{2e4}
\end{align}
where we denote $F'(x_j)$ as $\displaystyle{F'(x_j)=\frac{\dd F(x)}{\dd x}\Big|_{x=x_j}}$. 
Thus we have
\begin{align}
\frac{x^k}{F(x)}=\sum^g_{i=1} \frac{x_i^k}{F'(x_i)} \frac{1}{(x-x_i)},\ (0\le k \le g-1) .
\label{2e5}
\end{align}
Further, we multiply $x$ and take the limit of $x \rightarrow \infty$ in Eq.(\ref{2e5}),
we obtain the useful formula
\begin{align}
\sum^g_{i=1} \frac{x_i^k}{F'(x_i)} =\delta_{k, g-1},  \ (0\le k \le g-1) .
\label{2e6}
\end{align}

Next, we introduce $\chi_{g-j}(x_i; x_1,x_2,\cdots, x_g)$ in the form
\begin{align}
\chi_{g-j}(x;x_1, \cdots, x_g) = &\ x^{g-j}-h_1(x_1, \cdots, x_g) x^{g-j-1} \nonumber\\
&+h_2(x_1, x_2, \cdots, x_g) x^{g-j-2}+\cdots+(-1)^{g-j} h_{g-j}(x_1, \cdots, x_g)  ,
\label{2e7}
\end{align}
As the special case, $\chi_{g}(x;x_1, \cdots, x_g)=F(x)$ and $\chi_{0}(x;x_1, \cdots, x_g)=1$.  
By the definition of $\chi_{g-i}$, we can show 
\begin{equation}
\sum^{g}_{k=1} (x^k-x^{k-1} x_r) \chi_{g-k}(x_r; x_1, x_2, \cdots, x_g)=F(x) . 
\label{2e8}
\end{equation}
Then we obtain
\begin{equation}
\frac{F(x)}{x-x_r}=\sum^{g}_{k=1} x^{k-1}\chi_{g-k}(x_r; x_1, x_2, \cdots, x_g)  , 
\label{2e9}
\end{equation}
which will be used in Eq.(\ref{3e8}).

Further we obtain 
\begin{equation}
\frac{\partial }{\partial x} \chi_{g-i} (x_1;x, x_2, x_3, \cdots,x_g)
=-\chi_{g-i-1} (x_1;x_2,x_3, \cdots,x_g) , 
\label{2e10}
\end{equation}
by using 
$\displaystyle{ \frac{\partial }{\partial x} h_i (x, x_2, x_3, \cdots, x_g)
=h_{i-1}(x_2, x_3, \cdots, x_g) }$.  This will be used in Eq.(\ref{2e60}).
 
Next, we prove the following relation, which is useful in the Jacobi's inversion problem
\begin{equation}
\sum^{g}_{i=1} \frac{x_i^{g-j} \chi_n(x_i; x_1,x_2,\cdots, x_g)}{F'(x_i)}=\delta_{n, j-1},\ 
(1\le j \le g, \ 0 \le n \le g) .
\label{2e11}
\end{equation}
[Proof]\\
We fix $j$ in the region $1\le j \le g$.\\
i)\ $0 \le n \le j-1$ case:\\
\begin{align}
&\sum^{g}_{i=1} \frac{x_i^{g-j} \chi_n(x_i; x_1,x_2,\cdots, x_g)}{F'(x_i)}
=\sum^{g}_{i=1} \frac{x_i^{g-j} (x_i^n-h_1 x_i^{n-1}+\cdots+(-1)^n h_n)}{F'(x_i)}
\nonumber\\
&=\sum^{g}_{i=1} \frac{x_i^{g-(j-n)} -h_1 x_i^{g-(j-n+1)}+\cdots+(-1)^n h_n x_i^{g-j})}{F'(x_i)}
=\delta_{n, j-1},\ (1\le j \le g), 
\label{2e12}
\end{align}
where we used Eq.(\ref{2e6}).\\
ii)\ $g \ge n \ge j \ge 1$ case:\\
\begin{align}
&\sum^{g}_{i=1} \frac{x_i^{g-j} \chi_n(x_i; x_1,x_2,\cdots, x_g)}{F'(x_i)}
=\sum^{g}_{i=1} \frac{x_i^{g-(j-n)} -h_1 x_i^{g-(j-n+1)}+\cdots+(-1)^n h_n x_i^{g-j})}{F'(x_i)}
\nonumber\\
&=\sum^{g}_{i=1} \frac{-(-h_1 x_i^{g-(j-n+1)}+\cdots+(-1)^{g} h_g x_i ^{n-j})
-h_1 x_i^{g-(j-n+1)}+\cdots+(-1)^n h_n x_i^{g-j})}{F'(x_i)}
\nonumber\\
&=\left\{ \begin{array}{cc}
  0,   & \mbox{(if $n=g$)} , \\
 \displaystyle{-\sum^{g}_{i=1} \frac{(-1)^{n+1} h_{n+1} x_i^{g-(j+1)}+\cdots+(-1)^g h_g x_i ^{n-j})}{F'(x_i)}=0}, 
& \mbox{( if $n\le g-1$)} ,  
\end{array}
\right.
\label{2e13}
\end{align}
where we used Eq.(\ref{2e2}) and Eq.(\ref{2e6}). This completes the proof of Eq.(\ref{2e11}).

Further, we give the useful relation to obtain the $\wp_{g j}$$(1\le j \le g)$ of the form
\begin{equation}
\sum^{g}_{i=1} \frac{x_i^{g} \chi_n(x_i; x_1,x_2,\cdots, x_g)}{F'(x_i)}=(-1)^n h_{n+1}(x_1, x_2, \dots, x_g), \
(0 \le n \le g-1) .
\label{2e14}
\end{equation}
[Proof]\\
Using Eq.(\ref{2e2}), we have
\begin{align}
x_i^n-h_1 x^{n-1} +\cdots+(-1)^{n} h_{n} +(-1)^{n+1} \frac{h_{n+1}}{x_i}
+\cdots+(-1)^{g} \frac{h_g}{x_i^{g-n}} =0 . 
\label{2e15}
\end{align}
Then we have
\begin{align}
&\chi_n(x_i; x_1, x_2, \cdots, x_g)=x_i^n-h_1 x^{n-1} +\cdots+(-1)^{n} h_{n}
\nonumber\\
&=-\Big(\frac{(-1)^{n+1}h_{n+1}}{x_i}+\cdots+(-1)^{g} \frac{h_g}{x_i^{g-n}} \Big) .
\label{2e16}
\end{align}
Thus, we obtain 
\begin{align}
&\sum^{g}_{i=1} \frac{ x_i^g \chi_n(x_i; x_1, x_2, \cdots, x_g)}{F'(x_i)}
=-\sum^{g}_{i=1} \frac{(-1)^{n+1} h_{n+1} x_i^{g-1}+\cdots+(-1)^{g} h_g x_i^{g-(g-n)}}{F'(x_i)}
\nonumber\\
&=(-1)^n h_{n+1}, \quad (0 \le n \le g-1), 
\label{2e17}
\end{align}
where we used Eq.(\ref{2e6}). This completes the proof of Eq.(\ref{2e14}).
%%%%%%%%%%%%%%%%%%%%%%%%%%%%%%%%%%%%%%%%%%%%%%%%%%%%%%%
\subsection{Jacobi's inversion problem, $(-\zeta_j)$ function and $\wp_{ij}$ function}
We summarize the formulation of hyperelliptic $\wp_{ij}$ function according to 
Baker's work~\cite{Baker1,Baker2,Baker4}. We consider the genus $g$ hyperelliptic curve
\begin{equation}
C:\quad y_i^2=\sum_{k=0}^{2g+1} \lambda_k x_i^k , 
\quad  (i=0, 1, 2, \cdots, g). 
\label{2e18}
\end{equation}
The Jacobi's inversion problem consists of solving following equations
\begin{equation}
\dd u_1    =\sum_{i=1}^g \frac{         \dd x_i}{y_i}, \quad 
\dd u_2    =\sum_{i=1}^g \frac{x_i      \dd x_i}{y_i}, \quad \cdots,\quad  
\dd u_{g-1}=\sum_{i=1}^g \frac{x_i^{g-2}\dd x_i}{y_i}, \quad 
\dd u_g    =\sum_{i=1}^g \frac{x_i^{g-1}\dd x_i}{y_i}.
\label{2e19}
\end{equation} 
From these equations, we have 
\begin{align}
&\sum_{i=1}^g \frac{x_i^{k-1}}{y_i} \frac{\partial x_i}{\partial u_j}
=\frac{\partial u_k}{\partial u_j}=\delta_{k,j} .
\label{2e20}
\end{align}
If we compare the expression, which is obtained from Eq.(\ref{2e11}) by replacing
$j \rightarrow g-k+1$ and $n \rightarrow g-j$ in the form
\begin{equation}
\sum_{i=1}^g \frac{x_i^{k-1} \chi_{g-j}(x_i; x_1, x_2, \cdots, x_g)}{F'(x_i)}
=\delta_{k, j} \quad (1\leq j \leq g) ,
\label{2e21}
\end{equation}
we obtain the useful expression
\begin{equation}
\frac{\partial x_i}{\partial u_j}=\frac{y_i \chi_{g-j}\left(x_i; x_1,x_2,\cdots, x_g\right)}{F'(x_i)} .
\label{2e22}
\end{equation}
The $\zeta_j$ functions are given from the hyperelliptic curve in the following way~\cite{Baker1} 
\begin{equation}
\hskip -7mm \dd(-\zeta_j)=
\frac{1}{4}\sum_{i=1}^g\frac{\dd x_i}{y_i}
\sum_{k=j}^{2g+1-j} (k+1-j) \lambda_{k+1+j} x_i^k
-\frac{1}{2} \dd\left(\sum_{i=1}^g\frac{y_i \chi_{g-j-1}(x_i; x_1, \cdots, \widecheck{x}_i, \cdots,x_g)}{F'(x_i)}\right),
\label{2e23}
\end{equation} 
where $\widecheck{x}_j$ denotes that the $x_j$ variable is missing. 
These $\zeta_j(u_1,u_2,\cdots, u_g)$ satisfy the integrability condition 
\begin{equation}
 \frac{\partial \left(-\zeta_j(u_1,u_2,\cdots, u_g)\right)}{\partial u_k}
=\frac{\partial \left(-\zeta_k(u_1,u_2,\cdots, u_g)\right)}{\partial u_j}  .
\label{2e24}
\end{equation} 
In the Baker's textbook~\cite{Baker1}, the expression of the second term of the r.h.s of Eq.(\ref{2e23}) is misleading.
$\wp_{jk}(u_1,u_2,\cdots, u_g)$ functions are given from the above $\zeta_j(u_1,u_2,\cdots, u_g)$ functions in the form 
\begin{equation}
\wp_{jk}(u_1,u_2,\cdots, u_g)=\wp_{kj}(u_1,u_2,\cdots, u_g)
=\dfrac{\partial \left(-\zeta_j(u_1,u_2,\cdots, u_g)\right)}{\partial u_k}.
\label{2e25}
\end{equation}
These $\zeta_j$, $\wp_{jk}$ and $\wp_{j k l m}$ are given by the hyperelliptic $\sigma$ function in the form%
$$
(-\zeta_j)=\dfrac{\partial (-\log \sigma )}{\partial u_j}, \quad  
\wp_{jk}=\dfrac{\partial^2 (-\log \sigma )}{\partial u_j \partial u_k}, \quad  
\textrm{and} \quad 
\wp_{j k l m}=\dfrac{\partial^4(-\log \sigma )}
{\partial u_j \partial u_k \partial u_l \partial u_m}, \qquad \textrm{etc.}.
$$

For the Weierstrass type, i.e.,\,$\lambda_{2g+2}=0$, we have 
$\displaystyle{\dd(-\zeta_g)=\frac{\lambda_{2g+1}}{4}\sum_{i=1}^{g} \frac{x_i^g\dd x_i }{y_i}}$, which gives
\begin{align}
&\wp_{g  }(u_1,u_2,\cdots, u_g)=\frac{\lambda_{2g+1}}{4} h_1(x_1,x_2,\cdots,x_g), 
\label{2e26}\\
&\wp_{g,g-1}(u_1,u_2,\cdots, u_g)=-\frac{\lambda_{2g+1}}{4} h_2(x_1,x_2,\cdots,x_g), 
\label{2e27}\\
&\ \vdots \nonumber\\
&\wp_{g 1  }(u_1,u_2,\cdots, u_g)=  (-1)^{g-1} \frac{\lambda_{2g+1}}{4} h_g(x_1,x_2,\cdots,x_g), 
\label{2e28}
\end{align}
by using Eq.(\ref{2e22}) and Eq.(\ref{2e14}) in the form
\begin{equation}
\sum_{i=1}^g \frac{x_i^{g} \chi_{g-j}(x_i;x_1, x_2, \cdots, x_g)}{F'(x_i)} 
=(-1)^{g-j} h_{g-j+1}(x_1, x_2,\cdots,x_g)   .
\label{2e29}
\end{equation}
Then, by using Eq.(\ref{2e2}),  we obtain 
\begin{equation}
\frac{\lambda_{2g+1}}{4} x_i^g=\sum_{j=1}^{g} \wp_{g j} x_i^{j-1} 
     =\wp_{gg   } x_i^{g-1}
     +\wp_{g,g-1} x_i^{g-2}
     +\cdots
     +\wp_{g2} x_i
     +\wp_{g1}.
\label{2e30}
\end{equation}
In this way, we have $\displaystyle{\dd(-\zeta_g)=\sum_{j=1}^{g} \wp_{g j}\dd u_j}$. For other $\wp_{i j}$, we 
must use $\zeta_j$, which satisfies the integrability condition Eq.(\ref{2e24}).  

%%%%%%%%%%%%%%%%%%%%%%%%%%%%%%%%%%%%%%%%%%%%%%%%%%%%
\subsection{Derivation of the fundamental relation}
Hyperelliptic $\wp$ differential equations are considered as the 
higher dimensional generalization of the KdV equation, that is, 
integrable differential equations. Then, there exists the potential, the sigma
function, in the form 
$\displaystyle{\wp_{ij}(u_1,u_2, \cdots, u_g)=-\frac{\partial }{\partial u_i}
\frac{\partial}{\partial u_i} \log \sigma(u_1,u_2, \cdots, u_g)}$,
The sigma function corresponds to the $\tau$ function in the 
KdV equation.

In order to obtain this $\sigma$ function, it is useful to use the 
differential form of the third kind.
We introduce  the hyperelliptic functions 
$\displaystyle{y^2=f(x)=\sum_{i=0}^{2g+1} \lambda_i x^i }$ and 
$\displaystyle{s^2=f(z)=\sum_{i=0}^{2g+1} \lambda_i z^i }$.
%%%%%%%%%%%%%%%%%%%%%%%%%%%%%%%%%%%%%%%%%%%%%%%%%%
%
Using 
\begin{align}
&[(x, z)]=\frac{y+s}{x-z} \frac{1}{2y},
\label{2e31}\\
&\frac{\partial}{\partial z} [(x, z)]
=\frac{2 y s+2f(z)+f'(z)(x-z)}{4(x-z)^2 y s}
\sim \frac{1}{(x-z)^2}, 
\label{2e32}
\end{align}
we start from the differential form of the third kind~\cite{Baker1,Baker2}
\begin{align}
&\int^x_a \dd x  \int^z_c \dd z \ \frac{\partial}{\partial z}[x, z]
=\int_a^x ( [(x, z)]-[(x,c)] )\dd x  .
\label{2e33}
\end{align}
Near the singularity, we obtain
\begin{align}
&\int^x_a \dd x  \int^z_c \dd z \ \frac{\partial}{\partial z}[x, z]
\sim \int^x_a \dd x  \int^z_c \dd z \frac{1}{(x-z)^2}=\log\left( \frac{(x-z) (a-c)}{(x-c)(a-z)}\right) .
\label{2e34}
\end{align}
 In order to obtain the standard differential form of the third kind,
we add the sum of the product of differential forms of 
the first and the second kinds in the form
\begin{align}
&\sum_{i=1}^g u_i ^{x, a} \eta_i^{z, c}, 
\label{2e35}\\
&\text{where} \nonumber\\ 
&u_i^{x, a}=\int^x_a \frac{x^{i-1} \dd x}{y}, 
\label{2e36}\\
&\eta_i^{z, c}=\int^z_c \frac{\dd z}{4s} \sum_{k=i}^{2g+1-i}
\lambda_{k+1+i} (k+1-i) z^k .
\label{2e37}
\end{align}
%
%%%%%%%%%%%%%%%%%%%%%%%%
%
Thus, the standard differential of the third kind is given in the form
\begin{align}
&\int_a^x ( [(x, z)]-[(x, c)] )\dd x+\sum_{i=1}^g u_i ^{x, a} \eta_i^{z, c}. 
\label{2e38}
\end{align}
By using some function $f_1(x,z)$, which is determined in the following, we can 
express this standard differential of the third kind is given in the form
\begin{align}
&\int_a^x ( [(x, z)]-[(x, c)] )\dd x+\sum_{i=1}^g u_i ^{x, a} \eta_i^{z, c}
=\int^x_a  \frac{\dd x}{y} \int^z_c \frac{\dd z}{s} \frac{f_1(x, z)+2 y s}{4(x-z)^2} .
\label{2e39}
\end{align}
Then we have
\begin{align}
&\frac{\partial}{\partial z} [(x, z)]+\frac{1}{4 ys } \sum_{i=1}^g x^{i-1}
\sum_{k=i}^{2g+1-i}
\lambda_{k+1+i} (k+1-i) z^k=\frac{f_1(x, z)+2 y s}{4(x-z)^2 y s } ,
\label{2e40}
\end{align}
which gives
\begin{align}
&\frac{2 s y+2f(z)+f'(z)(x-z)}{4(x-z)^2}+\frac{1}{4} \sum_{i=1}^g x^{i-1}
\sum_{k=i}^{2g+1-i}
\lambda_{k+1+i} (k+1-i) z^k=\frac{f_1(x, z)+2 y s}{4(x-z)^2} .
\label{2e41}
\end{align}
By using $\displaystyle{y^2=f(x)=\sum_{i=0}^{2g+1} \lambda_i x^i }$,
we obtain
\begin{align}
f_1(x, z)=\sum_{i=0}^g x^i z^i [2 \lambda_{2i} +\lambda_{2i+1}(x+z)] ,
\label{2e42}
\end{align}
with $f_1(x,x)=2f(x)$. 
By adding the special combination of the first and the second differentials of the form
$\displaystyle{\sum_{i=1}^g u_i ^{x, z} \eta_i^{t, c}}$, we obtain the special form of 
$f_1 (x, z)$. 
If $f_1 (x, z)$ satisfies following conditions, 
\begin{align}
&{\rm i)}\ f_1(x, z)=f_1(z, x) , 
\label{2e43}\\
&{\rm ii)}\ f_1(x, x)=2 f(x)  ,
\label{2e44}\\
&{\rm iii)}\ \frac{\partial }{\partial x} f_1(x, z) \biggr|_{x=z}=\frac{\partial }{\partial z} f(z) , 
\label{2e45}
\end{align}
we obtain the general $f_1(x, z)$.
Then we construct hyperelliptic $\vartheta$ functions with $g$ zero points by 
replacing $x \rightarrow x_j$,  
$a \rightarrow a_j$ and take the sum over $j$. Further, we replace $z \rightarrow x$, $c \rightarrow a$
in the form 
\begin{align}
&\sum^g_{j=1} \int_{a_j}^{x_j} ( [(x_j, x)]-[(x_j, a)]) \dd x_j 
+\sum_{i=1}^g u_i \eta_i^{x, a}
%\nonumber\\
=\sum^g_{j=1} \int^{x_j}_{a_j} \frac{\dd x_j}{y_j} \int^x_a \frac{\dd x}{y} \frac{(f_1(x_j, x)+2 y_j y)}{4(x_j-x)^2} 
\nonumber\\
&\sim \log \vartheta ({\bf u}^{x, a}-\sum_{k=1}^g {\bf u}^{x_k, a_k})
=\log \vartheta ({\bf u}^{x, a}-{\bf u}), 
\label{2e46}\\
%\end{align}
%
%
&\text{where}\nonumber\\
%
%\begin{align}
&({\bf u})_i=u_i=\sum^g_{j=1} u_i^{{x_j}, {a_j}}
=\sum^g_{j=1}\int^{x_j}_{a_j} \frac{x_j^{i-1} \dd x_j}{y_j} .
\label{2e47}
\end{align}
Riemann constants are absorbed into parameters $a, a_1, a_2, \cdots, a_g$,
that is, the constant shift of ${\bf u}^{x, a}, {\bf u}^{x_1, a_1}, \cdots, {\bf u}^{x_g, a_g}$.
In order to make ${\bf u}^{x, a}$ and ${\bf u}^{x_j, a_j}$ symmetric, we replace 
$y_i \rightarrow -y_i$, which gives
\begin{align}
&- \log \vartheta ({\bf u}^{x, a}+\sum_{k=1}^g {\bf u}^{x_k, a_k})=- \log \vartheta ({\bf u}^{x, a}+{\bf u})
=- \log \vartheta ({\bf \hat{u}})=-\log \vartheta(\hat{u}_1, \cdots, \hat{u}_g) 
\nonumber\\
&\sim \sum^g_{j=1} \int^{x_j}_{a_j}  \frac{\dd x_j}{y_j} \int^x_a \frac{\dd x}{y} \frac{(f_1(x_j, x)-2 y_j y)}{4(x_j-x)^2} , 
\label{2e48}
\end{align}
where $\hat{u}_i=u^{x, a}_i +u_i$.
In order to construct the $\sigma( \hat{u}_1, \hat{u}_2, \cdots, \hat{u}_g)$ function, which is 
the potential of the $\wp_{ij}(\hat{u}_1, \hat{u}_2, \cdots, \hat{u}_g)$ 
function, we multiply the 
prefactor $\exp(-\sum c_{ij} \hat{u}_i \hat{u}_j)$ to the $\wp_{ij}(\hat{u}_1, \hat{u}_2, \cdots, \hat{u}_g)$ 
function, which play the role of the 
constant shift of $\wp_{ij}(\hat{u}_1, \hat{u}_2, \cdots, \hat{u}_g)$ function.
Thus, we have 
\begin{align}
&- \log \sigma (\hat{u}_1,  \cdots, \hat{u}_g)
=\sum^g_{j=1} \int^{x_j}_{z_j}  \frac{\dd x_j}{y_j} \int^x_a \frac{\dd x}{y} \frac{(f_1(x_j, x)-2 y_j y)}{4(x_j-x)^2} .
\label{2e49}
\end{align}
By using the relation,
\begin{align}
&\frac{\partial}{\partial x_j}=\sum_{\ell=1}^{g}  \frac{\partial \hat{u}_{\ell}} {\partial x_j} \frac{\partial}{\partial u_{\ell}}
=\sum_{\ell=1}^{g}  \frac{\partial u^{x_j, a_j}_{\ell}} {\partial x_j} \frac{\partial}{\partial u_{\ell}}
=\sum_{\ell=1}^{g} \frac{x_j^{\ell-1}}{y_j}\frac{\partial}{\partial u_{\ell}}  ,
\label{2e50}\\
&\frac{\partial}{\partial x}=\sum_{m=1}^{g} \frac{\partial \hat{u}_m } {\partial x} \frac{\partial}{\partial u_m}
=\sum_{m=1}^{g} \frac{\partial u^{x, a}_m } {\partial x} \frac{\partial}{\partial u_m}
=\sum_{m=1}^{g} \frac{x^{m-1}}{y} \frac{\partial}{\partial u_m} ,
\label{2e51}
\end{align}
and applying  $\displaystyle{\frac{\partial^2}{\partial x \partial x_j}}$ to Eq.(\ref{2e49}), we obtain 
\begin{align}
&& \sum^g_{\ell,m=1} x_j^{\ell-1} x^{m-1} \wp_{\ell m}(\hat{u}_1, \cdots, \hat{u}_g)
=\frac{f_1(x_j, x)-2 y_j y}{4(x_j-x)^2}, \ (1 \le j \le g) ,  
\label{2e52}
\end{align}
where 
\begin{align}
\wp_{\ell m}(\hat{u}_1, \cdots, \hat{u}_g)=-\frac{\partial^2 }{ \partial \hat{u}_m \partial \hat{u}_{\ell}} 
\log \sigma (\hat{u}_1, \cdots, \hat{u}_g) .
\label{2e53}
\end{align}
As $\hat{u}^{x, a}$ and $\hat{u}^{x_k, a_k}$ are symmetric, by replacing $x,y \rightarrow x_k, y_k$,
the above equation is equivalent to the fundamental relation of the form
\begin{align}
&\sum^g_{\ell,m=1} x_j^{\ell-1} x_k^{m-1}  \wp_{\ell m}(\hat{u}_1, \cdots, \hat{u}_g)
=\frac{f_1(x_j, x_k)-2 y_j y_k}{4(x_j-x_k)^2} ,  \ (1 \le j, k \le g) .
\label{2e54}
\end{align}
Then we take the limit of $x \rightarrow \infty$, which gives $u_i^{x, a} \rightarrow {\rm (const.)}$, 
that is, $\hat{u}_i \rightarrow u_i +\text{(const).}$, 
and this constant is absorbed into the Riemann's constant. Hence, we obtain the fundamental relation 
of the form
\begin{align}
&\sum^g_{\ell,m=1} x_j^{\ell-1} x_k^{m-1}  \wp_{\ell m}(u_1, \cdots, u_g)
=\frac{f_1(x_j, x_k)-2 y_j y_k}{4(x_j-x_k)^2}  ,  \ (1 \le j, k \le g) , 
\label{2e55}\\
&\text{where}
\nonumber\\
&u_i=\sum_{j=1}^g \int^{x_j}_{a_j} \frac{x_j^{i-1} \dd x_j }{y_j} .
\nonumber
\end{align}
%

%%%%%%%%%%%%%%%%%%%%%%%%%%%%%%%%%%%%%%%%%%%%%%%%%%%%
\subsection{Derivation of the formula of $(-\zeta_j)$ function}
Next, we derive the formula of the $\zeta$ function. By replacing $c \rightarrow z$ in Eq.(\ref{2e46}), 
we start from the relation
\begin{align}
&\sum^g_{j=1}\int_{a_j}^{x_j} ( [(x_j,x)]-[(x_j, z)] )\dd x_j
+\sum_{i=1}^g u_i \eta_i^{x, z}
=\sum^g_{j=1} \int^{x_j}_{a_j}  \frac{\dd x_j}{y_j} \int^x_z \frac{\dd x}{y} \frac{(f_1(x_j, x)-2 y_j y)}{4(x_j-x)^2} 
\nonumber\\
&=\log \left( \frac{\sigma({\bf u}^{x, a}-\sum_{k=1}^g {\bf u}^{x_k, a_k})}
{\sigma ({\bf u}^{z, a}-\sum_{k=1}^g {\bf u}^{x_k, a_k})} \right) .
\label{2e56}
\end{align}
Then we have
\begin{align}
&\sum^g_{j=1} \frac{\partial x_j}{ \partial u_{\ell}}( [(x_j,x)]-[(x_j,z)])
+\eta_{\ell}^{x,z}
= \frac{\partial }{ \partial u_{\ell}}
\log \frac{\sigma ({\bf u}^{x, a}-\sum_{k=1}^g {\bf u}^{x_k, a_k})}
{\sigma ({\bf u}^{z, a}-\sum_{k=1}^g {\bf u}^{x_k, a_k})}
\nonumber\\
&=-\zeta_{\ell}({\bf u}^{x,a}-\sum_{k=1}^g {\bf u}^{x_k, a_k})
+\zeta_{\ell}({\bf u}^{z,a}-\sum_{k=1}^g {\bf u}^{x_k, a_k}) ,
\label{2e57}
\end{align}
which gives 
\begin{align}
&\sum^g_{j=1} \frac{ y_j \chi_{g-\ell}(x_j;x_1,x_2,\cdots,x_g)}{F'(x_j)} 
\Big(\frac{y_j+y}{x_j-x} \frac{1}{2y_j} -\frac{y_j+s}{x_j-z} \frac{1}{2y_j}\Big)
+\eta_{\ell}^{x,z}
\nonumber\\
&=-\zeta_{\ell}({\bf u}^{x,a}-\sum_{k=1}^g {\bf u}^{x_k, a_k})
+\zeta_{\ell}({\bf u}^{z,a}-\sum_{k=1}^g {\bf u}^{x_k, a_k})  .
\label{2e58}
\end{align}
Then, by replacing $y_j \rightarrow  -y_j$, which cause $u^{x_j, a_j} \rightarrow -u^{x_j, a_j} $, we obtain 
\begin{align}
&\zeta_{\ell}({\bf u}^{x,a}+\sum_{k=1}^g {\bf u}^{x_k, a_k})
+\eta_{\ell}^{x, a}+\sum_{k=1}^g \eta_{\ell}^{x_k, a_k}
-\frac{1}{2} \sum^g_{j=1} \frac{ \chi_{g-\ell}(x_j;x_1,x_2,\cdots,x_g)}{F'(x_j)} 
\Big(\frac{y-y_j}{x-x_i}\Big)
\nonumber\\
&=\zeta_{\ell}({\bf u}^{z, a}+\sum_{k=1}^g  {\bf u}^{x_k, a_k})
+\eta_{\ell}^{z, a}+\sum_{k=1}^g \eta_{\ell}^{x_k, a_k}
-\frac{1}{2} \sum^g_{j=1} \frac{ \chi_{g-\ell}(x_j;x_1,x_2,\cdots,x_g)}{F'(x_j)} 
\Big(\frac{s-y_j}{z-x_i}\Big) ,
\label{2e59}
\end{align}
where we add  $\sum_{k=1}^g \eta_{\ell}^{x_k, a_k}$ in both sides of Eq.(\ref{2e59}) in order to make $\eta_{\ell}^{z, a}$ and $\eta_{\ell}^{x_k, a_k}$ are symmetric.
Therefore,  the left-hand side of Eq.(\ref{2e59}) becomes independent of $x$. Next, we introduce
\begin{align}
R(t)=(t-x)(t-x_1)\cdots (t-x_g)=(t-x)F(t), 
\nonumber
\end{align}
and define
\begin{align}
g_{g-i}(x, x_1, \cdots, x_g)=&\frac{y \chi_{g-i}(x; x_1, \cdots, x_g)}{R'(x)}+\frac{y_1 \chi_{g-i}(x_1; x, x_2, \cdots, x_g)}{R'(x_1)}
+\cdots
\nonumber\\
&+\frac{y_g \chi_{g-i}(x_g; x_1, \cdots, x_{g-1}, x)}{R'(x_g)} , 
\nonumber\\
g_{g-i-1}(x_1, \cdots, x_g)=&\frac{y_1 \chi_{g-i-1}(x_1; x_2, \cdots, x_g)}{F'(x_1)}+\frac{y_2 \chi_{g-i-1}(x_2; x_1, \cdots, x_g)}{F'(x_2)}
+\cdots
\nonumber\\
&+\frac{y_g \chi_{g-i-1}(x_g; x_1, \cdots, x_{g-1})}{R'(x_g)} .
\nonumber
\end{align}
Then we can prove that
\begin{align}
\sum^g_{j=1} \frac{ \chi_{g-\ell}(x_j;x_1,x_2,\cdots,x_g)}{F'(x_j)} 
\Big(\frac{y-y_j}{x-x_i}\Big)=g_{g-i}(x, x_1, \cdots, x_g)-g_{g-i-1}(x_1, \cdots, x_g) .
\label{2e60}
\end{align}
First, we consider $I(x)=\chi_{g-i}(x_1; x, x_2, \cdots, x_g)-\chi_{g-i}(x_1; x_1, x_2, \cdots, x_g)$, which is the linear 
function of $x$. By using $I(x_1)=0$ and 
$\displaystyle{\frac{\partial \chi_{g-i} (x_1; x, x_2, \cdots, x_g) }{\partial x} 
=-\chi_{g-i-1}(x_1; x_2, \cdots, x_g)}$ of Eq.(\ref{2e10}), we have
\begin{align}
\chi_{g-i}(x_1; x, x_2, \cdots, x_g)-\chi_{g-i}(x_1; x_1, x_2, \cdots, x_g)=(x_1-x)\chi_{g-i-1}(x_1; x_2, \cdots, x_g) .
\nonumber
\end{align}
Then we obtain
\begin{align}
&\frac{\chi_{g-i}(x_1; x, x_2, \cdots, x_g)}{R'(x_1)}
=\frac{\chi_{g-i}(x_1; x_1, x_2, \cdots, x_g)+(x_1-x)\chi_{g-i-1}(x_1;  x_2, \cdots, x_g)}{(x_1-x) F'(x_1)}
\nonumber\\
&=-\frac{\chi_{g-i}(x_1; x_1, x_2, \cdots, x_g)}{F'(x_1)} \frac{1}{x-x_1}
+\frac{\chi_{g-i-1}(x_1;  x_2, \cdots, x_g)}{F'(x_1)} .
\nonumber
\end{align}
By using the similar result for $\displaystyle{\frac{\chi_{g-i}(x_k; x, x_2, \cdots, x_g)}{R'(x_k)}}$, we obtain
\begin{align}
&g_{g-i}(x, x_1, \cdots, x_g)-g_{g-i-1}(x_1, \cdots, x_g)=\frac{y \chi_{g-i}(x; x_1, \cdots, x_g)}{R'(x)}
-\sum_{k=1}^g \frac{\chi_{g-i}(x_k; x_1, \cdots, x_g)}{F(x_k)}\frac{y_k}{x-x_k}
\nonumber\\
&=\sum_{k=1}^g \frac{\chi_{g-i}(x_k; x_1, \cdots, x_g)}{F'(x_k)}\frac{y-y_k}{x-x_k}, 
\nonumber
\end{align}
where we used Eq.(\ref{2e5}) and $R'(x)=F(x)$. This completes the 
proof of Eq.(\ref{2e60}). 
Substituting this relation into Eq.(\ref{2e59}), we obtain
\begin{align}
&\zeta_{\ell}({\bf u}^{x,a}+\sum_{k=1}^g {\bf u}^{x_k, a_k})+\eta_{\ell}^{x, a}
+\sum_{k=1}^g \eta_{\ell}^{x_k, a_k}
-\frac{1}{2} g_{g-i}(x, x_1, \cdots, x_g) 
=\text{($x$ independent)}, 
\label{2e61}
\end{align}
which gives
\begin{align}
&-\zeta_{\ell}({\bf u}^{x,a}+\sum_{k=1}^g {\bf u}^{x_k, a_k})
=\left(\eta_{\ell}^{x, a} -\frac{y \chi_{g-\ell}(x; x_1, \cdots, x_g)}{2R'(x)}\right)  
\nonumber\\
&+\sum_{k=1}^g \eta_{\ell}^{x_k, a_k}
-\frac{1}{2} \sum_{k=1}^g \frac{y_k \chi_{g-\ell}(x_k; x_1, \cdots, x_{k-1}, x, x_{k+1}, \cdots, x_g)}{2(x_k-x) F'(x_k)}+C.
\label{2e62}
\end{align}
Then we take the limit of $x \rightarrow \infty$. Thus, we have $u^{x,a} \rightarrow \text{(const.)}$ and this constant is absorbed into the Riemann's 
constant. 
The first term of the right-hand side of Eq.(\ref{2e62}) gives 
$\displaystyle{\left(\eta_{\ell}^{x, a} -\frac{y \chi_{g-\ell}(x; x_1, \cdots, x_g)}{2R'(x)}\right) \rightarrow \text{(const.)}}$. 
Thus, by using Eq.(\ref{2e37}), we obtain
 Eq.(\ref{2e23}) in the form
\begin{align}
&-\zeta_{\ell}(\sum_{k=1}^g {\bf u}^{x_k, a_k})
=\sum_{k=1}^g \eta_{\ell}^{x_k, a_k}
-\frac{1}{2} \sum_{k=1}^g\frac{y_k \chi_{g-\ell-1}(x_k; x_1, \cdots, \widecheck{x}_k, \cdots,x_g)}{F'(x_i)}+C. 
\label{2e63}
\end{align}
%
%\newpage
%%%%%%%%%%%%%%%%%%%%%%%%%%%%%%%%%%%%%%%%%%%%%%%%%%%%%%%%%%%
%%%%%%%%%%%%%%%%%%%%%%%%%%%%%%%%%%%%%%%%%%%%%%%%%%%%%%
\subsection{Bolza's formula (I)}
The fundamental relation, which generates all differential equations, is 
given by the Baker~\cite{Baker1}.
We replace $x \rightarrow x_0$, $a \rightarrow a_0$ and we denote 
$\widehat{\bf u}
={\bf u}^{x,a}+\sum_{k=1}^g {\bf u}^{x_k,a_k} \rightarrow \widehat{\bf u}
=\sum_{r=0}^g {\bf u}^{x_r, a_r}$. 
Then we start from Eq.(\ref{2e52}) and Eq.(\ref{2e54}) in the form
\begin{align}
&\sum_{i, j=1}^g \wp_{i j} (\widehat{\bf u}) x_0^{i-1} x_r^{j-1}
=\frac{f_1(x_0,x_r)-2 y_0 y_s}{4(x_0-x_r)^2}, \quad 
(1\le r \le g),
\label{2e64}\\
&\sum_{i, j=1}^g \wp_{i j} (\widehat{\bf u}) x_r^{i-1} x_s^{j-1}
=\frac{f_1(x_r,x_s)-2 y_r y_s}{4(x_r-x_s)^2}, \quad 
(1\le r < s \le g) .
\label{2e65}
\end{align}
Then, for arbitrary variables $e_1$, $e_2$, we construct the following quantity
\begin{align}
\sum_{i, j=1}^g \wp_{i j} (\widehat{\bf u}) e_1^{i-1} e_2^{j-1} .
\label{2e66}
\end{align}
The right-hand side of Eq.(\ref{2e66}) must satisfies\\
\noindent
i)\ it is the $(g-1)$-th order polynomial function with respect to $e_1$, $e_2$ respectively,
\\
ii)\ in the limit of $e_1 \rightarrow x_p$, $e_2 \rightarrow x_q$, 
it must become $\displaystyle{\frac{f_1(x_p,x_q)-2y(x_p) y(x_q)}{4(x_p-x_q)^2}}$. \\
Then the right-hand side of Eq.(\ref{2e66}) takes the form 
\begin{align}
&\text{(R.H.S.)}= \sum_{\substack{r, s=0 \\r<s}}^g \frac{R(e_1; \widecheck{x_r}, \widecheck{x_s}) 
R(e_2; \widecheck{x_r}, \widecheck{x_s})}{R(x_r; \widecheck{x_r}, \widecheck{x_s}) 
R(x_s; \widecheck{x_r}, \widecheck{x_s})}
\frac{f_1(x_r,x_s)-2y(x_r) y(x_s)}{4(x_r-x_s)^2}, 
\label{2e67}\\
&\text{where}
\nonumber\\
&R(x)=(x-x_0) (x-x_1)(x-x_2) \cdots (x-x_{g}) .
\label{2e68}
\end{align}
In the limit of $e_1 \rightarrow x_p$, $e_2 \rightarrow x_q$,
only $x_r=x_p$, $x_s=x_q$ or $x_r=x_q$, $x_s=x_p$ gives non-zero contribution.
Thus, we obtain
\begin{align}
&\text{(R.H.S.)}=\sum_{\substack{r, s=0 \\r<s}}^g \frac{R(e_1)R(e_2)}
{(e_1-x_r)(e_1-x_s)(e_2-x_r)(e_2-x_s)} 
\nonumber\\
&\times \frac{(x_r-x_s)}{R(x_r; \widecheck{x_r})}\frac{(x_s-x_r)}{R(x_s; \widecheck{x_s})}
\frac{f_1(x_r,x_s)-2y(x_r) y(x_s)}{4(x_r-x_s)^2}
\nonumber\\
&=R(e_1)R(e_2) \sum_{\substack{r, s=0 \\r<s}}^g \frac{2y(x_r) y(x_s)-f_1(x_r,x_s)}
{4(e_1-x_r)(e_1-x_s)(e_2-x_r)(e_2-x_s)R'(x_r) R'(x_s)} .
\label{2e69}
\end{align}
Hence, we obtain
\begin{align}
\hskip -10mm \sum_{i, j=1}^g \wp_{i j} (\widehat{\bf u}) e_1^{i-1} e_2^{j-1}
=R(e_1)R(e_2) \sum_{\substack{r, s=0 \\r<s}}^g \frac{2y(x_r) y(x_s)-f_1(x_r,x_s)}
{4(e_1-x_r)(e_1-x_s)(e_2-x_r)(e_2-x_s) R'(x_r) R'(x_s)} .
\label{2e70}
\end{align}
Then we take the limit $x_0 \rightarrow \infty$, and we obtain  
$\widehat{\bf u} \rightarrow {\bf u}+\text{(const.)}$ and this constant 
is absorbed into the Riemann's constant. Thus, we obtain the Bolza's  
formula (I)~\cite{Baker1,Bolza1,Bolza2}
\begin{align}
&\sum_{i, j=1}^g \wp_{i j} ({\bf u}) e_1^{i-1} e_2^{j-1}
=\lim_{x_0 \rightarrow \infty}R(e_1)R(e_2) \sum_{\substack{r, s=0 \\r<s}}^g \frac{2y(x_r) y(x_s)-f_1(x_r,x_s)}
{4(e_1-x_r)(e_1-x_s)(e_2-x_r)(e_2-x_s) R'(x_r) R'(x_s)} 
\nonumber\\
&=\lim_{x_0 \rightarrow \infty}R(e_1)R(e_2) \sum_{\substack{r, s=0 \\r<s}}^g \frac{2y(x_r) y(x_s)-f_1(x_r,x_s)}
{4G'(x_r) G'(x_s)} , 
\label{2e71}\\
\text{where}
\nonumber\\
&G(x)=(x-e_1)(x-e_2) (x-x_0)(x-x_1)(x-x_2) \cdots (x-x_{g})  .
\label{2e72}
\end{align}
In the following, we explicitly give explicit forms of $\wp_{ij}$ functions for genus two and three cases 
by using this Bolza's formula (I). \\
%%%%%%%%%%%%%%%%%%%%%%%%%%%%%%%%%%%%%%%%%%%%%
%%%%%%%%%%%%%%%%%%%%%%%%%%%%%%%%%%%%%%%%%%%%%%%%%%%%%%%%%
\vskip 3mm
\noindent
{\bf \underline{Genus two case}:}\\
In the genus two case, we obtain
\begin{align}
&\hskip -10 mm f_1(x_r, x_s)=\lambda_5 x_r^2 x_s^2 (x_r+x_s)+2 \lambda_4 x_r^2 x_s^2
+\lambda_3 x_r x_s (x_r+x_s)+2 \lambda_2 x_r x_s
+\lambda_1 (x_r+x_s)+2 \lambda_0,
\label{2e73}\\
&\hskip -10 mm G(x)=(x-e_1)(x-e_2) (x-x_0)(x-x_1)(x-x_2) ,
\label{2e74}\\
&\hskip -10 mm R(x)=(x-x_0) (x-x_1)(x-x_2) .
\label{2e75}
\end{align}
In this case, Bolza's formula (I) gives
\begin{align}
&\wp_{11} (u) +(e_1+e_2) \wp_{21}+e_1 e_2 \wp_{22}
=-\frac{(x_0-e_1)(x_0-e_2)}{(x_0-x_1)(x_0-x_2)} \frac{[2y_1 y_2-f_1(x_1,x_2)]}{4(x_1-x_2)^2}
\nonumber\\
&-\frac{(x_1-e_1)(x_1-e_2)}{(x_1-x_0)(x_1-x_2)} \frac{[2y_0 y_2-f_1(x_0,x_2)]}{4(x_0-x_2)^2}
-\frac{(x_2-e_1)(x_2-e_2)}{(x_2-x_0)(x_2-x_1)} \frac{[2y_0 y_1-f_1(x_0,x_1)]}{4(x_0-x_1)^2} .
\label{2e76}
\end{align}
This gives
\begin{align}
&\wp_{22} (u) =\frac{ [f_1(x_1,x_2)-2y_1 y_2]}{4(x_0-x_1)(x_0-x_2)(x_1-x_2)^2}
+\frac{[f_1(x_0,x_2)-2y_0 y_2]}{4(x_1-x_0)(x_1-x_2)(x_0-x_2)^2}
\nonumber\\
&+\frac{[f_1(x_0,x_1)-2y_0 y_1]}{4(x_2-x_0)(x_2-x_1)(x_0-x_1)^2} ,
\label{2e77}\\
&\wp_{21} (u) =-\frac{x_0 [f_1(x_1,x_2)-2y_1 y_2]}{4(x_0-x_1)(x_0-x_2)(x_1-x_2)^2}
-\frac{x_1 [f_1(x_0,x_2)-2y_0 y_2]}{4(x_1-x_0)(x_1-x_2)(x_0-x_2)^2}
\nonumber\\
&-\frac{x_2 [f_1(x_0,x_1)-2y_0 y_1]}{4(x_2-x_0)(x_2-x_1)(x_0-x_1)^2} ,
\label{2e78}\\
&\wp_{11} (u) =\frac{x_0^2 [f_1(x_1,x_2)-2y_1 y_2]}{4(x_0-x_1)(x_0-x_2)(x_1-x_2)^2}
+\frac{x_1^2[f_1(x_0,x_2)-2y_0 y_2]}{4(x_1-x_0)(x_1-x_2)(x_0-x_2)^2}
\nonumber\\
&+\frac{x_2^2 [f_1(x_0,x_1)-2y_0 y_1]}{4(x_2-x_0)(x_2-x_1)(x_0-x_1)^2}  .
\label{2e79}
\end{align}
In this way, by using $x_0$, we can write $\wp_{22}(u), \wp_{21}(u), \wp_{11}(u)$ in a symmetric way.
We must take the limit of $x_0 \rightarrow \infty$ in the end. Then we estimate 
$y_0 \sim  \sqrt{\lambda_5} x_0^{5/2}$, 
$f_1(x_0, x_1)\sim \lambda_5 (x_0+x_1) x_0^2 x_1^2\sim \lambda_5 x_0^3 x_1^2$ at large $x_0$, so that
$\displaystyle{ \frac{2y_0 y_2-f_1(x_0,x_2)}{4(x_0-x_2)^2} \sim
 \frac{-\lambda_5 x_0^3 x_1^2}{4x_0^2}\sim -\frac{\lambda_5}{4} x_0 x_1^2}$, which gives 
\begin{align}
&\wp_{11} (u) +(e_1+e_2) \wp_{21}+e_1 e_2 \wp_{22}
\nonumber\\
&=-\frac{2y_1 y_2-f_1(x_1,x_2)}{4(x_1-x_2)^2}
-\frac{\lambda_5}{4}\frac{(x_1-e_1)(x_1-e_2)x_2^2}{(x_1-x_2)}
-\frac{\lambda_5}{4}\frac{(x_2-e_1)(x_2-e_2)x_1^2}{(x_2-x_1)} 
\nonumber\\
&=\frac{f_1(x_1,x_2)-2y_1 y_2}{4(x_1-x_2)^2}-\frac{\lambda_5}{4} (e_1+e_2) x_1 x_2 
+\frac{\lambda_5}{4} e_1 e_2 (x_1+x_2) .
\label{2e80}
\end{align}
Hence, we obtain
%%%%%%%%%%%%%%%%%%%%%%%%%%%%%%%%%%%%%%%%
 ~\footnote{For the genus two hyperelliptic curve of the 
form $f(x)=x\prod_{i=1}^4(1-\mu_i x)$, Rosenhain\cite{Rosenhain1,Rosenhain2}
constructed 15 hyperelliptic functions, 
$x_1 x_2$, $(1-\mu_i x_1)(1-\mu_i x_2), (i=1,\cdots, 4)$, 
$\displaystyle{\frac{f_{\rho \sigma}(x_1,x_2) -2 y(x_1)y(x_2)}{4(x_1-x_2)^2}, 
(\rho < \sigma =0,1, \cdots, 4)}$
with $\displaystyle{f_{\rho \sigma}(x_1,x_2)
=\frac{g_{\rho \sigma}(x_2) f(x_1)}{g_{\rho \sigma}(x_1)}
+\frac{g_{\rho \sigma}(x_1) f(x_2)}{g_{\rho \sigma}(x_2)}}$,  
$g_{0 i}(x)=x(1-\mu_{i} x)$, 
$g_{i j}(x)=(1-\mu_{i} x)(1-\mu_{j} x)$. This $f_{\rho \sigma}(x_1,x_2)$
satisfies $f_{\rho \sigma}(x_1,x_1)=2f(x_1)$, 
$\displaystyle{\frac{\partial}{\partial x_1}f_{\rho \sigma}(x_1,x_2)|_{x_2=x_1}=f'(x_2)}$
of Eq(\ref{2e44}) and Eq.(\ref{2e45}).
%%%%%%%%%%%%%%%%%%%%%%%%%%%%%%%%%%%%%%%%%%%
}
\begin{align}
&\wp_{22}(u)=\frac{\lambda_5}{4} (x_1+x_2), 
\quad \wp_{21}(u)=-\frac{\lambda_5}{4} x_1 x_2, \quad
\wp_{11} (u)=\frac{f_1(x_1,x_2)-2y_1 y_2}{4(x_1-x_2)^2} .
\label{2e81}
\end{align}
%
%%%%%%%%%%%%%%%%%%%%%%%%%%%%%%%%%%%%
\vskip 3mm
\noindent
{\bf \underline{Genus three case}:}\\
In genus three case, we obtain
\begin{align}
&f_1(x_r, x_s)=\lambda_7 x_r^3 x_s^3 (x_r+x_s)+2 \lambda_6 x_r^3 x_s^3
+\lambda_5 x_r^2 x_s^2 (x_r+x_s)+2 \lambda_4 x_r^2 x_s^2
+\lambda_3 x_r x_s (x_r+x_s)
\nonumber\\
&+2 \lambda_2 x_r x_s
+\lambda_1 (x_r+x_s)+2 \lambda_0 ,
\label{2e82}\\
&G(x)=(x-e_1)(x-e_2) (x-x_0)(x-x_1)(x-x_2)(x-x_{3}) ,
\label{2e83}\\
&R(x)=(x-x_0) (x-x_1)(x-x_2)(x-x_{3}) .
\label{2e84}
\end{align}
In this case, after taking the limit $x_0 \rightarrow \infty$, we obtain 
\begin{align}
&\hskip -10 mm \wp_{11} (u) +(e_1+e_2) \wp_{21}+e_1 e_2\wp_{22}+(e_1^2+e_2^2)\wp_{31}
+e_1 e_2 (e_1+e_2) \wp_{32}+e_1^2 e_2^2 \wp_{33}
\nonumber\\
&\hskip -10 mm =-\frac{(x_3-e_1)(x_3-e_2)}{(x_3-x_1)(x_3-x_2)} \frac{[2y_1 y_2-f_1(x_1,x_2)]}{4(x_1-x_2)^2}
-\frac{(x_1-e_1)(x_1-e_2)}{(x_1-x_2)(x_1-x_3)} \frac{[2y_2 y_3-f_1(x_2, x_3)]}{4(x_2-x_3)^2}
\nonumber\\
&\hskip -10 mm -\frac{(x_2-e_1)(x_2-e_2)}{(x_2-x_1)(x_2-x_3)} \frac{[2y_3 y_1-f_1(x_3, x_1)]}{4(x_3-x_1)^2}
+\frac{\lambda_7}{4}\frac{(x_2-e_1)(x_2-e_2)(x_3-e_1)(x_3-e_2)}{(x_1-x_2)(x_1-x_3)}x_1^3
\nonumber\\
&\hskip -10 mm +\frac{\lambda_7}{4}\frac{(x_3-e_1)(x_3-e_2)(x_1-e_1)(x_1-e_2)}{(x_2-x_3)(x_2-x_1)}x_2^3
+\frac{\lambda_7}{4}\frac{(x_1-e_1)(x_1-e_2)(x_2-e_1)(x_2-e_2)}{(x_3-x_1)(x_3-x_2)}x_3^3 .
\label{2e85}
\end{align}
Comparing coefficients of the power of $e_1$ and $e_2$, we obtain 
\begin{align}
&\wp_{33}(u)=\frac{\lambda_7}{4} (x_1+x_2+x_3), 
\label{2e86}\\
&\wp_{32}(u)=-\frac{\lambda_7}{4} (x_1 x_2+x_2 x_3+x_3 x_1), 
\label{2e87}\\
&\wp_{31} (u)=\frac{\lambda_7}{4} x_1 x_2 x_3, 
\label{2e88}\\
&\wp_{22}=\frac{\lambda_7}{4} x_1 x_2 x_3 +\frac{f_1(x_1,x_2)-2 y_1 y_2}{4(x_3-x_1)(x_3-x_2)(x_1-x_2)^2} 
+\frac{f_1(x_2,x_3)-2 y_2 y_3}{4(x_1-x_2)(x_1-x_3)(x_2-x_3)^2} 
\nonumber\\
&+\frac{f_1(x_3,x_1)-2 y_3 y_1}{4(x_2-x_3)(x_2-x_1)(x_3-x_1)^2} , 
\label{2e89}\\
&\wp_{21}=-\frac{x_3 (f_1(x_1,x_2)-2 y_1 y_2)}{4(x_3-x_1)(x_3-x_2)(x_1-x_2)^2} 
-\frac{x_1(f_1(x_2,x_3)-2 y_2 y_3)}{4(x_1-x_2)(x_1-x_3)(x_2-x_3)^2}
\nonumber\\
&-\frac{x_2(f_1 (x_3,x_1)-2 y_3 y_1)}{4(x_2-x_3)(x_2-x_1)(x_3-x_1)^2}  ,
\label{2e90}\\
&\wp_{11}=\frac{x_3^2 (f_1(x_1,x_2)-2 y_1 y_2)}{4(x_3-x_1)(x_3-x_2)(x_1-x_2)^2} 
+\frac{x_1^2(f_1(x_2,x_3)-2 y_2 y_3)}{4(x_1-x_2)(x_1-x_3)(x_2-x_3)^2}
\nonumber\\
&+\frac{x_2^2(f_1(x_3,x_1)-2 y_3 y_1)}{4(x_2-x_3)(x_2-x_1)(x_3-x_1)^2} .
\label{2e91}
\end{align}
We have checked up to genus four case that 
$\displaystyle{\wp_{ij}(u)=-\frac{\partial \zeta_j}{\partial u_i} }$ by using $\dd (-\zeta_j)$ of Eq.(\ref{2e23})
and $\wp_{ij}(u)$ by using Bolza's formula (I) of Eq.(\ref{2e71}) give the same $\wp_{ij}(u)$ expression.

%%%%%%%%%%%%%%%%%%%%%%%%%%%%%%%%%%%%%%%%%%%%%%%%%%%%%%%%%%%%%
%%%%%%%%%%%%%%%%%%%%%%%%%%%%%%%%%%%%%%%%%%%%%%%%%%%%%%%%%%%%%
\subsection{Bolza's formula (II)}
We rewrite the Bolza's formula (I) into more useful form, which is the starting point 
of Baker's  paper~\cite{Baker4}.
By using Eq.(\ref{2e5}), we obtain\\
\begin{align}
\frac{x^k}{R(x)}=\sum_{i=0}^g \frac{x_i^k}{(x-x_i) R'(x_i)}, \quad (k=0,1, \cdots, g).
\label{2e92} 
\end{align}
Taking $k=g$ and multiply $x$ in both side of Eq.(\ref{2e92}), we have
\begin{align}
\frac{x^{g+1}}{R(x)}=\sum_{i=0}^g \frac{x_i^g}{R'(x_i)}\frac{(x-x_i+x_i) }{(x-x_i)}
=\sum_{i=0}^g \frac{x_i^g}{R'(x_i)}+\sum_{i=0}^g \frac{x_i^{g+1}}{(x-x_i) R'(x_i)}
=1+\sum_{i=0}^g \frac{x_i^{g+1}}{(x-x_i) R'(x_i)} .
\nonumber 
\end{align}
Then, for the $(g+1)$ polynomial function of the form 
$f(x)=a_{g+1} x^{g+1}+a_{g} x^{g}+\cdots+a_1 x+a_0$, we obtain
\begin{align}
\frac{f(x)}{R(x)}=a_{g+1}+\sum_{i=0}^g \frac{f(x_i)}{(x-x_i) R'(x_i)} .
\label{2e93} 
\end{align}
Later, we will use this formula in Eq.(\ref{2e98}).  
Then we start from Eq.(\ref{2e71}) without taking the 
$x_0 \rightarrow \infty$ limit in the form 
\begin{align}
&I=R(e_1) R(e_2) \sum_{\substack{r,s=0 \\ r<s} }^{g}\frac{[2 y_r y_s-f_1(x_r, x_s)]}{4G'(x_r) G'(x_s)} \quad 
=I_1+I_2,
\label{2e94}\\
&\text{where}
\nonumber\\
&I_1=R(e_1) R(e_2) \sum_{\substack{r,s=0 \\ r<s}}^{g}\frac{2 y_r y_s}{4G'(x_r) G'(x_s)}, 
\label{2e95}\\
&I_2=-R(e_1) R(e_2) \sum_{\substack{r,s=0 \\ r<s}}^{g}\frac{f_1(x_r, x_s)}{4G'(x_r) G'(x_s)} .
\label{2e96}
\end{align}
First, we rewrite $I_1$ in the form
\begin{align}
&I_1=R(e_1) R(e_2) \sum_{\substack{r,s=0 \\ r<s}}^{g}\frac{2 y_r y_s}{4G'(x_r) G'(x_s)}
\nonumber\\
&=R(e_1) R(e_2) \sum_{\substack{r,s=0 \\ r<s}}^{g}\frac{2 y_r y_s}{4(x_r-e_1)(x_r-e_2)(x_s-e_1)(x_s-e_2)R'(x_r) R'(x_s)}
\nonumber\\
&=\frac{R(e_1)R(e_2)}{4} \left[ 
\left( \sum_{r=0}^g \frac{y_r}{(x_r-e_1)(x_r-e_2)R'(x_r)} \right) 
\left( \sum_{s=0}^g \frac{y_s}{(x_s-e_1)(x_s-e_2)R'(x_s)} \right) \right. 
\nonumber\\
&\left. -\sum_{r=0}^g \left(\frac{y_r}{(x_r-e_1)(x_r-e_2)R'(x_r)} \right)^2 
\right] 
\nonumber\\
&=\frac{R(e_1)R(e_2)}{4}  
\left( \sum_{r=0}^g \frac{y_r}{(x_r-e_1)(x_r-e_2)R'(x_r)} \right)^2
-\frac{R(e_1)R(e_2)}{4} \sum_{r=0}^g \frac{f(x_r)}{(G'(x_r))^2}  
\nonumber\\
&=I_{11}+I_{12}.
\label{2e97}
\end{align}
Next, we rewrite $I_2$ in the form
\begin{align}
&I_2=-\frac{R(e_1) R(e_2)}{4 (e_1-e_2)^2}  \sum_{\substack{r,s=0 \\ r<s}}^{g}\frac{f_1(x_r, x_s)}{R'(x_r) R'(x_s)}
\left( \frac{1}{e_1-x_r}-\frac{1}{e_2-x_r} \right) \left( \frac{1}{e_1-x_s}-\frac{1}{e_2-x_s} \right)
\nonumber\\
&=-\frac{R(e_1) R(e_2)}{8 (e_1-e_2)^2} 
\left[  \sum_{s=0}^{g}\frac{1}{R'(x_s)} \left( \frac{1}{e_1-x_s}-\frac{1}{e_2-x_s} \right)
\sum_{r=0}^{g}\frac{1}{R'(x_r)} \left( \frac{1}{e_1-x_r}-\frac{1}{e_2-x_r} \right) f_1(x_r, x_s) \right.
\nonumber\\
&\left. -\sum_{r=0}^{g}\frac{2f(x_r)}{(R'(x_r))^2} \left( \frac{1}{e_1-x_r}-\frac{1}{e_2-x_r} \right)^2
\right]
\nonumber\\
&=-\frac{R(e_1) R(e_2)}{8 (e_1-e_2)^2} 
\left[  \sum_{s=0}^{g}\frac{1}{R'(x_s)} \left( \frac{1}{e_1-x_s}-\frac{1}{e_2-x_s} \right)
\left( \frac{f_1(e_1,x_s)}{R(e_1)}-\frac{f_1(e_2, x_s)}{R(e_2)} \right) \right.
\nonumber\\
&\left. -\sum_{r=0}^{g}\frac{2f(x_r)}{(R'(x_r))^2} \frac{(e_1-e_2)^2}{(e_1-x_r)^2(e_2-x_r)^2}
\right]
\nonumber\\
&=-\frac{R(e_1) R(e_2)}{8 (e_1-e_2)^2} 
\left[  \frac{1}{R(e_1)} \left( \frac{f_1(e_1,e_1)}{R(e_1}-\frac{f_1(e_2,e_1)}{R(e_2}\right)
-\frac{1}{R(e_2)} \left( \frac{f_1(e_1,e_2)}{R(e_1}-\frac{f_1(e_2,e_2)}{R(e_2}\right)
\right]
\nonumber\\
&+\frac{R(e_1) R(e_2)}{4} \sum_{r=0}^{g}\frac{f(x_r)}{(G'(x_r))^2} 
\nonumber\\
&=-\frac{R(e_1) R(e_2)}{8 (e_1-e_2)^2} 
\left( \frac{2 f(e_1)}{R(e_1)^2}+\frac{2 f(e_2)}{R(e_2)^2}
-\frac{2 f_1(e_1.e_2)}{R(e_1) R(e_2)}\right)
+\frac{R(e_1) R(e_2)}{4} \sum_{r=0}^{g}\frac{f(x_r)}{(G'(x_r))^2} 
\nonumber\\
&=I_{21}+I_{22} , 
\label{2e98}
\end{align}
where we used Eq.(\ref{2e44}), Eq.(\ref{2e93}) and the fact that $f_1(x_r , x_s)$ is the $(g+1)$-th order 
polynomial with respect to $x_r$ and $x_s$ respectively.
As $I_{12}$ and $I_{22}$ cancel out, by taking the $x_0 \rightarrow \infty$ limit, 
we finally obtain the Bolza's 
formula (II)~\cite{Baker1,Bolza1,Bolza2} of the form
\begin{align}
&\sum_{i, j=1}^g \wp_{i j} (u) e_1^{i-1} e_2^{j-1}
=\lim_{x_0 \rightarrow \infty}
\left[ \frac{1}{4} R(e_1) R(e_2)\left(\sum_{r=0}^{g}\frac{ y_r }{(e_1-x_r)(e_2-x_r)R'(x_r)} \right)^2 \right.
\nonumber\\
&\left. -\frac{f(e_1)R(e_2)}{4 (e_1-e_2)^2 R(e_1)}-\frac{f(e_2)R(e_1)}{4 (e_1-e_2)^2 R(e_2)}
+\frac{f_1(e_1,e_2)}{4 (e_1-e_2)^2} \right].
\label{2e99}
\end{align}
%

%%%%%%%%%%%%%%%%%%%%%%%%%%%%%%%%%%%%%%%%%%%%%%%%%%%%%%%%
%%%%%%%%%%%%%%%%%%%%%%%%%%%%%%%%%%%%%%%%%%%%%%%%%%%%%%%%%%%
\section{The formula of differential equations for hyperelliptic $\wp_{ij}$ functions} 
\setcounter{equation}{0}

%\subsection{Rewriting Bolza's formula (II) into Baker's starting relation}
We start from the Bolza's formula (II) Eq.(\ref{2e99}) in the form
\begin{align}
&4 (e_1-e_2)^2 \sum_{i, j=1}^g \wp_{i j} (u) e_1^{i-1} e_2^{j-1} -f_1(e_1, e_2)
=\lim_{x_0 \rightarrow \infty} R(e_1) R(e_2)\Omega_{12}, 
\label{3e1}\\
&\text{where}
\nonumber\\
&\Omega_{12}=(e_1-e_2)^2 \left(\sum_{r=0}^{g}\frac{ y_r }{(e_1-x_r)(e_2-x_r)R'(x_r)} \right)^2 
 -\frac{f(e_1)}{R(e_1)^2}-\frac{f(e_2)}{R(e_2)^2}
\label{3e2}
\end{align}
We introduce $\Delta_i\ (i=1, 2)$,  $\Delta_{12}$ and $\varphi_i\ (i=1, 2)$ in the form
\begin{align}
&\Delta_i=\sum_{r=0}^g \frac{y_r}{(e_i-x_r) R'(x_r)}, 
\label{3e3}\\
&\Delta_{12}=-\frac{\Delta_1-\Delta_2}{e_1-e_2}=\sum_{r=0}^g \frac{y_r}{(e_1-x_r)(e_2-x_r)R'(x_r)}
=\sum_{r=0}^g \frac{y_r}{G'(x_r)} ,
\label{3e4}\\
&\varphi_i=\frac{f(e_i)}{(R(e_i))^2}, 
\label{3e5}
\end{align}
Then we obtain $\Omega_{12}=(e_1-e_2)^2 \Delta_{12}^2-\varphi_1-\varphi_2$.
First, we obtain the expression in the limit $x_0 \rightarrow \infty$ 
\begin{align}
&R(e_i) \rightarrow -x_0 F(e_i), \ 
R'(x_r)  \rightarrow -x_0 F'(x_r), (1\le r \le g), \
\Delta_i \rightarrow -\frac{1}{x_0} \widehat{\Delta}_i+C_1, \ 
\nonumber\\
&\Delta_{12} \rightarrow -\frac{1}{x_0} \widehat{\Delta}_{12}, \ 
\varphi_i  \rightarrow \frac{1}{x_0^2} \widehat{\varphi}_i, 
\nonumber\\
&\text{where}
\nonumber\\
&\widehat{\Delta}_i=\sum_{r=1}^g \frac{y_r}{(e_i-x_r) F'(x_r)}, \ 
\widehat{\Delta}_{12}=\sum_{r=1}^g \frac{y_r}{(e_1-x_r)(e_2-x_r)F'(x_r)}, \
\widehat{\varphi}_i =\frac{f(e_i)}{(F(e_i))^2}.
\nonumber
\end{align}
Thus, the starting relation of Baker's paper~\cite{Baker4} is given by
\begin{align}
&4 (e_1-e_2)^2 \sum_{i, j=1}^g \wp_{i j} (u) e_1^{i-1} e_2^{j-1} -f_1(e_1, e_2)
=F(e_1) F(e_2)\left[(e_1-e_2)^2 \widehat{\Delta}_{12}^2-\widehat{\varphi}_1-\widehat{\varphi}_2\right]. 
\label{3e6}
\end{align}
By applying operators

\begin{align}
&\delta_3=\sum_{k=1}^g e_3^{k-1} \frac{\partial}{\partial u_k}, \quad
\delta_4=\sum_{k=1}^g e_4^{k-1} \frac{\partial}{\partial u_k}
\label{3e7}
\end{align}
to Eq.(\ref{3e6}), we obtain various differential equations.
For that purpose, by using Eq.(\ref{2e22}) and Eq.(\ref{2e9}), we obtain
the useful formula of $\delta_3$ in the form
\begin{align}
&\delta_3=\sum_{k=1}^g e_3^{k-1} \frac{\partial}{\partial u_k}
=\sum_{k=1}^g \sum_{r=1}^g e_3^{k-1} \frac{\partial x_r}{\partial u_k}\frac{\partial}{\partial x_r}
=\sum_{k=1}^g \sum_{r=1}^g e_3^{k-1} \frac{y_r \chi_{g-k}(x_r; x_1, \cdots, x_g) }{F'(x_r)}\frac{\partial}{\partial x_r}
\nonumber\\
&=\sum_{k=1}^g \left(\sum_{r=1}^g e_3^{k-1} \chi_{g-k}(x_r; x_1, \cdots, x_g) \right) \frac{y_r}{F'(x_r)} 
\frac{\partial}{\partial x_r}
=F(e_3)\sum_{k=1}^g  \frac{y_r}{(e_3-x_r)F'(x_r)} \frac{\partial}{\partial x_r}.
\label{3e8}
\end{align}
%
%
 
%%%%%%%%%%%%%%%%%%%%%%%%%%%%%%%%%%%%%%%%%%%%%%%%%%%
%%%%%%%%%%%%%%%%%%%%%%%%%%%%%%%%%%%%%%%%%%%%%%%%
\subsection{Step 1: Differential Equation of $\wp_{i j k}(u)$ }
First, we apply $\delta_3$ to Eq.(\ref{3e6}) and obtain
\begin{align}
&4 (e_1-e_2)^2 \sum_{i,  j, k=1}^g \wp_{i j k} (u) e_1^{i-1} e_2^{j-1} e_3^{k-1}
=\delta_3 \left(F(e_1) F(e_2)\left((e_1-e_2)^2 \widehat{\Delta}_{12}^2
-\widehat{\varphi}_1-\widehat{\varphi}_2\right) \right) .
\label{3e9}
\end{align}
In order to calculate the right-hand side of Eq.(\ref{3e9}), we must calculate
$\delta_3(F(e_i))$, $\delta_3(\widehat{\Delta}_i)$, $\delta_3(\widehat{\varphi}_i)$,
which gives
\begin{align}
&{\rm a)}\ \frac{1}{F(e_3)} \delta_3 (F(e_1))
=\frac{1}{F(e_3)} \times F(e_3) \sum_{r=1}^g \frac{y_r}{(e_3-x_r)F'(x_r)} \frac{\partial F(e_1)}{\partial x_r}
\nonumber\\
&=-F(e_1) \sum_{r=1}^g \frac{y_r}{(e_3-x_r)(e_1-x_r)F'(x_r)}=-F(e_1) \widehat{\Delta}_{13}, 
\label{3e10}\\
&{\rm b)}\  \frac{1}{F(e_3)} \delta_3 (\widehat{\varphi}_1)
=\frac{1}{F(e_3)} \times F(e_3) \sum_{r=1}^g \frac{y_r}{(e_3-x_r)F'(x_r)} 
\frac{\partial}{\partial x_r} \left( \frac{f(e_1)}{F(e_1)^2} \right)
\nonumber\\
&=2 \frac{f(e_1)}{F(e_1)^2} \sum_{r=1}^g \frac{y_r}{(e_3-x_r)(e_1-x_r)F'(x_r)}
=2 \widehat{\varphi}_1 \widehat{\Delta}_{12}  ,
\label{3e11}\\
&{\rm c)}\  \frac{1}{F(e_3)} \delta_3 (\widehat{\Delta}_1)
=\frac{1}{F(e_3)} \delta_3 \left(\sum_{r=1}^g \frac{y_r}{(e_1-x_r)F'(x_r)}\right)=J .
\label{3e12}
\end{align}
We calculate this $J$ in the following steps.
\begin{align}
&J=\frac{1}{F(e_3)} \sum_{r=1}^g \delta_3 \left(\frac{1}{e_1-x_r}\right) \frac{y_r}{F'(x_r)}
+\frac{1}{F(e_3)} \sum_{r=1}^g \frac{1}{e_1-x_r} \delta_3 \left(\frac{y_r}{F'(x_r)}\right)
=J_1+J_2, 
\label{3e13}\\
&J_1=\frac{1}{F(e_3)} \sum_{r=1}^g \delta_3 \left(\frac{1}{e_1-x_r}\right) \frac{y_r}{F'(x_r)}
=\sum_{r=1}^g \frac{f(x_r)}{(e_1-x_r)^2 (e_3-x_r) (F'(x_r))^2}, 
\label{3e14}\\
&J_2=\frac{1}{F(e_3)} \sum_{r=1}^g \frac{1}{e_1-x_r} \delta_3 \left(\frac{y_r}{F'(x_r)}\right)
\nonumber\\
&=\sum_{r=1}^g \frac{1}{e_1-x_r} \left(  \frac{y_r y'_r}{(e_3-x_r) (F'(x_r))^2 } 
-\frac{y_r^2 F''(x_r)}{(e_3-x_r) (F'(x_r))^3} 
+\sum_{\substack{s=1 \\ s \ne r}}^g \frac{y_r y_s }{(e_3-x_s)F'(x_r)F'(x_s)(x_r-x_s)} \right)
\nonumber\\
&=\frac{1}{2} \sum_{r=1}^g \frac{f'(x_r)}{(e_1-x_r)(e_3-x_r)(F'(x_r))^2}
-\sum_{r=1}^g \frac{f(x_r) F''(x_r)}{(e_1-x_r)(e_3-x_r) (F'(x_r))^3} 
\nonumber\\
&+\sum_{\substack{r, s=1 \\ r \ne s } }^g 
\frac{y_r y_s }{(e_1-x_r)(e_3-x_s)F'(x_r)F'(x_s)(x_r-x_s)}
\nonumber\\
&=\frac{1}{2} \sum_{r=1}^g  \frac{1}{(e_1-x_r)(e_3-x_r)}
\frac{\partial}{\partial x_r}\left( \frac{f(x_r)}{(F'(x_r))^2} \right)
\nonumber\\
&+\frac{(e_3-e_1)}{2}
\sum_{\substack{r, s=1 \\ r \ne s } }^g \frac{y_r }{(e_1-x_r)(e_3-x_r)F'(x_r)}
\frac{y_s }{(e_1-x_s)(e_3-x_s)F'(x_s)}
\nonumber\\
&=\frac{1}{2} \sum_{r=1}^g  \frac{1}{(e_1-x_r)(e_3-x_r)}
\frac{\partial}{\partial x_r}\left( \frac{f(x_r)}{(F'(x_r))^2} \right)
+\frac{(e_3-e_1)}{2}
\left[ \left(\sum_{r=1}^g \frac{y_r }{(e_1-x_r)(e_3-x_r)F'(x_r)} \right)^2 \right.
\nonumber\\
&\left. - \sum_{r=1) }^g \frac{f(x_r)}{(e_1-x_r)^2(e_3-x_r)^2 (F'(x_r))^2} \right]
=\frac{1}{2} \sum_{r=1}^g  \frac{1}{(e_1-x_r)(e_3-x_r)}
\frac{\partial}{\partial x_r}\left( \frac{f(x_r)}{(F'(x_r))^2} \right)
\nonumber\\
& - \frac{(e_3-e_1)}{2} \sum_{r=1 }^g \frac{f(x_r)}{(e_1-x_r)^2(e_3-x_r)^2 (F'(x_r))^2} 
+\frac{(e_3-e_1)}{2} {\widehat{\Delta}_{13}}^2   .
\label{3e15}
\end{align}
Then we obtain
\begin{align}
&J=J_1+J_2=\frac{1}{2} \sum_{r=1}^g  \frac{1}{(e_1-x_r)(e_3-x_r)}
\frac{\partial}{\partial x_r}\left( \frac{f(x_r)}{(F'(x_r))^2} \right)
+\frac{(e_3-e_1)}{2} {\widehat{\Delta}_{13}}^2 
\nonumber\\
&+\sum_{r=1}^g \frac{f(x_r)}{(e_1-x_r)^2 (e_3-x_r) (F'(x_r))^2}
 - \frac{(e_3-e_1)}{2} \sum_{r=1 }^g \frac{f(x_r)}{(e_1-x_r)^2(e_3-x_r)^2 (F'(x_r))^2} 
\nonumber\\
&=\frac{1}{2} \sum_{r=1}^g  \frac{1}{(e_1-x_r)(e_3-x_r)}
\frac{\partial}{\partial x_r}\left( \frac{f(x_r)}{(F'(x_r))^2} \right)
+\frac{(e_3-e_1)}{2} {\widehat{\Delta}_{13}}^2 
\nonumber\\
&+\frac{1}{2} \sum_{r=1}^g \frac{f(x_r)}{(e_1-x_r)^2 (e_3-x_r)^2 (F'(x_r))^2}
(2 e_3-2 x_r-e_3+e_1)
\nonumber\\
&=\frac{1}{2} \sum_{r=1}^g  \frac{1}{(e_1-x_r)(e_3-x_r)}
\frac{\partial}{\partial x_r}\left( \frac{f(x_r)}{(F'(x_r))^2} \right)
+\frac{(e_3-e_1)}{2} {\widehat{\Delta}_{13}}^2 
\nonumber\\
&+\frac{1}{2} \sum_{r=1}^g \frac{f(x_r)}{(e_1-x_r)^2 (e_3-x_r)^2 (F'(x_r))^2}
((e_1-x_r)+(e_3-x_r))
\nonumber\\
&=\frac{1}{2} \sum_{r=1}^g  \frac{1}{(e_1-x_r)(e_3-x_r)}
\frac{\partial}{\partial x_r}\left( \frac{f(x_r)}{(F'(x_r))^2} \right)
+\frac{(e_3-e_1)}{2} {\widehat{\Delta}_{13}}^2 
\nonumber\\
&+\frac{1}{2} \sum_{r=1}^g 
\frac{f(x_r)}{(F'(x_r))^2} \left( \frac{1}{(e_1-x_r) (e_3-x_r)^2}+\frac{1}{(e_1-x_r)^2 (e_3-x_r)}
\right)
\nonumber\\
&=\frac{1}{2} \sum_{r=1}^g  \frac{1}{(e_1-x_r)(e_3-x_r)}
\frac{\partial}{\partial x_r}\left( \frac{f(x_r)}{(F'(x_r))^2} \right)
+\frac{1}{2} \sum_{r=1}^g 
\frac{f(x_r)}{(F'(x_r))^2}  \frac{\partial}{\partial x_r}  \left( \frac{1}{(e_1-x_r) (e_3-x_r)}\right)
\nonumber\\
&+\frac{(e_3-e_1)}{2} {\widehat{\Delta}_{13}}^2
= \frac{1}{2} \sum_{r=1}^g  \frac{\partial}{\partial x_r}\left( \frac{f(x_r)}{(e_1-x_r)(e_3-x_r)(F'(x_r))^2} \right)+\frac{(e_3-e_1)}{2} {\widehat{\Delta}_{13}}^2 .
\label{3e16}
\end{align}
Hence, we have 
\begin{align}
&{\rm c)}\ \frac{1}{F(e_3)} \delta_3 (\widehat{\Delta}_1)
= \frac{1}{2} \sum_{r=1}^g  \frac{\partial}{\partial x_r}\left( \frac{f(x_r)}{(e_1-x_r)(e_3-x_r)(F'(x_r))^2} \right)+\frac{(e_3-e_1)}{2} {\widehat{\Delta}_{13}}^2
\nonumber\\
&=-\frac{(e_1-e_3)}{2} \widehat{\Delta}_{13}^2
-\frac{1}{2 (e_1-e_3)} (\widehat{\varphi}_1-\widehat{\varphi}_3)
+\frac{1}{2} \lambda_{2g+1} .
\label{3e17}
\end{align}
Here we use the formula
\begin{align}
\sum_{r=1}^g  \frac{\partial}{\partial x_r}\left( \frac{f(x_r)}{(e_1-x_r)(e_3-x_r)(F'(x_r))^2} \right)
=-\frac{1}{(e_1-e_3)} (\widehat{\varphi}_1-\widehat{\varphi}_3)
+\lambda_{2g+1}.
\label{3e18}
\end{align}
We have checked Eq,(\ref{3e18}) for genus two, three and four cases.
Next, we calculate $\displaystyle{\frac{1}{F(e_3)} \delta_3 (\widehat{\Delta}_{12})}$
in two ways, that is, by using Eq.(\ref{3e18}) and by not using Eq.(\ref{3e18})
in the form 
\begin{align}
&\frac{1}{F(e_3)} \delta_3 (\widehat{\Delta}_{12})=\frac{1}{2(e_1-e_2)}
\left((e_1-e_3)\widehat{\Delta}_{13}^2-(e_2-e_3)\widehat{\Delta}_{12}^2 \right)
\nonumber\\
&+\frac{1}{2} \sum_{r=1}^g  \frac{\partial}{\partial x_r}
\left( \frac{f(x_r)}{(e_1-x_r)(e_2-x_r)(e_3-x_r)(F'(x_r))^2} \right)
\label{3e19}\\
&=
\frac{1}{2(e_1-e_2)}
\left((e_1-e_3)\widehat{\Delta}_{13}^2-(e_2-e_3)\widehat{\Delta}_{12}^2 \right)
\nonumber\\
&+\frac{1}{2} \left( \frac{\widehat{\varphi}_1}{(e_1-e_2)(e_1-e_3)}
+\frac{\widehat{\varphi}_2}{(e_2-e_1)(e_2-e_3)}
+\frac{\widehat{\varphi}_3}{(e_3-e_1)(e_3-e_2)} \right) .
\label{3e20}
\end{align}
Then we obtain the following relation from Eq.(\ref{3e18})
\begin{align}
&\sum_{r=1}^g  \frac{\partial}{\partial x_r}
\left( \frac{f(x_r)}{(e_1-x_r)(e_2-x_r)(e_3-x_r)(F'(x_r))^2} \right)
\nonumber\\
&=\frac{\widehat{\varphi}_1}{(e_1-e_2)(e_1-e_3)}
+\frac{\widehat{\varphi}_2}{(e_2-e_1)(e_2-e_3)}
+\frac{\widehat{\varphi}_3}{(e_3-e_1)(e_3-e_2)}.
\label{3e21}
\end{align}
We have checked this relation for  genus two, three and four 
cases.~\footnote{In Baker's paper, the expression of the left-hand side of 
Eq.(\ref{3e21}) is in the incorrect form of 
$ \displaystyle{\sum_{r=1}^g  \frac{1}{F'(x_r)}\frac{\partial}{\partial x_r}
\left( \frac{f(x_r)}{(e_1-x_r)(e_2-x_r)(e_3-x_r)F'(x_r)} \right)} $.
}
Thus, we obtain
\begin{align}
&\frac{1}{F(e_3)} \delta_3 (\widehat{\Delta}_{12})=\frac{1}{2(e_1-e_2)}
\left((e_1-e_3)\widehat{\Delta}_{13}^2-(e_2-e_3)\widehat{\Delta}_{12}^2 \right)
\nonumber\\
&+\frac{1}{2} \left(\frac{\widehat{\varphi}_1}{(e_1-e_2)(e_1-e_3)}
+\frac{\widehat{\varphi}_2}{(e_2-e_1)(e_2-e_3)}
+\frac{\widehat{\varphi}_3}{(e_3-e_1)(e_3-e_2)} \right) .
\label{3e22}
\end{align}
Hence, using Eq.(\ref{3e10}) and Eq.(\ref{3e21}) , we obtain
\begin{align}
&\frac{1}{F(e_1) F(e_2) F(e_3)} \delta_3 
\left(F(e_1) F(e_2) (e_1-e_2)^2 \widehat{\Delta}_{12}^2 \right)
\nonumber\\
&=(e_1-e_2) \widehat{\Delta}_{12} 
\left( (e_1-e_3)\widehat{\Delta}_{13}^2-(e_2-e_3)\widehat{\Delta}_{23}^2 \right)
-(e_1-e_2)^2 \widehat{\Delta}_{12}^2(\widehat{\Delta}_{12}+\widehat{\Delta}_{23})
\nonumber\\
&+(e_1-e_2)^2 \widehat{\Delta}_{12} 
\left(\frac{\widehat{\varphi}_1}{(e_1-e_2)(e_1-e_3)}
+\frac{\widehat{\varphi}_2}{(e_2-e_1)(e_2-e_3)}
+\frac{\widehat{\varphi}_3}{(e_3-e_1)(e_3-e_2)} \right)
\nonumber\\
&=(e_1-e_2)\widehat{\Delta}_{12} \left( (e_1-e_3)\widehat{\Delta}_{13}^2-(e_2-e_3)
\widehat{\Delta}_{23}^2 -(e_1-e_2)\widehat{\Delta}_{12}
(\widehat{\Delta}_{13}+\widehat{\Delta}_{23}) \right)
\nonumber\\
&+(e_1-e_2)^2 \widehat{\Delta}_{12} 
\left(\frac{\widehat{\varphi}_1}{(e_1-e_2)(e_1-e_3)}
+\frac{\widehat{\varphi}_2}{(e_2-e_1)(e_2-e_3)}
+\frac{\widehat{\varphi}_3}{(e_3-e_1)(e_3-e_2)} \right)
\nonumber\\
&=-(e_1-e_2)^2 \widehat{\Delta}_{12} \widehat{\Delta}_{23} \widehat{\Delta}_{31} 
-(e_1-e_2)^2 \frac{(e_1-e_2) \widehat{\Delta}_{12} \widehat{\varphi}_3}
{(e_1-e_2)(e_2-e_3)(e_3-e_1)}
\nonumber\\
&+\text{( $\widehat{\varphi}_1$, $\widehat{\varphi}_2$ linear terms)}.
\label{3e23}
\end{align}
Here we use the identity
\begin{align}
&\left( (e_1-e_3)\widehat{\Delta}_{13}^2-(e_2-e_3)
\widehat{\Delta}_{23}^2 -(e_1-e_2)\widehat{\Delta}_{12}
(\widehat{\Delta}_{13}+\widehat{\Delta}_{23}) \right)
=-(e_1-e_2) \widehat{\Delta}_{23} \widehat{\Delta}_{31}, 
\nonumber
\end{align}
by using
$\displaystyle{\widehat{\Delta}_{ij}
=-\frac{ \widehat{\Delta}_i-\widehat{\Delta}_j}{e_i-e_j}}$.

Therefore, the right-hand side Eq.(\ref{3e9}) is given by
\begin{align}
&\frac{1}{F(e_1) F(e_2) F(e_3)} \delta_3 
\left(F(e_1) F(e_2)\left( (e_1-e_2)^2 \widehat{\Delta}_{12}^2 
-\widehat{\varphi}_1 -\widehat{\varphi}_2  \right)\right)
\nonumber\\
&=-(e_1-e_2)^2 \widehat{\Delta}_{12} \widehat{\Delta}_{23} \widehat{\Delta}_{31} 
-(e_1-e_2)^2 \frac{(e_1-e_2) \widehat{\Delta}_{12} \widehat{\varphi}_3}
{(e_1-e_2)(e_2-e_3)(e_3-e_1)}
\nonumber\\
&+\text{($\widehat{\varphi}_1$, $\widehat{\varphi}_2$ linear terms)}
\label{3e24}
\end{align}
where we used that $\delta_3 \left(F(e_1)F(e_2) (\widehat{\varphi}_1+\widehat{\varphi}_2)\right)$ is 
proportional to $\widehat{\varphi}_1$ or $\widehat{\varphi}_2$ from Eq.(\ref{3e11}).
But this term is symmetric in the index because $\delta_3 \delta_2 \delta_1( \log \sigma(u))$
is symmetric.
Then we finally obtain the right-hand side of Eq.(\ref{3e24}) in the form 
\begin{align}
&\hskip -11 mm \frac{1}{F(e_1) F(e_2) F(e_3)} \delta_3 
\left(F(e_1) F(e_2)\left( (e_1-e_2)^2 \widehat{\Delta}_{12}^2 
-\widehat{\varphi}_1 -\widehat{\varphi}_2  \right)\right)
\nonumber\\
&\hskip -11 mm =-(e_1-e_2)^2\left[
 \widehat{\Delta}_{12} \widehat{\Delta}_{23} \widehat{\Delta}_{31} 
+\frac{1}{\omega_{123}} \left( (e_1-e_2) \widehat{\Delta}_{12} \widehat{\varphi}_3
+(e_2-e_3) \widehat{\Delta}_{23} \widehat{\varphi}_1
+(e_3-e_1) \widehat{\Delta}_{31} \widehat{\varphi}_2
\right)\right].
\label{3e25}
\end{align}
with $\omega_{123}=(e_1-e_2)(e_2-e_3)(e_3-e_1)$.
Hence, we finally obtain the differential relation for $\wp_{ijk}(u)$ in the form  
\begin{align}
&-\frac{4}{F(e_1)F(e_2) F(e_3)} 
\sum_{i,  j, k=1}^g \wp_{i j k} (u) e_1^{i-1} e_2^{j-1} e_3^{k-1}
\nonumber\\
&=\left[
 \widehat{\Delta}_{12} \widehat{\Delta}_{23} \widehat{\Delta}_{31} 
+\frac{1}{\omega_{123}} \left( (e_1-e_2) \widehat{\Delta}_{12} \widehat{\varphi}_3
+(e_2-e_3) \widehat{\Delta}_{23} \widehat{\varphi}_1
+(e_3-e_1) \widehat{\Delta}_{31} \widehat{\varphi}_2 \right) \right].
\label{3e26}
\end{align} 
 
%%%%%%%%%%%%%%%%%%%%%%%%%%%%%%%%%%%%%%%%%%%%%%%%%%%
%%%%%%%%%%%%%%%%%%%%%%%%%%%%%%%%%%%%%%%%%%%%%%%%%%%
\subsection{Step 2: Differential Equations of $\wp_{i j k \ell}(u)$}
We apply $\delta_4$ to Eq.(\ref{3e26}) to obtain the differential equation.
By using the formulas of Eq.(\ref{3e10}), Eq.(\ref{3e11}) and Eq.(\ref{3e20}),
we can exactly obtain differential equations, but there comes 
many various terms. Then we use the simplified method to neglect the 
$\vvarphi_i(i=1, 2, 3, 4)$ terms. As $\delta_i( \vvarphi_j)$ is proportional to
$\vvarphi_j$, we can consistently neglect terms which is proportional to $\vvarphi_j$.
Then we start from the form
\begin{align}
&\hskip -5mm -4\sum_{i,  j, k=1}^g \wp_{i j k} (u) e_1^{i-1} e_2^{j-1} e_3^{k-1}
=F(e_1)F(e_2) F(e_3) \DDelta_{12} \DDelta_{23} \DDelta_{31} 
+\text{($\widehat{\varphi}_i$ linear terms)}.
\label{3e27}
\end{align} 
We apply $\delta_4$ to Eq.(\ref{3e27}) and obtain 
\begin{align}
&-4\sum_{i,  j, k=1}^g \wp_{i j k \ell } (u) e_1^{i-1} e_2^{j-1} e_3^{k-1} e_4^{\ell-1}
=\delta_4(F(e_1)F(e_2) F(e_3) \DDelta_{12} \DDelta_{23} \DDelta_{31} )
+\text{($\widehat{\varphi}_i$ linear terms)}
\nonumber\\
&=F(e_1)F(e_2)F(e_3)F(e_4) \left[
-(\DDelta_{14}+\DDelta_{24}+\DDelta_{34})\DDelta_{12}\DDelta_{23}\DDelta_{31} \right.
\nonumber\\
&\left.
+\frac{1}{2(e_1-e_2)}\left( (e_1-e_4)\DDelta_{14}^2-(e_2-e_4)\DDelta_{24}^2 \right)\DDelta_{23} \DDelta_{31} 
+\frac{1}{2(e_2-e_3)}\left( (e_2-e_4)\DDelta_{24}^2-(e_3-e_4)\DDelta_{34}^2 \right)\DDelta_{12} \DDelta_{31} 
\right.
\nonumber\\
& \left.
+\frac{1}{2(e_3-e_1)}\left( (e_3-e_4)\DDelta_{34}^2-(e_1-e_4)\DDelta_{14}^2 \right)\DDelta_{12} \DDelta_{23} 
\right]
+\text{($\widehat{\varphi}_i$ linear terms)}
\nonumber\\
&=-\frac{F(e_1)F(e_2)F(e_3)F(e_4) }{2(e_2-e_3)(e_3-e_1)(e_1-e_2)(e_4-e_1)(e_4-e_2)(e_4-e_3)}
\left( (e_2-e_3)^3 (e_4-e_1)^3 \DDelta_{23}^2 \DDelta_{41}^2 \right.
\nonumber\\
&\left. 
+(e_3-e_1)^3 (e_4-e_2)^3 \DDelta_{31}^2 \DDelta_{42}^2  
+(e_1-e_2)^3 (e_4-e_3)^3 \DDelta_{12}^2 \DDelta_{43}^2 \right) 
+\text{($\widehat{\varphi}_i$ linear terms)},
\label{3e28}
\end{align} 
where we used Eq.(\ref{3e10}), Eq.(\ref{3e20}) and the following identity
\begin{align}
&(e_2-e_3)^3 (e_4-e_1)^3 \DDelta_{23}^2 \DDelta_{41}^2 
+(e_3-e_1)^3 (e_4-e_2)^3 \DDelta_{31}^2 \DDelta_{42}^2  
+(e_1-e_2)^3 (e_4-e_3)^3 \DDelta_{12}^2 \DDelta_{43}^2
\nonumber\\
&=M \left[
2(\DDelta_{14}+\DDelta_{24}+\DDelta_{34})\DDelta_{12}\DDelta_{23}\DDelta_{31} 
-\DDelta_{23} \DDelta_{31} \frac{(e_1-e_4)\DDelta_{14}^2-(e_2-e_4)\DDelta_{24}^2 }{(e_1-e_2)}
\right.
\nonumber\\
&\left. -\DDelta_{12} \DDelta_{31} \frac{( (e_2-e_4)\DDelta_{24}^2-(e_3-e_4)\DDelta_{34}^2}{ (e_2-e_3)}
-\DDelta_{12} \DDelta_{23} \frac{(e_3-e_4)\DDelta_{34}^2-(e_1-e_4)\DDelta_{14}^2}{(e_3-e_1)}
\right], 
\label{3e29}\\
&\text{where}
\nonumber\\
&M=(e_2-e_3)(e_3-e_1)(e_1-e_2)(e_4-e_1)(e_4-e_2)(e_4-e_3).
\label{3e30}
\end{align} 
Thus we obtain
\begin{align}
&8M
\sum_{i,  j,  k, \ell =1}^g \wp_{i j k \ell} (u) e_1^{i-1} e_2^{j-1} e_3^{k-1} e_4^{\ell-1}
\nonumber\\
&=F(e_1)F(e_2)F(e_3)F(e_4)\left( (e_2-e_3)^3 (e_4-e_1)^3 \DDelta_{23}^2 \DDelta_{41}^2 \right.
\nonumber\\
&\left. 
+(e_3-e_1)^3 (e_4-e_2)^3 \DDelta_{31}^2 \DDelta_{42}^2  
+(e_1-e_2)^3 (e_4-e_3)^3 \DDelta_{12}^2 \DDelta_{43}^2 \right) 
+\text{($\widehat{\varphi}_i$ linear terms)}.
\label{3e31}
\end{align} 
Then we eliminate $\DDelta_{k \ell}$ by using 
\begin{align}
F(e_k)F(e_{\ell}) (e_k -e_{\ell})^2 \DDelta_{k \ell}^2
=4 (e_k -e_{\ell})^2 \sum_{i, j=1}^g \wp_{i j}(u) e_k^{i-1} e_{\ell}^{j-1}-f_1(e_k, e_{\ell})
+\text{($\widehat{\varphi}_i$ linear terms)}.
\label{3e32}
\end{align}
Hence, we obtain
\begin{align}
&\hskip -10mm 8M \sum_{i,  j, k, \ell =1}^g \wp_{i j k \ell} (u) e_1^{i-1} e_2^{j-1} e_3^{k-1} e_4^{\ell-1}
\nonumber\\
&\hskip -10mm =(e_2-e_3)(e_4-e_1)
\left[f_1(e_2, e_3)-4(e_2-e_3)^2\sum_{i, j=1}^g \wp_{i j}(u) e_2^{i-1} e_3^{j-1}\right]
\left[f_1(e_4, e_1)-4(e_4-e_1)^2\sum_{i, j=1}^g \wp_{i j}(u) e_4^{i-1} e_1^{j-1}\right]
\nonumber\\
&\hskip -10mm +(e_3-e_1)(e_4-e_2)
\left[f_1(e_3, e_1)-4(e_3-e_1)^2\sum_{i, j=1}^g \wp_{i j}(u) e_3^{i-1} e_1^{j-1}\right]
\left[f_1(e_4, e_2)-4(e_4-e_2)^2\sum_{i, j=1}^g \wp_{i j}(u) e_4^{i-1} e_2^{j-1}\right]
\nonumber\\
&\hskip -10mm +(e_1-e_2)(e_4-e_3)
\left[f_1(e_1, e_2)-4(e_1-e_2)^2\sum_{i, j=1}^g \wp_{i j}(u) e_1^{i-1} e_2^{j-1}\right]
\left[f_1(e_4, e_3)-4(e_4-e_3)^2\sum_{i, j=1}^g \wp_{i j}(u) e_4^{i-1} e_3^{j-1}\right] .
\label{3e33}
\end{align} 
By subtracting the following term from both side of Eq.(\ref{3e33})  
\begin{align}
&16M
\sum_{i,j,k,\ell=1}^g (\wp_{i \ell}(u)\wp_{j k}(u)+\wp_{i k}(u)\wp_{j \ell}(u)+\wp_{i j}(u)\wp_{k \ell}(u))
e_1^{i-1} e_2^{j-1}e_3^{k-1} e_4^{\ell-1}, 
\label{3e34}
\end{align}
we obtain
\begin{align}
&8M \sum_{i,  j, k, \ell=1}^g 
\Big(\wp_{i j k \ell} (u) -2(\wp_{i \ell}(u)\wp_{j k}(u)+\wp_{i k}(u)\wp_{j \ell}(u)+\wp_{i j}(u)\wp_{k \ell}(u))
\Big)
e_1^{i-1} e_2^{j-1} e_3^{k-1} e_4^{\ell-1}
\nonumber\\
&=16 \sum_{i,  j, k, \ell=1}^g 
\Big[ \left\{ (e_2-e_3)^3 (e_4-e_1)^3-M \right\} \wp_{i \ell}(u) \wp_{jk}(u)
+\left\{ (e_3-e_1)^3 (e_4-e_2)^3-M\right\} \wp_{i k}(u) \wp_{j \ell}(u)
\nonumber\\
&+\left\{ (e_1-e_2)^3 (e_4-e_3)^3-M\right\} \wp_{i j}(u) \wp_{k \ell}(u)
\Big] e_1^{i-1} e_2^{j-1} e_3^{k-1} e_4^{\ell-1}
\nonumber\\
&+\Big[(e_2-e_3)(e_4-e_1)f_1(e_2,e_3)f_1(e_4, e_1)+(e_3-e_1)(e_4-e_2)f_1(e_3,e_1)f_1(e_4, e_2)
\nonumber\\
&+(e_1-e_2)(e_4-e_3)f_1(e_1,e_2)f_1(e_4, e_3) \Big]
\nonumber\\
&+\Big[-4(e_2-e_3)(e_4-e_1)^3 f_1(e_2, e_3) \sum_{i, \ell=1}^g \wp_{i \ell}(u) e_1^{i-1} e_4^{\ell-1}
-4(e_2-e_3)^3 (e_4-e_1) f_1(e_4, e_1) \sum_{j, k=1}^g \wp_{j k}(u) e_2^{j-1} e_3^{k-1}
\nonumber\\
&-4(e_3-e_1)(e_4-e_2)^3 f_1(e_3, e_1) \sum_{j, \ell=1}^g \wp_{j \ell}(u) e_2^{j-1}e_4^{\ell-1}
-4(e_3-e_1)^3 (e_4-e_2) f_1(e_4, e_2) \sum_{i, k=1}^g \wp_{i k}(u) e_1^{i-1} e_3^{k-1}
\nonumber\\
&-4(e_1-e_2)(e_4-e_3)^3 f_1(e_1, e_2) \sum_{k, \ell=1}^g \wp_{k \ell}(u) e_3^{k-1} e_4^{\ell-1}
-4(e_1-e_2)^3 (e_4-e_3) f_1(e_4, e_3) \sum_{i, j=1}^g \wp_{i j}(u) e_1^{i-1} e_2^{j-1} \Big]
\nonumber\\
&=K_1+K_2+K_3 .
\label{3e35}
\end{align} 
From here, we change the notation in the form $e_1, e_2, e_3, e_4 \rightarrow x, y, z, t$.
Further we denote
$f_1(x,y)=f_{12}=f_{21}$, $f_1(y,z)=f_{23}=f_{32}$, $f_1(z,x)=f_{31}=f_{13}$, $f_1(x,t)=f_{14}=f_{41}$,
$f_1(y,t)=f_{24}=f_{42}$, $f_1(z,t)=f_{34}=f_{43}$ and we denote 
\begin{align}
&\sum_{i, j=1}^g \wp_{i j}(u)e_1^{i-1} e_2^{j-1}=\sum_{i, j=1}^g \wp_{i j}(u) x^{i-1} y^{j-1}=P_{12}=P_{21},
\label{3e36}\\
&\sum_{i, j=1}^g \wp_{i j}(u)e_1^{i-1} e_3^{j-1}=\sum_{i, j=1}^g \wp_{i j}(u) x^{i-1} z^{j-1}=P_{13}=P_{31},
\label{3e37}\\
&\sum_{i, j=1}^g \wp_{i j}(u)e_1^{i-1} e_4^{j-1}=\sum_{i, j=1}^g \wp_{i j}(u) x^{i-1} t^{j-1}=P_{14}=P_{41},
\label{3e38}\\
&\sum_{i, j=1}^g \wp_{i j}(u)e_2^{i-1} e_3^{j-1}=\sum_{i, j=1}^g \wp_{i j}(u) y^{i-1} z^{j-1}=P_{23}=P_{32},
\label{3e39}\\
&\sum_{i, j=1}^g \wp_{i j}(u)e_2^{i-1} e_4^{j-1}=\sum_{i, j=1}^g \wp_{i j}(u) y^{i-1} t^{j-1}=P_{24}=P_{42},
\label{3e40}\\
&\sum_{i, j=1}^g \wp_{i j}(u)e_3^{i-1} e_4^{j-1}=\sum_{i, j=1}^g \wp_{i j}(u) z^{i-1} t^{j-1}=P_{34}=P_{43},
\label{3e41}\\
&\text{where}
\nonumber\\
&M=(x-y)(y-z)(z-x)(t-x)(t-y)(t-z)
=\left| \begin{array}{@{\,}cccc@{\,}}
 x^3  & y^3  & z^3 & t^3     \\
 x^2 & y^2  & z^2  & t^2    \\
 x  & y  & z & t  \\
 1  & 1 & 1 & 1 \\
\end{array} \right| .
\label{3e42}
\end{align}
We denote the totally symmetric combination of $x^{\alpha} y^{\beta} z^{\gamma} t^{\delta}$
with respect to $x, y, z, t$, as 
\begin{align}
S(x^{\alpha} y^{\beta} z^{\gamma} t^{\delta})
=(\alpha \beta \gamma \delta) .
\nonumber
\end{align}
For example,
$S(x^2)=x^2+y^2+z^2+t^2=(2000)$, $S(xy)=xy+yz+zx+xt+yt+zt=(1100)$ etc. 
The first term $K_1$ in the right-hand side of Eq.(\ref{3e35}) gives
\begin{align}
& K_1=16(x-y)(t-z)\left\{ (x-y)^2 (t-z)^2-(y-z)(z-x)(t-x)(t-y) \right\} P_{12} P_{34}
\nonumber\\
&+16(y-z)(t-x)\left\{ (y-z)^2 (t-x)^2-(z-x)(x-y)(t-y)(t-z) \right\} P_{23} P_{14}
\nonumber\\
&+16(z-x)(t-y)\left\{ (z-x)^2 (t-y)^2-(x-y)(y-z)(t-z)(t-x) \right\} P_{31} P_{24}   
\nonumber\\
&=16H \Big( (x-y)(t-z) P_{12} P_{34}+(y-z)(t-x) P_{23} P_{14}+(z-x)(t-y) P_{31}  P_{24} \Big)
\nonumber\\
&=16H S((x-y)(t-z) P_{12} P_{34}), 
\label{3e43}
\end{align}
where
\begin{align}
&H=(x-y)^2 (t-z)^2-(y-z)(z-x)(t-x)(t-y)  
\nonumber\\
&=(x^2 y^2+y^2 z^2+z^2 x^2+t^2 x^2+t^2 y^2+t^2 z^2)-(x^2 yz +x^2 y t+x^2 z t+y^2 x z+y^2 x t+y^2 z t
\nonumber\\
&+z^2 x y +z^2 x t+z^2 y t+t^2 x y+t^2 y z+t^2 z x) +6 x y z t
\nonumber\\
&=(2200)-(2110)+6(1111) .
\label{3e44}
\end{align}
The second term $K_2$ in the right-hand side of Eq.(\ref{3e35}) gives
\begin{align}
&K_2=(x-y)(t-z) f_{12} f_{34}+(y-z)(t-x) f_{23} f_{14}+(z-x)(t-y) f_{31} f_{24}
\nonumber\\
&=S((x-y)(t-z) f_{12} f_{34} ).
\label{3e45}
\end{align}
The third term $K_3$ in the right-hand side of Eq.(\ref{3e35}) gives
\begin{align}
&K_3=-4(x-y)(t-z) \left\{ (t-z)^2 f_{12} P_{34}+(x-y)^2 f_{34} P_{12} \right\}
\nonumber\\
&-4(y-z)(t-x) \left\{ (t-x)^2 f_{23} P_{14}+(y-z)^2 f_{14} P_{23} \right\}
\nonumber\\
&-4(z-x)(t-y) \left\{ (t-y)^2 f_{31} P_{24}+(z-x)^2 f_{24} P_{31} \right\}
\nonumber\\
&=-4 S((x-y)(t-z) \left\{ (x-y)^2 P_{12}  f_{34}\right\}).
\label{3e46}
\end{align}
Up to this stage, we obtain 
\begin{align}
&8M \sum_{\alpha,  \beta, \gamma, \delta=1}^g (\alpha-1, \beta-1, \gamma-1, \delta-1)
 [\alpha \beta \gamma \delta]
\nonumber\\
&=16 H S((x-y) P_{12} (t-z) P_{34})
+S((x-y) f_{12} (t-z) f_{34} )
\nonumber\\
&-4 S((x-y)(t-z) \left\{ (x-y)^2 P_{12} f_{34} \right\}) ,
\label{3e47}
\end{align}
where
\begin{align}
&[\alpha \beta \gamma \delta]=[\wp_{\alpha \beta \gamma \delta}]
=\wp_{\alpha \beta \gamma \delta}(u)
-2(\wp_{\alpha \beta}(u) \wp_{\gamma \delta}(u)+\wp_{\alpha \gamma}(u) \wp_{\beta \delta}(u)
+\wp_{\alpha \delta}(u) \wp_{\beta \gamma}(u)),
\label{3e48}\\
&(\alpha-1,\beta-1, \gamma-1, \delta-1)=S(x^{\alpha-1} y^{\beta-1} z^{\gamma-1} t^{\delta-1}), 
\label{3e49}\\
&H=(2200)-(2110)+6(1111) .
\label{3e50}
\end{align}
We further simplify each term in the right-hand side of Eq.(\ref{3e47}). We put 
\begin{align}
&P_{12}=P_{21}=\sum_{\alpha, \beta=1}^g \wp_{\alpha+1, \beta+1}(u) x^{\alpha} y^{\beta}
=\sum_{\alpha, \beta}a_{\alpha \beta} x^{\alpha} y^{\beta}, 
\label{3e51}\\
&\text{where}
\nonumber\\
&a_{\alpha \beta}=\wp_{\alpha+1, \beta+1}. 
\label{3e52}
\end{align}
Then we obtain
\begin{align}
&(x-y) P_{12}=-(y-x)P_{21}=\sum_{\alpha, \beta} (x-y) \wp_{\alpha+1, \beta+1} x^{\alpha} y^{\beta}
=\sum_{\alpha, \beta}(\wp_{\alpha, \beta+1} -\wp_{\alpha+1, \beta} ) x^{\alpha} y^{\beta}
\nonumber\\
&=-\sum_{\alpha, \beta} \langle \alpha \beta \rangle  x^{\alpha} y^{\beta}, 
\label{3e53}\\
&\text{where}
\nonumber\\
&\langle \alpha \beta \rangle = -\langle \beta \alpha \rangle 
=\wp_{\alpha+1, \beta}(u)-\wp_{\alpha, \beta+1}(u) .
\label{3e54}
\end{align}
Similarly, we obtain
\begin{align}
&(t-z) P_{34}=-\sum_{\gamma, \delta} (z-t) \wp_{\gamma+1, \delta+1} z^{\gamma} t^{\delta}
=\sum_{\gamma, \delta} \langle \gamma \delta \rangle  z^{\gamma} t^{\delta}, 
\label{3e55}\\ 
&\text{where}
\nonumber\\
& \langle \gamma \delta \rangle = -\langle \delta \gamma \rangle 
=\wp_{\gamma+1, \delta}(u)-\wp_{\gamma, \delta+1}(u) .
\label{3e56}
\end{align}
Then the first term in the right-hand side of Eq.(\ref{3e47}) gives
\begin{align}
& 16 H S((x-y) P_{12} (t-z) P_{34})=-16H \sum S(\langle \alpha \beta \rangle \langle \gamma \delta \rangle
x^{\alpha} y^{\beta} z^{\gamma} t^{\delta})
\nonumber\\
&=-16H \sum S\left((\langle \alpha \beta \rangle \langle \gamma \delta \rangle
-\langle \alpha \gamma \rangle \langle \beta \delta \rangle
+\langle \alpha \delta \rangle \langle \beta \gamma \rangle)
x^{\alpha} y^{\beta} z^{\gamma} t^{\delta} \right)
\nonumber\\
& =-16 M H \sum \langle \alpha \beta \gamma \delta \rangle |\alpha \beta \gamma \delta|, 
\label{3e57}\\
&\text{where}
\nonumber\\
&\langle \alpha \beta \gamma \delta \rangle
=\langle \alpha \beta \rangle \langle \gamma \delta \rangle
-\langle \alpha \gamma \rangle \langle \beta \delta \rangle
+\langle \alpha \delta \rangle \langle \beta \gamma \rangle
\nonumber\\
&=(\wp_{\alpha+1,\beta}(u)- \wp_{\alpha, \beta+1}(u))
(\wp_{\gamma+1,\delta}(u)- \wp_{\gamma, \delta+1}(u))
\nonumber\\
&-(\wp_{\alpha+1,\gamma}(u)- \wp_{\alpha, \gamma+1}(u))
(\wp_{\beta+1,\delta}(u)- \wp_{\beta, \delta+1}(u))
\nonumber\\
&+(\wp_{\alpha+1,\delta}(u)- \wp_{\alpha, \delta+1}(u))
(\wp_{\beta+1,\gamma}(u)- \wp_{\beta, \gamma+1}(u)) ,
\label{3e58}
\end{align}
is totally anti-symmetric with respect to $\alpha, \beta, \gamma, \delta$ and 
\begin{align}
|\alpha \beta \gamma \delta|
=\frac{1}{M}\left| \begin{array}{@{\,}cccc@{\,}}
 x^{\alpha}  & y^{\alpha}  & z^{\alpha}  & t^{\alpha} \\
 x^{\beta}  & y^{\beta}  & z^{\beta}  & t^{\beta}   \\
 x^{\gamma}  & y^{\gamma}  & z^{\gamma}  & t^{\gamma}    \\
 z^{\delta}  & y^{\delta}  & z^{\delta}  & t^{\delta}     \\
\end{array} \right|  .
\label{3e59}
\end{align}
The second term in the right-hand side of Eq.(\ref{3e47}) has the same structure.
Then we define
\begin{align}
&f_{12}=\sum_{\alpha=0}^g \left[ 2 \lambda_{2 \alpha} +\lambda_{2 \alpha+1} (x+y)\right]
x^{\alpha} y^{\alpha}
\nonumber\\
&=\sum_{\alpha=0}^g \left[ 2 \lambda_{2 \alpha} x^{\alpha} y^{\alpha}
 +\lambda_{2 \alpha+1} (x^{\alpha+1}y^{\alpha}+x^{\alpha} y^{\alpha+1} )
\right]
=\sum_{\alpha, \beta } b_{\alpha \beta} x^{\alpha} y^{\beta}, 
\label{3e60}\\
&\text{where}
\nonumber\\
&b_{\alpha \alpha}=2 \lambda_{2 \alpha}, b_{\alpha, \alpha+1}=b_{\alpha+1, \alpha}=\lambda_{2 \alpha+1},
\ \text{Others}=0 .
\label{3e61}
\end{align}
In this case, we construct
\begin{align}
&(x-y) f_{12}= \sum_{\alpha, \beta } (x-y) b_{\alpha \beta} x^{\alpha} y^{\beta}
=\sum_{\alpha, \beta } \left[  b_{\alpha \beta} x^{\alpha+1} y^{\beta}
-b_{\alpha \beta} x^{\alpha} y^{\beta+1} \right]
=-\sum_{\alpha, \beta} \left\{ \alpha \beta \right\} x^{\alpha} y^{\beta},
\label{3e62}\\
&\text{where}
\nonumber\\
&\left\{ \alpha \beta \right\}=b_{\alpha, \beta-1}-b_{\alpha-1, \beta},
\label{3e63}\\
&\left\{ \alpha, \alpha+1 \right\}=-\left\{ \alpha+1, \alpha \right\}=2 \lambda_{2\alpha},
\left\{ \alpha, \alpha+2 \right\}=-\left\{ \alpha+2, \alpha \right\}= \lambda_{2\alpha+1},
\ \text{Others}=0 .
\label{3e64}
\end{align}
Then the second term in the right-hand side of Eq.(\ref{3e47}) becomes
\begin{align}
&S((x-y)f_{12} (t-z)f_{34})=-S\Big(\left\{ \alpha \beta \right\}\left\{ \gamma \delta \right\}
x^{\alpha}y^{\beta} z^{\gamma} t^{\delta}\Big)
\nonumber\\
&=-S\left( (\left\{ \alpha \beta \right\}\left\{ \gamma \delta \right\}
-\left\{ \alpha \gamma \right\}\left\{ \beta \delta \right\}
+\left\{ \alpha \delta \right\}\left\{ \beta \gamma \right\}) x^{\alpha}y^{\beta} z^{\gamma} t^{\delta}
\right)
=-M \sum  \left\{ \alpha \beta \gamma \delta \right\} |\alpha \beta \gamma \delta|, 
\label{3e65}\\
&\text{where}
\nonumber\\
& \left\{ \alpha \beta \gamma \delta \right\} 
=\left\{ \alpha \beta \right\}\left\{ \gamma \delta \right\}
-\left\{ \alpha \gamma \right\}\left\{ \beta \delta \right\}
+\left\{ \alpha \delta \right\}\left\{ \beta \gamma \right\},
\label{3e66}\\
&\left\{ \alpha, \alpha+1 \right\}=-\left\{ \alpha+1, \alpha \right\}=2 \lambda_{2\alpha},
\left\{ \alpha, \alpha+2 \right\}=-\left\{ \alpha+2, \alpha \right\}= \lambda_{2\alpha+1},
\ \text{Others}=0 .
\label{3e67}
\end{align}
In order to obtain the third term in the right-hand side of Eq.(\ref{3e47}), we must first 
calculate the symmetric function $(x-y)^2 P_{12}$ and obtain

\begin{align}
&(x-y)^2 P_{12}= \sum_{\alpha, \beta } (x-y)^2 a_{\alpha \beta} x^{\alpha} y^{\beta}
\nonumber\\
&=\sum_{\alpha, \beta } \left[  a_{\alpha \beta} x^{\alpha+2} y^{\beta}
+a_{\alpha \beta} x^{\alpha} y^{\beta+2} 
-2a_{\alpha \beta} x^{\alpha+1} y^{\beta+1}  \right]
=\sum_{\alpha, \beta} c_{\alpha \beta} x^{\alpha} y^{\beta}, 
\label{3e68}\\
&\text{where}
\nonumber\\
&c_{\alpha \beta}=\wp_{\alpha-1,\beta+1}(u)-2 \wp_{\alpha,\beta}(u)+\wp_{\alpha+1,\beta-1}(u).
\label{3e69}
\end{align}
Then we calculate
\begin{align}
&\hskip -5 mm (x-y) \cdot (x-y)^2 P_{12}
=-\sum_{\alpha, \beta} \langle\langle \alpha \beta \rangle \rangle x^{\alpha} y^{\beta}, 
\label{3e70}\\
&\hskip -5 mm \text{where}
\nonumber\\
&\hskip -5 mm\langle\langle \alpha \beta \rangle \rangle
=c_{\alpha, \beta-1}-c_{\alpha-1, \beta}
=-(\wp_{\alpha-2, \beta+1}(u)-3\wp_{\alpha-1, \beta}(u)
+3\wp_{\alpha, \beta-1}(u)-\wp_{\alpha+1, \beta-2}(u)). 
\label{3e71}
\end{align}
While we have 
\begin{align}
&(t-z)f_{34}=\sum \left\{ \gamma \delta \right\} z^{\gamma} t^{\delta}.
\label{3e72}
\end{align}
Using Eq.(\ref{3e70}) and Eq.(\ref{3e72}), we obtain the third term in the right-hand 
side of Eq.(\ref{3e47}) in the form
\begin{align}
&-4 S((x-y)(t-z) \Big((x-y)^2 P_{12} f_{34} \Big)
=4S \left( \langle\langle \alpha \beta \rangle \rangle 
\left\{ \gamma \delta \right\} x^{\alpha} y^{\beta} z^{\gamma} t^{\delta}
\right)
\nonumber\\
&=4M \sum \langle\langle \alpha \beta \gamma \delta \rangle \rangle 
|\alpha \beta \gamma \delta| , 
\label{3e73}\\
&\text{where}
\nonumber\\
&\langle\langle \alpha \beta \gamma \delta \rangle \rangle 
=\langle\langle \alpha \beta \rangle \rangle \left\{\gamma \delta \right\}
-\langle\langle \alpha \gamma \rangle \rangle \left\{\beta \delta \right\}
+\langle\langle \alpha \delta \rangle \rangle \left\{\beta \gamma \right\}
\nonumber\\
&+\left\{ \alpha \beta \right\}  \langle\langle \gamma \delta \rangle \rangle
-\left\{ \alpha \gamma \right\} \langle\langle \beta \delta \rangle \rangle
+\left\{ \alpha \delta \right\}  \langle\langle \beta \gamma \rangle \rangle .
\label{3e74}
\end{align}
Hence, we finally obtain the following relation
%%%%%%%%%%%%%%%%%%%%%%%%%%%%%%%%%%%%%%%%
%
\begin{align}
&\sum (\mu-1, \nu-1, \rho-1, \sigma-1)
 [\mu \nu \rho \sigma]
\nonumber\\
&=-2 H\sum \langle \alpha \beta \gamma \delta \rangle |\alpha \beta \gamma \delta|
+\frac{1}{2}\sum \langle \langle \alpha \beta \gamma \delta \rangle \rangle |\alpha \beta \gamma \delta|
-\frac{1}{8}\sum \left\{\alpha \beta \gamma \delta \right\}  |\alpha \beta \gamma \delta| .
\label{3e75}
\end{align}
%
%%%%%%%%%%%%%%%%%%%%%%%%%%%%%%%%%%%%%%%%%%%%%%%%
%%%%%%%%%%%%%%%%%%%%%%%%%%%%%%%%%%%%%%%%%%%%%
\section{Differential Equations of Genus Four Hyperelliptic $\wp$ Functions} 
\setcounter{equation}{0}
The differential equations for general genus $g$ hyperelliptic $\wp_{ij}(u)$ functions 
are obtained from the following relation
%%%%%%%%%%%%%%%%%%%%%%%%%%%%%%%%%%%%%%%%
%
\begin{align}
&\sum (\mu-1, \nu-1, \rho-1, \sigma-1)
 [\mu \nu \rho \sigma]
\nonumber\\
&=-2 H\sum \langle \alpha \beta \gamma \delta \rangle |\alpha \beta \gamma \delta|
+\frac{1}{2}\sum \langle \langle \alpha \beta \gamma \delta \rangle \rangle |\alpha \beta \gamma \delta|
-\frac{1}{8}\sum \left\{\alpha \beta \gamma \delta \right\}  |\alpha \beta \gamma \delta|, 
\label{4e1}
\end{align}
\begin{align}
&\text{where} \nonumber\\
& (\mu \nu \rho \sigma)=S(x^{\mu} y^{\nu} z^{\rho} t^{\sigma})
=\text{(cyclic symmetrization of $x$, $y$, $z$ and $t$)},
\label{4e2}\\
& [\mu \nu \rho \sigma]= [\mu \nu \rho \sigma]=\wp_{\mu \nu \rho \sigma}
-2(\wp_{\mu \nu} \wp_{\rho \sigma}+\wp_{\mu \rho} \wp_{\nu \sigma}
+\wp_{\mu \sigma} \wp_{\nu \rho}),
\label{4e3}\\
& \{ \alpha \beta \gamma \delta \}=\{\alpha \beta\} \{\gamma \delta \}
-\{\alpha \gamma\} \{\beta \delta\}+\{\alpha \delta \} \{ \beta \gamma \}, 
\label{4e4}\\
&\langle\langle \alpha \beta \gamma \delta \rangle \rangle 
=\langle\langle \alpha \beta \rangle \rangle \left\{\gamma \delta \right\}
-\langle\langle \alpha \gamma \rangle \rangle \left\{\beta \delta \right\}
+\langle\langle \alpha \delta \rangle \rangle \left\{\beta \gamma \right\}
\nonumber\\
&+\left\{ \alpha \beta \right\}  \langle\langle \gamma \delta \rangle \rangle
-\left\{ \alpha \gamma \right\} \langle\langle \beta \delta \rangle \rangle
+\left\{ \alpha \delta \right\}  \langle\langle \beta \gamma \rangle \rangle , 
\label{4e5}\\
& \langle \alpha \beta \gamma \delta \rangle
=\langle \alpha \beta \rangle  \langle \gamma \delta \rangle
-\langle \alpha \gamma\rangle \langle \beta \delta \rangle
+\langle \alpha \delta \rangle \langle \beta \gamma \rangle  ,
\label{4e6}\\
&\langle \alpha \beta \rangle=\wp_{\alpha+1,\beta}(u)-\wp_{\alpha,\beta+1}(u), 
\label{4e7}\\
&\left\{ \alpha, \alpha+1 \right\}=-\left\{ \alpha+1, \alpha \right\}=2 \lambda_{2 \alpha},
\left\{ \alpha, \alpha+2 \right\}=-\left\{ \alpha+2, \alpha \right\}= \lambda_{2 \alpha+1},
\text{Others}=0 ,
\label{4e8}\\
&\langle \langle \alpha \beta \rangle \rangle=-\Big(\wp_{\alpha-2, \beta+1}(u)
-3\wp_{\alpha-1, \beta}(u)
+3\wp_{\alpha, \beta-1}(u)-\wp_{\alpha+1, \beta-2}(u)  \Big) , 
\label{4e9}\\
&|\alpha \beta \gamma \delta|
=\frac{1}{M}\left| \begin{array}{@{\,}cccc@{\,}}
 x^{\alpha}  & y^{\alpha}  & z^{\alpha}  & t^{\alpha} \\
 x^{\beta}  & y^{\beta}  & z^{\beta}  & t^{\beta}   \\
 x^{\gamma}  & y^{\gamma}  & z^{\gamma}  & t^{\gamma}    \\
 z^{\delta}  & y^{\delta}  & z^{\delta}  & t^{\delta}     \\
\end{array} \right|  ,\quad 
M=\left| \begin{array}{@{\,}cccc@{\,}}
 x^3  & y^3  & z^3 & t^3     \\
 x^2 & y^2  & z^2  & t^2    \\
 x  & y  & z & t  \\
 1  & 1 & 1 & 1 \\
\end{array} \right| ,
\label{4e10}\\ 
&H=(2200)-(2110)+6(1111)  .
\label{4e11}
\end{align}
%%%%%%%%%%%%%%%%%%%%%%%%%%%%%%%%%%%%%%%%%%%%%%%%%%%%%%%%

For the genus $g$ differential equations, each differential equation 
contain terms $[\mu \nu \rho \sigma ]$$=$
$[ \wp_{\mu \nu \rho \sigma} ]$$=$
$\wp_{\mu \nu \rho \sigma}
-2(\wp_{\mu \nu} \wp_{\rho \sigma}+\wp_{\mu \rho} \wp_{\nu \sigma}
+\wp_{\mu \sigma} \wp_{\nu \rho})$
, which appear  in the left-hand side of 
Eq.(\ref{4e1}). Then $[\mu \nu \rho \sigma ]$$(g \ge \mu \ge \nu \ge \rho \ge \sigma \ge 1)$
 takes in the region $[gggg] \sim [1111]$.
Thus, we compare the coefficients of 
$\left(\mu-1, \nu-1, \rho-1, \sigma-1 \right)$ in the region 
$\left(g-1, g-1, g-1, g-1 \right) \sim \left(0 0 0 0 \right)$
in both side of Eq.(\ref{4e1}), which gives $_gH_4=_{g+3} C_4$ pieces of differential
equations. For $g=2$ case, we have $_5 C_4=5$ pieces of differential equations. 
For $g=3$ case, we have $_6 C_4=15$ pieces of differential equations. 
For $g=4$ case, we have $_7 C_4=35$ pieces of differential equations. 

In order to find the differential equations  for the general 
$|\alpha \beta \gamma \delta|$, we must know 
how many $\left(\mu-1, \nu-1, \rho-1, \sigma-1 \right)$ is contained 
in $|\alpha \beta \gamma \delta |$ in the right-hand side of Eq.(\ref{4e1}).
We know that $(g-1, g-1, g-1, g-1)$ is the $(4 g-4)$ degree polynomial of $x, y, z, t$. 
While $|\alpha \beta \gamma \delta|$$(\alpha > \beta >\gamma >\delta)$ is the 
$(\alpha+\beta+\gamma+\delta-6)$ degree
polynomial. Using $\alpha-1 \ge \beta$, $\alpha-2 \ge \gamma$, $\alpha-3 \ge \delta$,
we obtain $4 \alpha -6\ge (\alpha +\beta +\gamma+\delta)$.
Thus, if $|\alpha \beta \gamma \delta|$ contains $\left(\mu-1, \nu-1, \rho-1, \sigma-1 \right)$ 
, it must satisfy $ \alpha +\beta +\gamma+\delta-6 \le 4 \alpha -12 \le 4 g -4$, which
gives $\alpha \le g+2$.
Then it is sufficient to consider $|\alpha \beta \gamma \delta|$ only 
in the region $ g+2 \ge \alpha > \beta > \gamma > \delta \ge 0$.
That is, it is sufficient to calculate $_{g+3} C_4$ terms for $|\alpha \beta \gamma \delta|$
in the right-hand side of Eq.(\ref{4e1}).
  
For the first term in the right-hand side of Eq.(\ref{4e1}), it takes in the form
\begin{align}
&\langle \alpha \beta \gamma \delta \rangle
 =(\wp_{\alpha+1, \beta}(u)-\wp_{\alpha, \beta+1}(u))(\wp_{\gamma+1, \delta}(u)-\wp_{\gamma, \delta+1}(u))
\nonumber\\
&-(\wp_{\alpha+1, \gamma}(u)-\wp_{\alpha, \gamma+1}(u))(\wp_{\beta+1, \delta}(u)-\wp_{\beta, \delta+1}(u))
\nonumber\\
&+(\wp_{\alpha+1, \delta}(u)-\wp_{\alpha, \delta+1}(u))(\wp_{\beta+1, \gamma}(u)-\wp_{\beta, \gamma+1}(u)), 
\label{4e12}
\end{align}
where $\alpha > \beta >\gamma >\delta$.
If we put $\alpha=g+1$, we obtain $\langle \alpha \beta \gamma \delta \rangle=0$ by using $\wp_{g+1, *}(u)=0$.
Then it is sufficient
to consider $\langle \alpha \beta \gamma \delta \rangle$ in the range 
$g \ge \alpha > \beta > \gamma > \delta \ge 0$, that is, it is sufficient to calculate $_{g+1} C_4$ pieces of terms for the first term in the right-hand side of Eq.(\ref{4e1}).  
In $g=2$ case, we must calculate $_{3} C_4=0$ pieces of terms.  In $g=3$ case, we must calculate 
$_{4} C_4=1$ pieces of terms. In $g=4$ case, we must calculate $_{5} C_4=5$ pieces of terms. 

%%%%%%%%%%%%%%%%%%%%%%%%%%%%%%%%%%%%%%%%%%%%%%%%%%%%%%%%
\subsection{Demonstration to derive the genus four differential equation of 
hyperelliptic $\wp_{ij}(u)$ functions}

Here, we demonstrate to obtain differential equations for genus four hyperelliptic $\wp_{ij}(u)$ functions.\\
 {\bf {\underline{First term in the right-hand side of Eq.(\ref{4e1}) }}}\\
We first calculate the irregular term in the right-hand side of Eq.(\ref{4e1}) in the form
\begin{align}
&-2 H\sum \langle \alpha \beta \gamma \delta \rangle |\alpha \beta \gamma \delta|, 
\label{4e13}\\
&\text{where}
\nonumber\\
& \quad H=(2200)-(2110)+6(1111), 
\nonumber
\end{align}
There are five cases for non-zero $\langle \alpha \beta \gamma \delta \rangle$, which we 
denote 
$\Delta_2=-\langle 4 3 2 1 \rangle$, $\Delta_3=\langle 4 3 2 0 \rangle$, 
$\Delta_4=\langle 4 3 1 0 \rangle$, $\Delta_5=\langle 4 2 1 0 \rangle$,  
$\Delta_7=-\langle 3 2 1 0 \rangle$. They are explicitly given in the form 
\begin{align}
{\rm \bf a)}: \Delta_2=&-\langle 4 3 2 1 \rangle
=-\wp_{44}(u)\wp_{22}(u)+\wp_{42}(u)^2+\wp_{44}(u)\wp_{31}(u)-\wp_{43}(u)\wp_{41}(u)
\nonumber\\
&+\wp_{43}(u)\wp_{32}(u)-\wp_{42}\wp_{33}(u),  
\label{4e14}\\
{\rm \bf b)}: \Delta_3=&\langle 4 3 2 0 \rangle
=\wp_{44(u)}\wp_{21}(u)-\wp_{42}(u)\wp_{41}(u)-\wp_{43}(u)\wp_{31}(u)+\wp_{41}(u)\wp_{33},
\label{4e15}\\
{\rm \bf c)}: \Delta_4=&\langle 4 3 1 0 \rangle
=\wp_{44}(u)\wp_{11}(u)-\wp_{41}(u)^2-\wp_{42}(u)\wp_{31}(u)+\wp_{41}(u)\wp_{32}(u), 
\label{4e16}\\
{\rm \bf d)}: \Delta_5=&\langle 4 2 1 0 \rangle
=\wp_{43}(u)\wp_{11}(u)-\wp_{41}(u)\wp_{31}(u)-\wp_{42}(u)\wp_{21}(u)+\wp_{41}(u)\wp_{22}(u),  
\label{4e17}\\
{\rm \bf e)}: \Delta_7=&-\langle 3 2 1 0 \rangle
=-\wp_{33}(u)\wp_{11}(u)+\wp_{31}(u)^2+\wp_{42}(u)\wp_{11}(u)-\wp_{41}(u)\wp_{21}(u)
\nonumber\\
&+\wp_{32}(u)\wp_{21}(u)-\wp_{31}(u)\wp_{22}(u).
\label{4e18}
\end{align}
These irregular terms appear in the right-hand side of following differential equations.\\ 
\noindent
{\bf a)} $\Delta_2$ term appears in $-2H \langle 4321 \rangle |4321|=2 \Delta_2 H\cdot (1111)$ by using
$|4321|=(1111)$.
Then we obtain  
\begin{align}
&2 \Delta_2 H\cdot (1111)=2 \Delta_2 ((2200)-(2110)+6(1111))\cdot (1111) 
\nonumber\\
& =2 \Delta_2( (2200)\cdot (1111)-(2110)\cdot (1111)+6 (1111) \cdot (1111)) 
\nonumber\\
&=2 \Delta_2((3311)-(3221)+6(2222)), 
\label{4e19}
\end{align}
by using $(2200)\cdot (1111)=(3311)$, $(2110)\cdot (1111)=(3221)$, $(1111)\cdot (1111)=(2222)$.
Thus, in the right-hand side $[\wp_{4422}(u)]$, $[\wp_{4332}(u)]$, 
$[\wp_{3333}(u)]$, $\Delta_2$ appears
in the form 
\begin{align}
&{\rm \bf 8)}\ [\wp_{4422}]=2 \Delta_2+\cdots, 
\label{4e20}\\    
&{\rm \bf 12)}\ [\wp_{4332}]=-2 \Delta_2+\cdots, 
\label{4e21}\\  
& {\rm \bf 21)}\ [\wp_{3333}]=12 \Delta_2+\cdots.
\label{4e22}  
\end{align}
\noindent
{\bf b)} $\Delta_3$ term appears in $-2H \langle 4320 \rangle |4320|=-2 \Delta_3 H\cdot (1110)$ by using
$|4320|=(1110)$.
Then we obtain  
\begin{align}
&-2 \Delta_3 H\cdot (1110)=-2 \Delta_3 ((2200)-(2110)+6(1111))\cdot (1110)
\nonumber\\
& =-2 \Delta_3( (2200)\cdot (1110)-(2110)\cdot (1110)+6 (1111) \cdot (1110))
\nonumber\\
&=-2 \Delta_3((3310)-(3220)-(3211)+3(2221)) ,
\label{4e23}
\end{align}
by using $(2200)\cdot (1110)=(3310)+(3211)$, 
$(2110)\cdot (1110)=2(3211)+3(2221)+(3220)$, 
$(1111)\cdot (1110)=(2221)$.
Thus, in the right hand side $[\wp_{4421}(u)]$, $[\wp_{4331}(u)]$, 
$[\wp_{4322}(u)]$, $[\wp_{3332}(u)]$, $\Delta_3$ appears
in the form 
\begin{align}
&{\rm \bf 9)}\ [\wp_{4421}]=-2 \Delta_3+\cdots,
\label{4e24}\\    
&{\rm \bf 13)}\ [\wp_{4331}]=2 \Delta_3+\cdots,
\label{4e25}\\  
& {\rm \bf 14)}\ [\wp_{4322}]=2 \Delta_3+\cdots,
\label{4e26}\\ 
& {\rm \bf 22)}\ [\wp_{3332}]=-6 \Delta_3+\cdots.
\label{4e27} 
\end{align}
\noindent
{\bf c)} $\Delta_4$ term appears in $-2H \langle 4310 \rangle |4310|=-2 \Delta_4 H\cdot (1100)$ by using
$|4310|=(1100)$.
Then we obtain  
\begin{align}
&-2 \Delta_4 H\cdot (1100)=-2 \Delta_4 ((2200)-(2110)+6(1111))\cdot (1100)
\nonumber\\
& =-2 \Delta_4( (2200)\cdot (1100)-(2110)\cdot (1100)+6 (1111) \cdot (1100))
\nonumber\\
&=-2 \Delta_4((3300)-3(3111)-3(2220)+3(2211)), 
\label{4e28}
\end{align}
by using $(2200)\cdot (1100)=(3300)+(3210)+(2211)$, 
$(2110)\cdot (1100)=(3210)+3(2220)+4(2211)+3(3111)$, 
$(1111)\cdot (1100)=(2211)$.
Thus, in the right-hand side $[\wp_{4411}(u)]$, $[\wp_{4222}(u)]$, 
$[\wp_{3331}(u)]$, $[\wp_{3322}(u)]$, $\Delta_4$ appears
in the form 
\begin{align}
&{\rm \bf 10)}\ [\wp_{4411}]=-2 \Delta_4+\cdots,
\label{4e29}\\    
&{\rm \bf 17)}\ [\wp_{4222}]=6 \Delta_4+\cdots,
\label{4e30}\\  
& {\rm \bf 23)}\ [\wp_{3331}]=6 \Delta_4+\cdots,
\label{4e31}\\ 
& {\rm \bf 24)}\ [\wp_{3322}]=-6 \Delta_4+\cdots.
\label{4e32} 
\end{align}
\noindent
{\bf d)} $\Delta_5$ term appears in $-2H \langle 4210 \rangle |4210|=-2 \Delta_5 H\cdot (1000)$ by using
$|4210|=(1000)$.
Then we obtain  
\begin{align}
&-2 \Delta_5 H\cdot (1100)=-2 \Delta_4 ((2200)-(2110)+6(1111))\cdot (1000)
\nonumber\\
& =-2 \Delta_5( (2200)\cdot (1000)-(2110)\cdot (1000)+6 (1111) \cdot (1000))
\nonumber\\
&=-2 \Delta_5((3200)-(3110)-(2210)+3(2111)),
\label{4e33}
\end{align}
by using $(2200)\cdot (1000)=(3200)+(2210)$, 
$(2110)\cdot (1000)=(3110)+2(2210)+3(2111)$, 
$(1111)\cdot (1000)=(2111)$.
Thus, in the right-hand side $[\wp_{4311}(u)]$, $[\wp_{4221}(u)]$, 
$[\wp_{3321}(u)]$, $[\wp_{3222}(u)]$, $\Delta_5$ appears
in the form 
\begin{align}
&{\rm \bf 16)}\ [\wp_{4311}]=-2 \Delta_5+\cdots,
\label{4e34}\\    
&{\rm \bf 18)}\ [\wp_{4221}]=2 \Delta_5+\cdots,
\label{4e35}\\  
& {\rm \bf 25)}\ [\wp_{3321}]=2 \Delta_5+\cdots,
\label{4e36}\\ 
& {\rm \bf 27)}\ [\wp_{3222}]=-6 \Delta_5+\cdots.
\label{4e37} 
\end{align}
\noindent
{\bf e)} $\Delta_7$ term appears in $-2H \langle 3210 \rangle |3210|
=2 \Delta_7 H\cdot 1$ by using
$|3210|=1$.
Then we obtain  
\begin{align}
&2 \Delta_7 H\cdot 1=2 \Delta_7 ((2200)-(2110)+6(1111)) .
\label{4e38}
\end{align}
Thus, in the right-hand side $[\wp_{3311}(u)]$, $[\wp_{3211}(u)]$, 
$[\wp_{2222}(u)]$, $\Delta_7$ appears
in the form 
\begin{align}
&{\rm \bf 26)}\ [\wp_{3311}]=2 \Delta_7+\cdots,
\label{4e39}\\    
&{\rm \bf 28)}\ [\wp_{3221}]=-2 \Delta_7+\cdots,
\label{4e40}\\  
& {\rm \bf 31)}\ [\wp_{2222}]=12 \Delta_7+\cdots.
\label{4e41} 
\end{align}

Next, we give the list of expressions $|\alpha \beta \gamma \delta|$ with 
the combination of $(\mu \nu \rho \sigma )$.\\
\begin{align}
&{\rm 1)} |6543|=(3333), \ {\rm 2)} |6542|=(3332), \  {\rm 3)} |6541|=(3331)+(3322), 
\nonumber\\
& {\rm 4)} |6540|=(3330)+(3321)+(3222), \ {\rm 5)} |6532|=(3322), \ {\rm 6)} |6531|=(3321)+2(3222),
\nonumber\\ 
& {\rm 7)} |6530|=(3320)+(3311)+2(3221)+3(2222), 
\ {\rm 8)} |6521|=(3311)+(3221)+2(2222), 
\nonumber\\
& {\rm 9)} |6520|=(3310)+(3220)+2(3211)+3(2221),
\nonumber\\
&{\rm 10)} |6510|=(3300)+(3210)+(3111)+(2220)+2(2211),
\nonumber\\ 
& {\rm 11)} |6432|=(3222), \ {\rm 12)} |6431|=(3221)+3(2222), \ {\rm 13)} |6430|=(3220)+(3211)+3(2221),
\nonumber\\ 
& {\rm 14)} |6421|=(3211)+2(2221), \ {\rm 15)} |6420|=(3210)+2(3111)+2(2220)+4(2211), 
\nonumber\\
& {\rm 16)} |6410|=(3200)+(3110)+2(2210)+3(2111), \ {\rm 17)} |6321|=(3111)+(2211), 
\nonumber\\
& {\rm 18)} |6320|=(3110)+(2210)+3(2111), \  {\rm 19)} |6310|=(3100)+(2200)+2(2110)+3(1111), 
\nonumber\\
& {\rm 20)} |6210|=(3000)+(2100)+(1110), \ {\rm 21)} |5432|=(2222), 
\  {\rm 22)} |5431|=(2221), 
\nonumber\\
& {\rm 23)} |5430|=(2220)+(2211), \ {\rm 24)} |5421|=(2211), 
\ {\rm 25)} |5420|=(2210)+2(2111), 
\nonumber\\
& {\rm 26)} |5410|=(2200)+(2110)+2(1111), \ {\rm 27)} |5321|=(2111), 
\ {\rm 28)} |5320|=(2110)+3(1111),
\nonumber\\ 
& {\rm 29)} |5310|=(2100)+2(1110), \ {\rm 30)} |5210|=(2000)+(1100), \ {\rm 31)} |4321|=(1111), 
\nonumber\\
& {\rm 32)} |4320|=(1110), \ {\rm 33)} |4310|=(1100), \ {\rm 34)} |4210|=(1000), \ {\rm 35)} |3210|=1. 
\nonumber      
\end{align}
%%%%%%%%%%%%%%%%%%%%%%%%%%%%%%%%%%
We denote the left-hand side term, the first term in the right-hand side, 
the second term in the right-hand side, the third term in the right-hand side,  
in the form
\begin{align}
&{\rm (L.H.S.)}=\sum (\mu-1, \nu-1, \rho-1, \sigma-1) [\mu \nu \rho \sigma], 
\label{4e42}\\
&{\rm (F.R.H.S)}=-2 H\sum \langle \alpha \beta \gamma \delta \rangle |\alpha \beta \gamma \delta|, 
\label{4e43}\\
&{\rm (S.R.H.S)}=\frac{1}{2}\sum \langle \langle \alpha \beta \gamma \delta \rangle \rangle |\alpha \beta \gamma \delta|, 
\label{4e44}\\
&{\rm (T.R.H.S)}=-\frac{1}{8}\sum \left\{\alpha \beta \gamma \delta \right\} 
 |\alpha \beta \gamma \delta|.
\label{4e45}
\end{align}

Next, we calculate each differential equation step by step.\\
%%%%%%%%%%%%%%%%%%%%%%%%%%%%%%%
\vskip 3mm 
\noindent
{\bf  1)}\ coefficient of $(3333)$: \\
$|\alpha \beta \gamma \delta|$ which contains $(3333)$: $|6543|=\underline{(3333)}$, \\
(L.H.S.)=$[\wp_{4444}(u)]=\wp_{4444}(u)-6\wp_{44}(u)^2$, \\
(F.R.H.S.)=0, \\
(S.R.H.S.)=
$\displaystyle{\frac{1}{2} \langle \langle 6543 \rangle \rangle
=\lambda_9 \wp_{43}(u) +\lambda_8 \wp_{44}(u)}$, \\
(T.R.H.S.)=$\displaystyle{-\frac{1}{8} \left\{ 6543 \right\}
=\frac{1}{8} \lambda_7 \lambda_9}$, \\
(Diff. Eq.)\  
$\displaystyle{\wp_{4444}(u)-6\wp_{44}(u)^2=\lambda_9 \wp_{43}(u) +\lambda_8 \wp_{44}(u)
+\frac{1}{8} \lambda_7 \lambda_9}$. \\
%%%%%%%%%%%%%%%%%%%%%%%%%%%%%%
\vskip 3mm 
\noindent
{\bf  2)}\ coefficient of $(3332)$: \\
$|\alpha \beta \gamma \delta|$ which contains $(3332)$: $|6542|=\underline{(3332)}$, \\
(L.H.S.)=$[\wp_{4443}(u)]=\wp_{4443}(u)-6\wp_{44}(u) \wp_{43}(u)$, \\
(F.R.H.S.)=0, \\
(S.R.H.S.)=
$\displaystyle{\frac{1}{2} \langle \langle 6542 \rangle \rangle
=\frac{3}{2} \lambda_9 \wp_{42}(u) -\frac{1}{2}\lambda_9 \wp_{33}(u)
+\lambda_8 \wp_{43}(u)}$, \\
(T.R.H.S.)=$\displaystyle{-\frac{1}{8} \left\{ 6542 \right\}
=0}$, \\
(Diff. Eq.)\  
$\displaystyle{\wp_{4443}(u)-6\wp_{44}(u) \wp_{43}(u)=\frac{3}{2} \lambda_9 \wp_{42}(u) 
-\frac{1}{2}\lambda_9 \wp_{33}(u)+\lambda_8 \wp_{43}(u)}$. \\
%%%%%%%%%%%%%%%%%%%%%%%%%%%%%%%%%%%
%%%%%%%%%%%%%%%%%%%%%%%%%%%%%%
\vskip 3mm 
\noindent
{\bf  3)}\ coefficient of $(3331)$: \\
$|\alpha \beta \gamma \delta|$ which contains $(3331)$: $|6541|=\underline{(3331)}$+(3322), \\
(L.H.S.)=$[\wp_{4442}(u)]=\wp_{4442}(u)-6\wp_{44}(u) \wp_{42}(u)$, \\
(F.R.H.S.)=0, \\
(S.R.H.S.)=
$\displaystyle{\frac{1}{2} \langle \langle 6541 \rangle \rangle
=\frac{3}{2} \lambda_9 \wp_{41}(u) -\frac{1}{2}\lambda_9 \wp_{32}(u)
+\lambda_8 \wp_{42}(u)}$, \\
(T.R.H.S.)=$\displaystyle{-\frac{1}{8} \left\{ 6541 \right\}
=0}$, \\
(Diff. Eq.): 
$\displaystyle{\wp_{4442}(u)-6\wp_{44}(u) \wp_{42}(u)=\frac{3}{2} \lambda_9 \wp_{41}(u) 
-\frac{1}{2}\lambda_9 \wp_{32}(u)+\lambda_8 \wp_{42}(u)}$. \\
%%%%%%%%%%%%%%%%%%%%%%%%%%%%%%%%%%%
%%%%%%%%%%%%%%%%%%%%%%%%%%%%%%
\vskip 3mm 
\noindent
{\bf  4)}coefficient of $(3330)$: \\
$|\alpha \beta \gamma \delta|$ which contains $(3330)$: $|6540|=\underline{(3330)}$
+(3321)+(3222), \\
(L.H.S.)=$[\wp_{4441}(u)]=\wp_{4441}(u)-6\wp_{44}(u) \wp_{41}(u)$, \\
(F.R.H.S.)=0, \\
(S.R.H.S)=
$\displaystyle{\frac{1}{2} \langle \langle 6540 \rangle \rangle
=-\frac{1}{2} \lambda_9 \wp_{31}(u)+\lambda_8 \wp_{41}(u)}$, \\
(T.R.H.S)=$\displaystyle{-\frac{1}{8} \left\{ 6540 \right\} 
=0}$, \\
(Diff. Eq.): 
$\displaystyle{\wp_{4441}(u)-6\wp_{44}(u) \wp_{41}(u)=-\frac{1}{2} \lambda_9 \wp_{31}(u) 
+\lambda_8 \wp_{41}(u)}$. \\
%%%%%%%%%%%%%%%%%%%%%%%%%%%%%%%%%%%
%%%%%%%%%%%%%%%%%%%%%%%%%%%%%%
\vskip 3mm 
\noindent
{\bf  5)} coefficient of $(3322)$: \\
$|\alpha \beta \gamma \delta|$ which contains $(3322)$: $|6541|=(3331)+\underline{(3322)}$,
$|6532|=\underline{(3322)}$, \\
(L.H.S.)=$[\wp_{4433}(u)]=\wp_{4433}(u)-4\wp_{43}(u)^2-2\wp_{44}(u) \wp_{33}(u)$, \\
(F.R.H.S.)=0, \\
(S.R.H.S)=
$\displaystyle{\frac{1}{2} (\langle \langle 6541 \rangle \rangle+\langle \langle 6532 \rangle \rangle)
=\frac{3}{2} \lambda_9 \wp_{41}(u)-\frac{1}{2}\lambda_9\wp_{32}(u) +\lambda_8 \wp_{42}(u)
+\frac{1}{2}\lambda_7\wp_{43}(u)} $, \\
(T.R.H.S)=$\displaystyle{-\frac{1}{8}( \left\{ 6541 \right\}
+\left\{ 6532 \right\})=0}$, \\
(Diff. Eq.): $\displaystyle{\wp_{4433}(u)-4\wp_{43}(u)^2-2\wp_{44}(u) \wp_{33}(u)
=\frac{3}{2} \lambda_9 \wp_{41}(u)-\frac{1}{2}\lambda_9\wp_{32}(u) +\lambda_8 \wp_{42}(u)
+\frac{1}{2} \lambda_7 \wp_{43}(u) }$. \\
%%%%%%%%%%%%%%%%%%%%%%%%%%%%%%%%%%%
%%%%%%%%%%%%%%%%%%%%%%%%%%%%%%
\vskip 3mm 
\noindent
{\bf  6)} coefficient of $(3321)$: \\
$|\alpha \beta \gamma \delta|$ which contains $(3321)$: $|6540|=(3330)+\underline{(3321)}+(3222)$,
$|6531|=\underline{(3321)}+2(3222)$, \\
(L.H.S.)=$[\wp_{4432}(u)]=\wp_{4432}(u)-4\wp_{43}(u)\wp_{42}(u)-2\wp_{44}(u) \wp_{32}(u)$, \\
(F.R.H.S.)=0, \\
(S.R.H.S)=
$\displaystyle{\frac{1}{2} (\langle \langle 6540 \rangle \rangle+\langle \langle 6531 \rangle \rangle)
=-\frac{1}{2} \lambda_9 \wp_{31}(u)+\lambda_8\wp_{41}(u) +\frac{1}{2}\lambda_7\wp_{42}(u)} $, \\
(T.R.H.S)=$\displaystyle{-\frac{1}{8}( \left\{ 6540 \right\}
+\left\{ 6531 \right\})=0}$, \\
(Diff. Eq.): $\displaystyle{\wp_{4432}(u)-4\wp_{43}(u)\wp_{42}(u)-2\wp_{44}(u) \wp_{32}(u)
=-\frac{1}{2} \lambda_9 \wp_{31}(u)+\lambda_8\wp_{41}(u) +\frac{1}{2}\lambda_7\wp_{42}(u)} $. \\
%%%%%%%%%%%%%%%%%%%%%%%%%%%%%%%%%%%
%%%%%%%%%%%%%%%%%%%%%%%%%%%%%%
\vskip 3mm 
\noindent
{\bf  7)} coefficient of $(3320)$: \\
$|\alpha \beta \gamma \delta|$ which contains $(3320)$: $|6530|=\underline{(3320)}+(3310)
+2(3221)+3(2222)$, \\
(L.H.S.)=$[\wp_{4431}(u)]=\wp_{4431}(u)-4\wp_{43}(u)\wp_{41}(u)-2\wp_{44}(u) \wp_{31}(u)$, \\
(F.R.H.S.)=0, \\
(S.R.H.S)=
$\displaystyle{\frac{1}{2} \langle \langle 6530 \rangle \rangle
=\frac{1}{2} \lambda_7 \wp_{41}(u)}$, \\
(T.R.H.S)=$\displaystyle{-\frac{1}{8}\left\{ 6530 \right\}=0}$, \\
(Diff. Eq.): $\displaystyle{\wp_{4431}(u)-4\wp_{43}(u)\wp_{41}(u)-2\wp_{44}(u) \wp_{31}(u)
=\frac{1}{2} \lambda_7 \wp_{41}(u)}$. \\
%%%%%%%%%%%%%%%%%%%%%%%%%%%%%%%%%%%
%%%%%%%%%%%%%%%%%%%%%%%%%%%%%%
\vskip 3mm 
\noindent
{\bf  8)} coefficient of $(3311)$: \\
$|\alpha \beta \gamma \delta|$ which contains $(3311)$: $|6530|=(3320)+\underline{(3311)}
+2(3221)+3((2222)$,
$|6521|=\underline{(3311)}+(3221)+2(2222)$, \\
(L.H.S.)=$[\wp_{4422}(u)]=\wp_{4422}(u)-4\wp_{42}(u)^2-2\wp_{44}(u) \wp_{22}(u)$, \\
(F.R.H.S.)=$2 \Delta_2$, \\
(S.R.H.S)=
$\displaystyle{\frac{1}{2} (\langle \langle 6530 \rangle \rangle+\langle \langle 6521 \rangle \rangle)
=\frac{1}{2} \lambda_7 \wp_{41}(u)} $, \\
(T.R.H.S)=$\displaystyle{-\frac{1}{8}( \left\{ 6530 \right\}
+\left\{ 6521 \right\})=0}$, \\
(Diff. Eq.): $\displaystyle{\wp_{4422}(u)-4\wp_{42}(u)^2-2\wp_{44}(u) \wp_{22}(u)
=2 \Delta_2+\frac{1}{2} \lambda_7 \wp_{41}(u)} $. \\
%%%%%%%%%%%%%%%%%%%%%%%%%%%%%%%%%%%
%%%%%%%%%%%%%%%%%%%%%%%%%%%%%%
\vskip 3mm 
\noindent
{\bf  9)} coefficient of $(3310)$: \\
$|\alpha \beta \gamma \delta|$ which contains $(3310)$: $|6520|=\underline{(3310)}+(3220)
+2(3211)+3(2221)$, \\
(L.H.S.)=$[\wp_{4421}(u)]=\wp_{4421}(u)-4\wp_{42}(u)\wp_{41}(u)-2\wp_{44}(u) \wp_{21}(u)$, \\
(F.R.H.S.)=$-2 \Delta_3$, \\
(S.R.H.S)=
$\displaystyle{\frac{1}{2} \langle \langle 6520 \rangle \rangle =0}$, \\
(T.R.H.S)=$\displaystyle{-\frac{1}{8}\left\{ 6520 \right\}=0}$, \\
(Diff. Eq.): $\displaystyle{\wp_{4421}(u)-4\wp_{42}(u)\wp_{41}(u)-2\wp_{44}(u) \wp_{21}(u)
=-2 \Delta_3}$. \\
%%%%%%%%%%%%%%%%%%%%%%%%%%%%%%%%%%%
%%%%%%%%%%%%%%%%%%%%%%%%%%%%%%
\vskip 3mm 
\noindent
{\bf  10)} coefficient of $(3300)$: \\
$|\alpha \beta \gamma \delta|$ which contains $(3300)$: $|6510|=\underline{(3300)}+(3210)
+(3111)+(2220)+2(2211)$, \\
(L.H.S.)=$[\wp_{4411}(u)]=\wp_{4411}(u)-4\wp_{41}(u)^2-2\wp_{44}(u) \wp_{11}(u)$, \\
(F.R.H.S.)=$-2 \Delta_4$, \\
(S.R.H.S)=
$\displaystyle{\frac{1}{2} \langle \langle 6510 \rangle \rangle =0}$, \\
(T.R.H.S)=$\displaystyle{-\frac{1}{8}\left\{ 6510 \right\}=0}$, \\
(Diff. Eq.): $\displaystyle{\wp_{4411}(u)-4\wp_{41}(u)^2-2\wp_{44}(u) \wp_{11}(u)
=-2 \Delta_4}$. \\
%%%%%%%%%%%%%%%%%%%%%%%%%%%%%%%%%%%
%%%%%%%%%%%%%%%%%%%%%%%%%%%%%%
\vskip 3mm 
\noindent
{\bf  11)} coefficient of $(3222)$: \\
$|\alpha \beta \gamma \delta|$ which contains $(3222)$: $|6540|=(3330)+(3321)+\underline{(3222)}$,
$|6531|=(3321)+\underline{2(3222)}$, $|6432|=\underline{(3222)}$, \\
(L.H.S.)=$[\wp_{4333}(u)]=\wp_{4333}(u)-6\wp_{43}(u)\wp_{33}(u)$, \\
(F.R.H.S.)=0, \\
(S.R.H.S)=
$\displaystyle{\frac{1}{2} (\langle \langle 6540 \rangle \rangle+2\langle \langle 6531 \rangle \rangle
+\langle \langle 6432 \rangle \rangle)}$\\
$\displaystyle{=\frac{3}{2} \lambda_9 \wp_{31}(u)-\frac{3}{2} \lambda_9 \wp_{22}(u)+\lambda_8 \wp_{41}(u)
+\lambda_7 \wp_{42}(u) +\lambda_6 \wp_{43}(u)-\frac{1}{2} \lambda_5 \wp_{44}(u) } $, \\
(T.R.H.S)=$\displaystyle{-\frac{1}{8}( \left\{ 6540 \right\}+2\left\{ 6531 \right\}
+\left\{ 6432 \right\} )=-\frac{1}{4} \lambda_4 \lambda_9}$, \\
(Diff. Eq.): $\displaystyle{\wp_{4333}(u)-6\wp_{43}(u)\wp_{33}(u)
=\frac{3}{2} \lambda_9 \wp_{31}(u)-\frac{3}{2} \lambda_9 \wp_{22}(u)+\lambda_8 \wp_{41}(u)
+\lambda_7 \wp_{42}(u)}$\\
$\displaystyle{+\lambda_6 \wp_{43}(u)-\frac{1}{2} \lambda_5 \wp_{44}(u) 
-\frac{1}{4} \lambda_4 \lambda_9 } $. \\
%%%%%%%%%%%%%%%%%%%%%%%%%%%%%%%%%%%
%%%%%%%%%%%%%%%%%%%%%%%%%%%%%%
\vskip 3mm 
\noindent
{\bf  12)} coefficient of $(3221)$: \\
$|\alpha \beta \gamma \delta|$ which contains $(3221)$: $|6530|=(3320)+(3311)+\underline{2(3221)}+3(2222)$,
$|6521|=(3311)+\underline{(3221)}+2(2222)$, $|6431|=\underline{(3221)}+3(2222)$, \\
(L.H.S.)=$[\wp_{4332}(u)]=\wp_{4332}(u)-4\wp_{43}(u)\wp_{32}(u)-2\wp_{42}(u)\wp_{33}(u)$, \\
(F.R.H.S.)=$-2 \Delta_2$, \\
(S.R.H.S)=
$\displaystyle{\frac{1}{2} (2 \langle \langle 6530 \rangle \rangle+\langle \langle 6521\rangle \rangle
+\langle \langle 6431 \rangle \rangle)}$
$\displaystyle{=-\lambda_9 \wp_{21}(u)+\lambda_7 \wp_{41}(u)+\lambda_6 \wp_{42}(u)  } $, \\
(T.R.H.S)=$\displaystyle{-\frac{1}{8}(2 \left\{ 6530 \right\}+\left\{ 6521 \right\}
+\left\{ 6431 \right\}) =-\frac{1}{8} \lambda_3 \lambda_9}$, \\
(Diff. Eq.): $\displaystyle{\wp_{4332}(u)-4\wp_{43}(u)\wp_{32}(u)-2\wp_{42}(u)\wp_{33}(u)
=-2 \Delta_2 -\lambda_9 \wp_{21}(u)+\lambda_7 \wp_{41}(u)}$\\
$\displaystyle{+\lambda_6 \wp_{42}(u) 
-\frac{1}{8} \lambda_3 \lambda_9 } $. \\
%%%%%%%%%%%%%%%%%%%%%%%%%%%%%%%%%%%
%%%%%%%%%%%%%%%%%%%%%%%%%%%%%%
\vskip 3mm 
\noindent
{\bf  13)} coefficient of $(3220)$: \\
$|\alpha \beta \gamma \delta|$ which contains $(3220)$: $|6520|=(3310)+\underline{(3220)}
+2(3211)+3(2221)$,
$|6430|=\underline{(3220)}+(3211)+3(2221)$, \\
(L.H.S.)=$[\wp_{4331}(u)]=\wp_{4331}(u)-4\wp_{43}(u)\wp_{31}(u)-2\wp_{41}(u)\wp_{33}(u)$, \\
(F.R.H.S.)=$2 \Delta_3$, \\
(S.R.H.S)=
$\displaystyle{\frac{1}{2} (\langle \langle 6520 \rangle \rangle+\langle \langle 6430\rangle \rangle)}$
$\displaystyle{=\frac{1}{2}\lambda_9 \wp_{11}(u)+\lambda_6 \wp_{41}(u) } $, \\
(T.R.H.S)=$\displaystyle{-\frac{1}{8}(\left\{ 6520 \right\}+\left\{ 6430 \right\}) =0}$, \\
(Diff. Eq.): $\displaystyle{\wp_{4331}(u)-4\wp_{43}(u)\wp_{31}(u)-2\wp_{41}(u)\wp_{33}(u)
=2 \Delta_3 +\frac{1}{2}\lambda_9 \wp_{11}(u)+\lambda_6 \wp_{41}(u) } $.  \\
%%%%%%%%%%%%%%%%%%%%%%%%%%%%%%%%%%%
%%%%%%%%%%%%%%%%%%%%%%%%%%%%%%
\vskip 3mm 
\noindent
{\bf  14)} coefficient of $(3211)$: \\
$|\alpha \beta \gamma \delta|$ which contains $(3211)$: $|6520|=(3310)+(3220)+\underline{2(3211)}+3(2221)$,
$|6421|=\underline{(3211)}+2(2221)$, $|6430|=(3220)+\underline{(3211)}+3(2221)$, \\
(L.H.S.)=$[\wp_{4322}(u)]=\wp_{4322}(u)-4\wp_{42}(u)\wp_{32}(u)-2\wp_{43}(u)\wp_{22}(u)$, \\
(F.R.H.S.)=$2 \Delta_3$, \\
(S.R.H.S)=
$\displaystyle{\frac{1}{2} (2 \langle \langle 6520 \rangle \rangle+\langle \langle 6421\rangle \rangle
+\langle \langle 6430 \rangle \rangle)}$
$\displaystyle{=-\lambda_9 \wp_{11}(u)+\lambda_6 \wp_{41}(u)+\frac{1}{2}\lambda_5 \wp_{42}(u)  } $, \\
(T.R.H.S)=$\displaystyle{-\frac{1}{8}(2 \left\{ 6520 \right\}+\left\{ 6421 \right\}
+\left\{ 6430 \right\}) =-\frac{1}{4} \lambda_2 \lambda_9}$, \\
(Diff. Eq.): $\displaystyle{\wp_{4322}(u)-4\wp_{42}(u)\wp_{32}(u)-2\wp_{43}(u)\wp_{22}(u)
=2 \Delta_3 -\lambda_9 \wp_{11}(u)+\lambda_6 \wp_{41}(u)}$\\
$\displaystyle{+\frac{1}{2}\lambda_5 \wp_{42}(u) 
-\frac{1}{4} \lambda_2 \lambda_9 } $. \\
%%%%%%%%%%%%%%%%%%%%%%%%%%%%%%%%%%%
%%%%%%%%%%%%%%%%%%%%%%%%%%%%%%
\vskip 3mm 
\noindent
{\bf  15)} coefficient of $(3210)$: \\
$|\alpha \beta \gamma \delta|$ which contains $(3210)$: $|6510|=(3300)+\underline{(3210)}
+(3111)+(2220)+2(2211)$,
$|6420|=\underline{(3210)}+2(3111)+2(2220)+4(2211)$, \\
(L.H.S.)=$[\wp_{4321}(u)]=\wp_{4321}(u)-2\wp_{43}(u)\wp_{21}(u)-2\wp_{42}(u)\wp_{31}(u)
-2\wp_{41}(u)\wp_{32}(u)$, \\
(F.R.H.S.)=0, \\
(S.R.H.S)=
$\displaystyle{\frac{1}{2} (\langle \langle 6510 \rangle \rangle+\langle \langle 6420\rangle \rangle)
=\frac{1}{2}\lambda_5 \wp_{41}(u)  } $, \\
(T.R.H.S)=$\displaystyle{-\frac{1}{8}(\left\{ 6510 \right\}+\left\{ 6420 \right\}) 
=-\frac{1}{8} \lambda_1 \lambda_9}$, \\
(Diff. Eq.): $\displaystyle{\wp_{4321}(u)-2\wp_{43}(u)\wp_{21}(u)-2\wp_{42}(u)\wp_{31}(u)
-2\wp_{41}(u)\wp_{32}(u)=\frac{1}{2}\lambda_5 \wp_{41}(u) -\frac{1}{8} \lambda_1 \lambda_9 } $. \\
%%%%%%%%%%%%%%%%%%%%%%%%%%%%%%%%%%%
%%%%%%%%%%%%%%%%%%%%%%%%%%%%%%
\vskip 3mm 
\noindent
{\bf  16)} coefficient of $(3200)$: \\
$|\alpha \beta \gamma \delta|$ which contains $(3200)$: $|6410|=\underline{(3200)}
+(3110)+2(2210)+3(2111)$, \\
(L.H.S.)=$[\wp_{4311}(u)]=\wp_{4311}(u)-4\wp_{41}(u)\wp_{31}(u)-2\wp_{43}(u)\wp_{11}(u)$, \\
(F.R.H.S.)=$-2 \Delta_5$, \\
(S.R.H.S)=
$\displaystyle{\frac{1}{2} \langle \langle 6410 \rangle \rangle=0 } $, \\
(T.R.H.S)=$\displaystyle{-\frac{1}{8}\left\{ 6410 \right\} 
=-\frac{1}{4} \lambda_0 \lambda_9}$, \\
(Diff. Eq.): $\displaystyle{\wp_{4311}(u)-4\wp_{41}(u)\wp_{31}(u)-2\wp_{43}(u)\wp_{11}(u)
=-2 \Delta_5-\frac{1}{4} \lambda_0 \lambda_9 } $. \\
%%%%%%%%%%%%%%%%%%%%%%%%%%%%%%%%%%%
%%%%%%%%%%%%%%%%%%%%%%%%%%%%%%
\vskip 3mm 
\noindent
{\bf  17)} coefficient of $(3111)$: \\
$|\alpha \beta \gamma \delta|$ which contains $(3111)$: $|6510|=(3300)+(3210)+\underline{(3111)}
+(2222)+2(2211)$,
$|6420|=(3210)+\underline{2(3111)}+2(2220)+4(2211)$, $|6321|=\underline{(3111)}+(2211)$, \\
(L.H.S.)=$[\wp_{4222}(u)]=\wp_{4222}(u)-6\wp_{42}(u)\wp_{22}(u)$, \\
(F.R.H.S.)=$6 \Delta_4$, \\
(S.R.H.S)=
$\displaystyle{\frac{1}{2} (\langle \langle 6510 \rangle \rangle+2 \langle \langle 6420\rangle \rangle
+\langle \langle 6321 \rangle \rangle)}$
$\displaystyle{=\lambda_5 \wp_{41}(u)+\lambda_4 \wp_{42}(u)-\frac{1}{2}\lambda_3 \wp_{43}(u)
+\lambda_2\wp_{44}(u)  } $, \\
(T.R.H.S)=$\displaystyle{-\frac{1}{8}(\left\{ 6510 \right\}+2 \left\{ 6420 \right\}
+\left\{ 6421 \right\}) =-\frac{1}{4} \lambda_1 \lambda_9}$, \\
(Diff. Eq.): $\displaystyle{\wp_{4222}(u)-6\wp_{42}(u)\wp_{22}(u)}$\\
$\displaystyle{=6 \Delta_4 +\lambda_5 \wp_{41}(u)+\lambda_4 \wp_{42}(u)-\frac{1}{2}\lambda_3 \wp_{43}(u)
+\lambda_2 \wp_{44}(u) 
-\frac{1}{4} \lambda_1 \lambda_9 } $. \\
%%%%%%%%%%%%%%%%%%%%%%%%%%%%%%%%%%%
%%%%%%%%%%%%%%%%%%%%%%%%%%%%%%
\vskip 3mm 
\noindent
{\bf  18)} coefficient of $(3110)$: \\
$|\alpha \beta \gamma \delta|$ which contains $(3110)$: $|6410|=(3200)+\underline{(3110)}
+2(2210)+3(2111)$,
$|6320|=\underline{(3110)}+(2210)+3(2111)$, \\
(L.H.S.)=$[\wp_{4221}(u)]=\wp_{4221}(u)-4\wp_{42}(u)\wp_{21}(u)-2\wp_{41}(u)\wp_{22}(u)$, \\
(F.R.H.S.)=$2 \Delta_5$, \\
(S.R.H.S)=
$\displaystyle{\frac{1}{2} (\langle \langle 6410 \rangle \rangle+ \langle \langle 6320\rangle \rangle)
=\lambda_4 \wp_{41}(u)+\frac{1}{2}\lambda_1 \wp_{44}(u)  } $, \\
(T.R.H.S)=$\displaystyle{-\frac{1}{8}(\left\{ 6410 \right\}+\left\{ 6320 \right\} ) 
=-\frac{1}{4} \lambda_0 \lambda_9}$, \\
(Diff. Eq.): $\displaystyle{\wp_{4221}(u)-4\wp_{42}(u)\wp_{21}(u)-2\wp_{41}(u)\wp_{22}(u)  }$
$\displaystyle{=2 \Delta_5 +\lambda_4 \wp_{41}(u)+\frac{1}{2}\lambda_1 \wp_{44}(u)
-\frac{1}{4} \lambda_0 \lambda_9 } $. \\
%%%%%%%%%%%%%%%%%%%%%%%%%%%%%%%%%%%
%%%%%%%%%%%%%%%%%%%%%%%%%%%%%%
\vskip 3mm 
\noindent
{\bf  19)} coefficient of $(3100)$: \\
$|\alpha \beta \gamma \delta|$ which contains $(3100)$: $|6310|=\underline{(3100)}
+(2200)+2(2110)+3(1111)$, \\
(L.H.S.)=$[\wp_{4211}(u)]=\wp_{4211}(u)-4\wp_{41}(u)\wp_{21}(u)-2\wp_{42}(u)\wp_{11}(u)$, \\
(F.R.H.S.)=0, \\
(S.R.H.S)=
$\displaystyle{\frac{1}{2} \langle \langle 6310 \rangle \rangle
=\frac{1}{2}\lambda_3 \wp_{41}(u)+\lambda_0 \wp_{44}(u) } $, \\
(T.R.H.S)=$\displaystyle{-\frac{1}{8}\left\{ 6310 \right\}=0 }$, \\
(Diff. Eq.): $\displaystyle{\wp_{4211}(u)-4\wp_{41}(u)\wp_{21}(u)-2\wp_{42}(u)\wp_{11}(u)
= \frac{1}{2}\lambda_3 \wp_{41}(u)+\lambda_0 \wp_{44}(u) } $. \\
%%%%%%%%%%%%%%%%%%%%%%%%%%%%%%%%%%%
%%%%%%%%%%%%%%%%%%%%%%%%%%%%%%
\vskip 3mm 
\noindent
{\bf  20)} coefficient of $(3000)$: \\
$|\alpha \beta \gamma \delta|$ which contains $(3000)$: $|6210|=\underline{(3000)}
+(2100)+(1110)$, \\
(L.H.S.)=$[\wp_{4111}(u)]=\wp_{4111}(u)-6\wp_{41}(u)\wp_{11}(u)$, \\
(F.R.H.S.)=0, \\
(S.R.H.S)=
$\displaystyle{\frac{1}{2} \langle \langle 6210 \rangle \rangle
=\lambda_2 \wp_{41}(u)-\frac{1}{2}\lambda_1 \wp_{42}(u)+\lambda_0 \wp_{43}(u) } $, \\
(T.R.H.S)=$\displaystyle{-\frac{1}{8}\left\{ 6210 \right\}=0 }$, \\
(Diff. Eq.): $\displaystyle{\wp_{4111}(u)-6\wp_{41}(u)\wp_{11}(u)
=\lambda_2 \wp_{41}(u)-\frac{1}{2}\lambda_1 \wp_{42}(u)+\lambda_0 \wp_{43}(u) } $. \\
%%%%%%%%%%%%%%%%%%%%%%%%%%%%%%%%%%%
%%%%%%%%%%%%%%%%%%%%%%%%%%%%%%
\vskip 3mm 
\noindent
{\bf  21)} coefficient of $(2222)$: \\
$|\alpha \beta \gamma \delta|$ which contains $(2222)$: $|6530|=(3320)+(3311)
+2(3221)+\underline{3(2222)}$,
$|6521|=(3311)+(3221)+\underline{2(2222)}$, $|6431|=(3221)+\underline{3(2222)}$,
$|5432|=\underline{(2222)}$, \\
(L.H.S.)=$[\wp_{3333}(u)]=\wp_{3333}(u)-6\wp_{33}(u)^2$, \\
(F.R.H.S.)=$12 \Delta_2$, \\
(S.R.H.S)=
$\displaystyle{\frac{1}{2} (3 \langle \langle 6530 \rangle \rangle+2 \langle \langle 6521\rangle \rangle
+3 \langle \langle 6431 \rangle \rangle+\langle \langle 5432 \rangle \rangle)}$\\
$\displaystyle{=-3\lambda_9 \wp_{21}(u)+4 \lambda_8 \wp_{31}(u)-3\lambda_8 \wp_{22}(u)  
+\lambda_7 \wp_{32}(u)+\lambda_6 \wp_{33}(u)+\lambda_5 \wp_{43}(u)-3\lambda_4 \wp_{44}(u)}$, \\
(T.R.H.S)=$\displaystyle{-\frac{1}{8}(3 \left\{ 6530 \right\}+2 \left\{ 6521 \right\}
+3 \left\{ 6431 \right\}+\left\{ 5432 \right\}) =\frac{1}{8} \lambda_5 \lambda_7
-\frac{1}{2} \lambda_4 \lambda_8-\frac{3}{8} \lambda_3 \lambda_9}$, \\
(Diff. Eq.): $\displaystyle{\wp_{3333}(u)-6\wp_{33}(u)^2=12 \Delta_2
-3\lambda_9 \wp_{21}(u)+4 \lambda_8 \wp_{31}(u)-3\lambda_8 \wp_{22}(u)  
+\lambda_7 \wp_{32}(u)}$\\
$\displaystyle{+\lambda_6 \wp_{33}(u)+\lambda_5 \wp_{43}(u)-3\lambda_4 \wp_{44}(u)
+\frac{1}{8} \lambda_5 \lambda_7
-\frac{1}{2} \lambda_4 \lambda_8-\frac{3}{8} \lambda_3 \lambda_9} $. \\
%%%%%%%%%%%%%%%%%%%%%%%%%%%%%%%%%%%
%%%%%%%%%%%%%%%%%%%%%%%%%%%%%%
\vskip 3mm 
\noindent
{\bf  22)} coefficient of $(2221)$: \\
$|\alpha \beta \gamma \delta|$ which contains $(2221)$: $|6520|=(3310)+(3220)
+2(3211)+\underline{3(2221)}$,
$|6430|=(3220)+(3211)+\underline{3(2221)}$, $|6421|=(3211)+\underline{2(2221)}$,
$|5431|=\underline{(2221)}$, \\
(L.H.S.)=$[\wp_{3332}(u)]=\wp_{3332}(u)-6\wp_{33}(u)\wp_{32}(u)$, \\
(F.R.H.S.)=$-6 \Delta_3$, \\
(S.R.H.S)=
$\displaystyle{\frac{1}{2} (3 \langle \langle 6520 \rangle \rangle+3 \langle \langle 6430 \rangle \rangle
+2 \langle \langle 6421 \rangle \rangle+\langle \langle 5431 \rangle \rangle)}$\\
$\displaystyle{=-\frac{3}{2}\lambda_9 \wp_{11}(u)-2 \lambda_8 \wp_{21}(u)+\frac{3}{2}\lambda_7 \wp_{31}(u)  
-\frac{1}{2}\lambda_7 \wp_{22}(u)+\lambda_6 \wp_{32}(u)+\lambda_5 \wp_{42}(u)
-\frac{3}{2}\lambda_3 \wp_{44}(u)}$, \\
(T.R.H.S)=$\displaystyle{-\frac{1}{8}(3 \left\{ 6520 \right\}+3 \left\{ 6430 \right\}
+2 \left\{ 6421 \right\}+\left\{ 5431 \right\}) =-\frac{1}{2} \lambda_2 \lambda_9
-\frac{1}{4} \lambda_3 \lambda_8}$, \\
(Diff. Eq.): $\displaystyle{\wp_{3332}(u)-6\wp_{33}(u)\wp_{32}(u)=-6 \Delta_3
-\frac{3}{2}\lambda_9 \wp_{11}(u)-2 \lambda_8 \wp_{21}(u)+\frac{3}{2}\lambda_7 \wp_{31}(u) }$\\
$\displaystyle{ 
-\frac{1}{2}\lambda_7 \wp_{22}(u)+\lambda_6 \wp_{32}(u)+\lambda_5 \wp_{42}(u)
-\frac{3}{2}\lambda_3 \wp_{44}(u)
-\frac{1}{2} \lambda_2 \lambda_9-\frac{1}{4} \lambda_3 \lambda_8} $. \\
%%%%%%%%%%%%%%%%%%%%%%%%%%%%%%%%%%%
%%%%%%%%%%%%%%%%%%%%%%%%%%%%%%
\vskip 3mm 
\noindent
{\bf  23)} coefficient of $(2220)$: \\
$|\alpha \beta \gamma \delta|$ which contains $(2220)$: $|6510|=(3300)+(3210)
+(3111)+\underline{(2220)}+2(2211)$,
$|6420|=(3210)+2(3111)+\underline{2(2220)}+4(2211)$, $|5430|=\underline{(2220)}+(2211)$, \\
(L.H.S.)=$[\wp_{3331}(u)]=\wp_{3331}(u)-6\wp_{33}(u)\wp_{31}(u)$, \\
(F.R.H.S.)=$6 \Delta_4$, \\
(S.R.H.S)=
$\displaystyle{\frac{1}{2} (\langle \langle 6510 \rangle \rangle+2 \langle \langle 6420 \rangle \rangle
+\langle \langle 5430 \rangle \rangle)}$
$\displaystyle{=\lambda_8 \wp_{11}(u)-\frac{1}{2} \lambda_7 \wp_{21}(u)+\lambda_6 \wp_{31}(u)  
+\lambda_5 \wp_{41}(u)}$, \\
(T.R.H.S)=$\displaystyle{-\frac{1}{8}(\left\{ 6510 \right\}+2 \left\{ 6420 \right\}
+\left\{ 5430 \right\}) =-\frac{1}{4} \lambda_1 \lambda_9}$, \\
(Diff. Eq.): $\displaystyle{\wp_{3331}(u)-6\wp_{33}(u)\wp_{31}(u)=6 \Delta_4
+\lambda_8 \wp_{11}(u)-\frac{1}{2} \lambda_7 \wp_{21}(u)+\lambda_6 \wp_{31}(u) }$\\ 
$\displaystyle{+\lambda_5 \wp_{41}(u)
-\frac{1}{4} \lambda_1 \lambda_9}$. \\
%%%%%%%%%%%%%%%%%%%%%%%%%%%%%%%%%%%
%%%%%%%%%%%%%%%%%%%%%%%%%%%%%%
\vskip 3mm 
\noindent
{\bf  24)} coefficient of $(2211)$: \\
$|\alpha \beta \gamma \delta|$ which contains $(2211)$: $|6510|=(3300)+(3210)
+(3111)+(2220)+\underline{2(2211)}$,
$|6420|=(3210)+2(3111)+2(2220)+\underline{4(2211)}$, 
$|6321|=(3100)+\underline{(2211)}$, 
$|5430|=(2220)+\underline{(2211)}$, $|5421|=\underline{(2211)}$, \\
(L.H.S.)=$[\wp_{3322}(u)]=\wp_{3322}(u)-4\wp_{32}(u)^2-2\wp_{33}(u)\wp_{22}(u)$, \\
(F.R.H.S.)=$-6 \Delta_4$, \\
(S.R.H.S)=
$\displaystyle{\frac{1}{2} (2\langle \langle 6510 \rangle \rangle+4 \langle \langle 6420 \rangle \rangle
+\langle \langle 6321 \rangle \rangle+\langle \langle 5430 \rangle \rangle
+\langle \langle 5421 \rangle \rangle)}$\\
$\displaystyle{=-2\lambda_8 \wp_{11}(u)-\frac{1}{2} \lambda_7 \wp_{21}(u)+\lambda_6 \wp_{31}(u)  
+\frac{1}{2}\lambda_5 \wp_{41}(u)+\frac{1}{2}\lambda_5 \wp_{32}(u)+\lambda_4 \wp_{42}(u)}$\\
$\displaystyle{-\frac{1}{2} \lambda_3 \wp_{43}(u)-2\lambda_2 \wp_{44}(u)}$, \\
(T.R.H.S)=$\displaystyle{-\frac{1}{8}(2 \left\{ 6510 \right\}+4 \left\{ 6420 \right\}
+\left\{ 6321 \right\}+\left\{ 5430 \right\}+\left\{ 5421 \right\}) 
=-\frac{1}{2} \lambda_1 \lambda_9-\frac{1}{2} \lambda_2 \lambda_8}$, \\
(Diff. Eq.): $\displaystyle{\wp_{3322}(u)-4\wp_{32}(u)^2-2  \wp_{33}(u) \wp_{22}(u)
=-6 \Delta_4
-2\lambda_8 \wp_{11}(u)-\frac{1}{2} \lambda_7 \wp_{21}(u)}$\\
$\displaystyle{+\lambda_6 \wp_{31}(u)  
+\frac{1}{2}\lambda_5 \wp_{41}(u)+\frac{1}{2}\lambda_5 \wp_{32}(u)+\lambda_4 \wp_{42}(u)
-\frac{1}{2} \lambda_3 \wp_{43}(u)-2\lambda_2 \wp_{44}(u)}$\\
$\displaystyle{-\frac{1}{2} \lambda_1 \lambda_9-\frac{1}{2} \lambda_2 \lambda_8}$. \\
%%%%%%%%%%%%%%%%%%%%%%%%%%%%%%%%%%%
%%%%%%%%%%%%%%%%%%%%%%%%%%%%%%
\vskip 3mm 
\noindent
{\bf  25)} coefficient of $(2210)$: \\
$|\alpha \beta \gamma \delta|$ which contains $(2210)$: $|6320|=(3110)+\underline{(2210)}
+3(2111)$, $|5420|=\underline{(2210)}+2(2111)$, $|6410|=(3200)+(3110)+\underline{2(2210)}+3(2111)$, \\
(L.H.S.)=$[\wp_{3321}(u)]=\wp_{3321}(u)-4\wp_{32}(u)\wp_{31}(u)-2\wp_{33}(u)\wp_{21}(u)$, \\
(F.R.H.S.)=$2 \Delta_5$, \\
(S.R.H.S)=
$\displaystyle{\frac{1}{2} (\langle \langle 6320 \rangle \rangle+\langle \langle 5420 \rangle \rangle
+2\langle \langle 6410 \rangle \rangle)}$
$\displaystyle{=\frac{1}{2}\lambda_5 \wp_{31}(u)+\lambda_4 \wp_{41}(u)-\lambda_1 \wp_{44}(u) }$, \\
(T.R.H.S)=$\displaystyle{-\frac{1}{8}( \left\{ 6320 \right\}+\left\{ 5420 \right\}
+2\left\{ 6410 \right\}) =-\frac{1}{4} \lambda_1 \lambda_8-\frac{1}{2} \lambda_0 \lambda_9}$, \\
(Diff. Eq.): $\displaystyle{\wp_{3321}(u)-4\wp_{32}(u)\wp_{31}(u)-2\wp_{33}(u)\wp_{21}(u)
=2 \Delta_5
+\frac{1}{2}\lambda_5 \wp_{31}(u)+\lambda_4 \wp_{41}(u)-\lambda_1 \wp_{44}(u)}$\\
$\displaystyle{-\frac{1}{4} \lambda_1 \lambda_8-\frac{1}{2} \lambda_0 \lambda_9}$. \\
%%%%%%%%%%%%%%%%%%%%%%%%%%%%%%%%%%%
%%%%%%%%%%%%%%%%%%%%%%%%%%%%%%
\vskip 3mm 
\noindent
{\bf  26)} coefficient of $(2200)$: \\
$|\alpha \beta \gamma \delta|$ which contains $(2200)$: $|6310|=(3100)+\underline{(2200)}
+2(2110)+3(1111)$, $|5410|=\underline{(2200)}+(2110)+2(1111)$, \\
(L.H.S.)=$[\wp_{3311}(u)]=\wp_{3311}(u)-4\wp_{31}(u)^2-2\wp_{33}(u)\wp_{11}(u)$, \\
(F.R.H.S.)=$2 \Delta_7$, \\
(S.R.H.S)=
$\displaystyle{\frac{1}{2} (\langle \langle 6310 \rangle \rangle+\langle \langle 5410 \rangle \rangle)
=\frac{1}{2}\lambda_3 \wp_{41}(u)-2\lambda_0 \wp_{44}(u) }$, \\
(T.R.H.S)=$\displaystyle{-\frac{1}{8}( \left\{ 6310 \right\}+\left\{ 5410 \right\})
 =-\frac{1}{2} \lambda_0 \lambda_8}$, \\
(Diff. Eq.): $\displaystyle{\wp_{3311}(u)-4\wp_{31}(u)^2-2\wp_{33}(u)\wp_{11}(u)
=2 \Delta_7
+\frac{1}{2}\lambda_3 \wp_{41}(u)-2\lambda_0 \wp_{44}(u)-\frac{1}{2} \lambda_0 \lambda_8}$. \\
%%%%%%%%%%%%%%%%%%%%%%%%%%%%%%%%%%%
%%%%%%%%%%%%%%%%%%%%%%%%%%%%%%
\vskip 3mm 
\noindent
{\bf  27)} coefficient of $(2111)$: \\
$|\alpha \beta \gamma \delta|$ which contains $(2111)$: $|6410|=(3200)+(3110)
+2(2210)+\underline{3(2111)}$,
$|6320|=(3110)+(2210)+\underline{3(2111)}$, 
$|5420|=(2210)+\underline{2(2111)}$, 
$|5321|=\underline{(2111)}$, \\
(L.H.S.)=$[\wp_{3222}(u)]=\wp_{3222}(u)-6\wp_{32}(u)\wp_{22}(u)$, \\
(F.R.H.S.)=$-6 \Delta_5$, \\
(S.R.H.S)=
$\displaystyle{\frac{1}{2} (3\langle \langle 6410 \rangle \rangle+3 \langle \langle 6320 \rangle \rangle
+2\langle \langle 5420 \rangle \rangle+\langle \langle 5321 \rangle \rangle)}$\\
$\displaystyle{=-\frac{3}{2}\lambda_7 \wp_{11}(u)+ \lambda_5 \wp_{31}(u)+\lambda_4 \wp_{32}(u)  
+\frac{3}{2}\lambda_3 \wp_{42}(u)-\frac{1}{2}\lambda_3 \wp_{33}(u)-2\lambda_2 \wp_{43}(u)
-\frac{3}{2} \lambda_1 \wp_{44}(u)}$, \\
(T.R.H.S)=$\displaystyle{-\frac{1}{8}(3 \left\{ 6410 \right\}+3 \left\{ 6320 \right\}
+2\left\{ 5420 \right\}+\left\{ 5321 \right\}) 
=-\frac{3}{4} \lambda_0 \lambda_9-\frac{1}{2} \lambda_1 \lambda_8-\frac{1}{4} \lambda_2 \lambda_7}$, \\
(Diff. Eq.): $\displaystyle{\wp_{3222}(u)-6\wp_{32}(u)\wp_{22}(u)
=-6 \Delta_5
-\frac{3}{2}\lambda_7 \wp_{11}(u)+ \lambda_5 \wp_{31}(u)+\lambda_4 \wp_{32}(u)  }$\\
$\displaystyle{+\frac{3}{2}\lambda_3 \wp_{42}(u)-\frac{1}{2}\lambda_3 \wp_{33}(u)-2\lambda_2 \wp_{43}(u)
-\frac{3}{2} \lambda_1 \wp_{44}(u)
-\frac{3}{4} \lambda_0 \lambda_9-\frac{1}{2} \lambda_1 \lambda_8-\frac{1}{4} \lambda_2 \lambda_7}$. \\
%%%%%%%%%%%%%%%%%%%%%%%%%%%%%%%%%%%
%%%%%%%%%%%%%%%%%%%%%%%%%%%%%%
\vskip 3mm 
\noindent
{\bf  28)} coefficient of $(2110)$: \\
$|\alpha \beta \gamma \delta|$ which contains $(2110)$: $|6310|=(3100)+(2200)
+\underline{2(2110)}+3(1111)$,
$|5410|=(2200)+\underline{(2110)}+2(1111)$, 
$|5320|=\underline{(2110)}+3(1111)$, \\
(L.H.S.)=$[\wp_{3221}(u)]=\wp_{3221}(u)-4\wp_{32}(u)\wp_{21}(u)-2\wp_{31}(u)\wp_{22}(u)$, \\
(F.R.H.S.)=$-2 \Delta_7$, \\
(S.R.H.S)=
$\displaystyle{\frac{1}{2} (2\langle \langle 6310 \rangle \rangle+\langle \langle 5410 \rangle \rangle
+\langle \langle 5320 \rangle \rangle)}$
$\displaystyle{=\lambda_4 \wp_{31}(u)+ \lambda_3 \wp_{41}(u)-\lambda_1 \wp_{43}(u)  
-\lambda_0 \wp_{44}(u)}$, \\
(T.R.H.S)=$\displaystyle{-\frac{1}{8}(2 \left\{ 6310 \right\}+\left\{ 5410 \right\}
+\left\{ 5320 \right\}) 
=-\frac{1}{2} \lambda_0 \lambda_8-\frac{1}{8} \lambda_1 \lambda_7}$, \\
(Diff. Eq.): $\displaystyle{\wp_{3221}(u)-4\wp_{32}(u)\wp_{21}(u)-2\wp_{31}(u)\wp_{22}(u)
=-2 \Delta_7
+\lambda_4 \wp_{31}(u)+ \lambda_3 \wp_{41}(u)}$\\
$\displaystyle{-\lambda_1 \wp_{43}(u)  -\lambda_0 \wp_{44}(u)
-\frac{1}{2} \lambda_0 \lambda_8-\frac{1}{8} \lambda_1 \lambda_7}$. \\
%%%%%%%%%%%%%%%%%%%%%%%%%%%%%%%%%%%
%%%%%%%%%%%%%%%%%%%%%%%%%%%%%%
\vskip 3mm 
\noindent
{\bf  29)} coefficient of $(2100)$: \\
$|\alpha \beta \gamma \delta|$ which contains $(2100)$: $|6210|=(3000)+\underline{(2100)}+(1110)$,
$|5310|=\underline{(2100)}+2(1110)$, \\
(L.H.S.)=$[\wp_{3211}(u)]=\wp_{3211}(u)-4\wp_{31}(u)\wp_{21}(u)-2\wp_{32}(u)\wp_{11}(u)$, \\
(F.R.H.S.)=0, \\
(S.R.H.S)=
$\displaystyle{\frac{1}{2} (\langle \langle 6210 \rangle \rangle+\langle \langle 5310 \rangle \rangle)}$
$\displaystyle{=\frac{1}{2}\lambda_3 \wp_{31}(u)+ \lambda_2 \wp_{41}(u)-\frac{1}{2}\lambda_1 \wp_{42}(u)  
-\lambda_0 \wp_{43}(u)}$, \\
(T.R.H.S)=$\displaystyle{-\frac{1}{8}(\left\{ 6210 \right\}+\left\{ 5310 \right\}) 
=-\frac{1}{4} \lambda_0 \lambda_7}$, \\
(Diff. Eq.): $\displaystyle{\wp_{3211}(u)-4\wp_{31}(u)\wp_{21}(u)-2\wp_{32}(u)\wp_{11}(u)
=\frac{1}{2}\lambda_3 \wp_{31}(u)+ \lambda_2 \wp_{41}(u)-\frac{1}{2}\lambda_1 \wp_{42}(u) }$\\
$\displaystyle{-\lambda_0 \wp_{43}(u)-\frac{1}{4} \lambda_0 \lambda_7}$. \\
%%%%%%%%%%%%%%%%%%%%%%%%%%%%%%%%%%%
%%%%%%%%%%%%%%%%%%%%%%%%%%%%%%
\vskip 3mm 
\noindent
{\bf  30)} coefficient of $(2000)$: \\
$|\alpha \beta \gamma \delta|$ which contains $(2000)$: $|52310|=\underline{(2000)}+(1100)$, \\
(L.H.S.)=$[\wp_{3111}(u)]=\wp_{3111}(u)-6\wp_{31}(u)\wp_{11}(u)$, \\
(F.R.H.S.)=0, \\
(S.R.H.S)=
$\displaystyle{\frac{1}{2} \langle \langle 5210 \rangle \rangle
=\lambda_2 \wp_{31}(u)+\frac{3}{2} \lambda_1 \wp_{41}(u)-\frac{1}{2}\lambda_1 \wp_{32}(u)  
-3\lambda_0 \wp_{42}(u)+\lambda_0 \wp_{33}(u)}$, \\
(T.R.H.S)=$\displaystyle{-\frac{1}{8} \left\{ 5210 \right\} =0}$, \\
(Diff. Eq.): $\displaystyle{\wp_{3111}(u)-6\wp_{31}(u)\wp_{11}(u)
=\lambda_2 \wp_{31}(u)+\frac{3}{2} \lambda_1 \wp_{41}(u)-\frac{1}{2}\lambda_1 \wp_{32}(u)  
-3\lambda_0 \wp_{42}(u)+\lambda_0 \wp_{33}(u) }$. \\
%%%%%%%%%%%%%%%%%%%%%%%%%%%%%%%%%%%
%%%%%%%%%%%%%%%%%%%%%%%%%%%%%%
\vskip 3mm 
\noindent
{\bf  31)} coefficient of $(1111)$: \\
$|\alpha \beta \gamma \delta|$ which contains $(1111)$: $|6310|=(3100)+(2200)
+2(2110)+\underline{3(1111)}$,
$|5410|=(2200)+(2110)+\underline{2(1111)}$, 
$|5320|=(2110)+\underline{3(1111)}$, 
$|4321|=\underline{(1111)}$, \\
(L.H.S.)=$[\wp_{2222}(u)]=\wp_{2222}(u)-6\wp_{22}(u)^2$, \\
(F.R.H.S.)=$12 \Delta_7$, \\
(S.R.H.S)=
$\displaystyle{\frac{1}{2} (3\langle \langle 6310 \rangle \rangle+2 \langle \langle 5410 \rangle \rangle
+3\langle \langle 5320 \rangle \rangle+\langle \langle 4321 \rangle \rangle)}$\\
$\displaystyle{=-3\lambda_6 \wp_{11}(u)+ \lambda_5 \wp_{21}(u)+\lambda_4 \wp_{22}(u)  
+\lambda_3 \wp_{32}(u)-3\lambda_2 \wp_{33}(u)+4\lambda_2 \wp_{42}(u)}$\\
$\displaystyle{-3 \lambda_1 \wp_{43}(u)-3 \lambda_0 \wp_{44}(u)}$, \\
(T.R.H.S)=$\displaystyle{-\frac{1}{8}(3 \left\{ 6310 \right\}+2 \left\{ 5410 \right\}
+3\left\{ 5320 \right\}+\left\{ 4321 \right\}) 
=-\lambda_0 \lambda_8-\frac{3}{8} \lambda_1 \lambda_7-\frac{1}{2} \lambda_2 \lambda_6
+\frac{1}{8} \lambda_3 \lambda_5}$, \\
(Diff. Eq.): $\displaystyle{\wp_{2222}(u)-6\wp_{22}(u)^2
=12 \Delta_7
-3\lambda_6 \wp_{11}(u)+ \lambda_5 \wp_{21}(u)+\lambda_4 \wp_{22}(u)  
+\lambda_3 \wp_{32}(u)}$\\
$\displaystyle{+4\lambda_2 \wp_{42}(u)-3\lambda_2 \wp_{33}(u)-3 \lambda_1 \wp_{43}(u)
-3 \lambda_0 \wp_{44}(u)
-\lambda_0 \lambda_8-\frac{3}{8} \lambda_1 \lambda_7-\frac{1}{2} \lambda_2 \lambda_6
+\frac{1}{8} \lambda_3 \lambda_5}$. \\
%%%%%%%%%%%%%%%%%%%%%%%%%%%%%%%%%%%
%%%%%%%%%%%%%%%%%%%%%%%%%%%%%%
\vskip 3mm 
\noindent
{\bf  32)} coefficient of $(1110)$: \\
$|\alpha \beta \gamma \delta|$ which contains $(1110)$: $|6210|=(3000)+(2100)+\underline{(1110)}$,
$|5310|=(2100)+\underline{2(1110)}$, $|4320|=\underline{(1110)}$, \\
(L.H.S.)=$[\wp_{2221}(u)]=\wp_{2221}(u)-6\wp_{22}(u)\wp_{21}(u)$, \\
(F.R.H.S.)=0, \\
(S.R.H.S)=
$\displaystyle{\frac{1}{2} (\langle \langle 6210 \rangle \rangle+2 \langle \langle 5310 \rangle \rangle
+\langle \langle 4320 \rangle \rangle)}$\\
$\displaystyle{=-\frac{1}{2}\lambda_5 \wp_{11}(u)+ \lambda_4 \wp_{21}(u)+\lambda_3 \wp_{31}(u)  
+\lambda_2 \wp_{41}(u)+\frac{3}{2}\lambda_1 \wp_{42}(u)-\frac{3}{2}\lambda_1 \wp_{33}(u)
-3\lambda_0 \wp_{43}(u)}$, \\
(T.R.H.S)=$\displaystyle{-\frac{1}{8}( \left\{ 6210 \right\}+2 \left\{ 5310 \right\}
+\left\{ 4320 \right\}) 
=-\frac{1}{2} \lambda_0 \lambda_7-\frac{1}{4} \lambda_1 \lambda_6}$, \\
(Diff. Eq.): $\displaystyle{\wp_{2221}(u)-6\wp_{22}(u)\wp_{21}(u)
=-\frac{1}{2}\lambda_5 \wp_{11}(u)+ \lambda_4 \wp_{21}(u)+\lambda_3 \wp_{31}(u)  
+\lambda_2 \wp_{41}(u)}$\\
$\displaystyle{+\frac{3}{2}\lambda_1 \wp_{42}(u)-\frac{3}{2}\lambda_1 \wp_{33}(u)
-3\lambda_0 \wp_{43}(u)
-\frac{1}{2} \lambda_0 \lambda_7-\frac{1}{4} \lambda_1 \lambda_6}$. \\
%%%%%%%%%%%%%%%%%%%%%%%%%%%%%%%%%%%
%%%%%%%%%%%%%%%%%%%%%%%%%%%%%%
\vskip 3mm 
\noindent
{\bf  33)} coefficient of $(1100)$: \\
$|\alpha \beta \gamma \delta|$ which contains $(1100)$: $|5210|=(2000)+\underline{(1100)}$,
 $|4310|=\underline{(1100)}$, \\
(L.H.S.)=$[\wp_{2211}(u)]=\wp_{2211}(u)-4\wp_{21}(u)^2-2\wp_{22}(u)\wp_{11}(u)$, \\
(F.R.H.S.)=0, \\
S.R.H.S)=
$\displaystyle{\frac{1}{2} (\langle \langle 5210 \rangle \rangle+\langle \langle 4310 \rangle \rangle)
=\frac{1}{2}\lambda_3 \wp_{21}(u)+ \lambda_2 \wp_{31}(u)+\frac{3}{2}\lambda_1 \wp_{41}(u)  
-\frac{1}{2}\lambda_1 \wp_{32}(u)}$\\
$\displaystyle{+\lambda_0 \wp_{42}(u)-2\lambda_0 \wp_{33}(u)}$, \\
(T.R.H.S)=$\displaystyle{-\frac{1}{8}( \left\{ 5210 \right\}+ \left\{ 4310 \right\}) 
=-\frac{1}{2} \lambda_0 \lambda_6}$, \\
(Diff. Eq.): $\displaystyle{\wp_{2211}(u)-4\wp_{21}(u)^2-2\wp_{22}(u)\wp_{11}(u)
=\frac{1}{2}\lambda_3 \wp_{21}(u)+ \lambda_2 \wp_{31}(u)}$\\
$\displaystyle{+\frac{3}{2}\lambda_1 \wp_{41}(u) -\frac{1}{2}\lambda_1 \wp_{32}(u)
+\lambda_0 \wp_{42}(u)-2\lambda_0 \wp_{33}(u)
-\frac{1}{2} \lambda_0 \lambda_6}$. \\
%%%%%%%%%%%%%%%%%%%%%%%%%%%%%%%%%%%
%%%%%%%%%%%%%%%%%%%%%%%%%%%%%%
\vskip 3mm 
\noindent
{\bf  34)} coefficient of $(1000)$: \\
$|\alpha \beta \gamma \delta|$ which contains $(1000)$: $|4210|=\underline{(1000)}$, \\
(L.H.S.)=$[\wp_{2111}(u)]=\wp_{2111}(u)-6\wp_{21}(u)\wp_{11}(u)$, \\
(F.R.H.S.)=0, \\
S.R.H.S)=
$\displaystyle{\frac{1}{2} \langle \langle 4210 \rangle \rangle
=\lambda_2 \wp_{21}(u)+\frac{3}{2}\lambda_1 \wp_{31}(u)-\frac{1}{2}\lambda_1 \wp_{22}(u)  
+3\lambda_0 \wp_{41}(u)-2\lambda_0 \wp_{32}(u)}$, \\
(T.R.H.S)=$\displaystyle{-\frac{1}{8} \left\{ 4210 \right\} 
=-\frac{1}{4} \lambda_0 \lambda_5}$, \\
(Diff. Eq.): $\displaystyle{\wp_{2111}(u)-6\wp_{21}(u)\wp_{11}(u)
=\lambda_2 \wp_{21}(u)+\frac{3}{2}\lambda_1 \wp_{31}(u)-\frac{1}{2}\lambda_1 \wp_{22}(u)  
+3\lambda_0 \wp_{41}(u)}$\\
$\displaystyle{-2\lambda_0 \wp_{32}(u)
-\frac{1}{4} \lambda_0 \lambda_5}$. \\
%%%%%%%%%%%%%%%%%%%%%%%%%%%%%%%%%%%
%%%%%%%%%%%%%%%%%%%%%%%%%%%%%%
\vskip 3mm 
\noindent
{\bf  35)} coefficient of $(0000)=1$: \\
$|\alpha \beta \gamma \delta|$ which contains $(0000)=1$: $|3210|=\underline{(0000)}=1$, \\
(L.H.S.)=$[\wp_{1111}(u)]=\wp_{1111}(u)-6\wp_{11}(u)^2$, \\
(F.R.H.S.)=0, \\
S.R.H.S)=
$\displaystyle{\frac{1}{2} \langle \langle 3210 \rangle \rangle
=\lambda_2 \wp_{11}(u)+\lambda_1 \wp_{21}(u)+4\lambda_0 \wp_{31}(u)  
-3\lambda_0 \wp_{22}(u)}$, \\
(T.R.H.S)=$\displaystyle{-\frac{1}{8} \left\{ 3210 \right\} 
=-\frac{1}{2} \lambda_0 \lambda_4+\frac{1}{8} \lambda_1 \lambda_3}$, \\
(Diff. Eq.): $\displaystyle{\wp_{1111}(u)-6\wp_{11}(u)^2
=\lambda_2 \wp_{11}(u)+\lambda_1 \wp_{21}(u)+4\lambda_0 \wp_{31}(u)  
-3\lambda_0 \wp_{22}(u)}$\\
$\displaystyle{-\frac{1}{2} \lambda_0 \lambda_4+\frac{1}{8} \lambda_1 \lambda_3}$. \\
%%%%%%%%%%%%%%%%%%%%%%%%%%%%%%%%%%%
%%%%%%%%%%%%%%%%%%%%%%%%%%%%%%%%%%%%%%%%%%%%%%%%%%%%%%%%%%%%%%%%%
%\newpage
%
\vskip 3mm
Therefore, we have the following differential equations
\begin{align}
 {\rm \bf 1)}\quad & \wp_{4444}(u)-6\wp_{44}(u)^2=\lambda_9 \wp_{43}(u) +\lambda_8 \wp_{44}(u)
+\frac{1}{8} \lambda_7 \lambda_9 ,
\label{4e46}\\
{\rm \bf  2)}\quad & \wp_{4443}(u)-6\wp_{44}(u) \wp_{43}(u)=\frac{3}{2} \lambda_9 \wp_{42}(u) 
-\frac{1}{2}\lambda_9 \wp_{33}(u)+\lambda_8 \wp_{43}(u),
\label{4e47}\\
{\rm \bf  3)}\quad & \wp_{4442}(u)-6\wp_{44}(u) \wp_{42}(u)=\frac{3}{2} \lambda_9 \wp_{41}(u) 
-\frac{1}{2}\lambda_9 \wp_{32}(u)+\lambda_8 \wp_{42}(u) ,
\label{4e48}\\
{\rm \bf  4)}\quad &\wp_{4441}(u)-6\wp_{44}(u) \wp_{41}(u)=-\frac{1}{2} \lambda_9 \wp_{31}(u) 
+\lambda_8 \wp_{41}(u) ,
\label{4e49}\\
{\rm \bf  5)}\quad & \wp_{4433}(u)-4\wp_{43}(u)^2-2\wp_{44}(u) \wp_{33}(u)
=\frac{3}{2} \lambda_9 \wp_{41}(u)-\frac{1}{2}\lambda_9\wp_{32}(u) +\lambda_8 \wp_{42}(u)
\nonumber\\
&+\frac{1}{2} \lambda_7 \wp_{43}(u) ,
\label{4e50}\\
{\rm \bf  6)}\quad &\wp_{4432}(u)-4\wp_{43}(u)\wp_{42}(u)-2\wp_{44}(u) \wp_{32}(u)
=-\frac{1}{2} \lambda_9 \wp_{31}(u)+\lambda_8\wp_{41}(u) 
\nonumber\\
&+\frac{1}{2}\lambda_7\wp_{42}(u) ,
\label{4e51}\\
{\rm \bf  7)}\quad  & \wp_{4431}(u)-4\wp_{43}(u)\wp_{41}(u)-2\wp_{44}(u) \wp_{31}(u)
=\frac{1}{2} \lambda_7 \wp_{41}(u) ,
\label{4e52}\\
{\rm \bf  8)}\quad &\wp_{4422}(u)-4\wp_{42}(u)^2-2\wp_{44}(u) \wp_{22}(u)
=2 \Delta_2+\frac{1}{2} \lambda_7 \wp_{41}(u) ,
\label{4e53}\\
{\rm \bf  9)}\quad &\wp_{4421}(u)-4\wp_{42}(u)\wp_{41}(u)-2\wp_{44}(u) \wp_{21}(u)
=-2 \Delta_3 ,
\label{4e54}\\
{\rm \bf 10)}\quad & \wp_{4411}(u)-4\wp_{41}(u)^2-2\wp_{44}(u) \wp_{11}(u)
=-2 \Delta_4 ,
\label{4e55}\\
{\rm \bf 11)}\quad & \wp_{4333}(u)-6\wp_{43}(u)\wp_{33}(u)
=\frac{3}{2} \lambda_9 \wp_{31}(u)-\frac{3}{2} \lambda_9 \wp_{22}(u)+\lambda_8 \wp_{41}(u)
+\lambda_7 \wp_{42}(u)
\nonumber\\
&+\lambda_6 \wp_{43}(u)-\frac{1}{2} \lambda_5 \wp_{44}(u) 
-\frac{1}{4} \lambda_4 \lambda_9 ,
\label{4e56}\\
{\rm \bf 12)}\quad &\wp_{4332}(u)-4\wp_{43}(u)\wp_{32}(u)-2\wp_{42}(u)\wp_{33}(u)
=-2 \Delta_2 -\lambda_9 \wp_{21}(u)+\lambda_7 \wp_{41}(u)
\nonumber\\
&+\lambda_6 \wp_{42}(u) 
-\frac{1}{8} \lambda_3 \lambda_9 ,
\label{4e57}\\
{\rm \bf 13)}\quad & \wp_{4331}(u)-4\wp_{43}(u)\wp_{31}(u)-2\wp_{41}(u)\wp_{33}(u)
=2 \Delta_3 +\frac{1}{2}\lambda_9 \wp_{11}(u)+\lambda_6 \wp_{41}(u) ,
\label{4e58}\\
{\rm \bf 14)}\quad &\wp_{4322}(u)-4\wp_{42}(u)\wp_{32}(u)-2\wp_{43}(u)\wp_{22}(u)
=2 \Delta_3 -\lambda_9 \wp_{11}(u)+\lambda_6 \wp_{41}(u)
\nonumber\\
&+\frac{1}{2}\lambda_5 \wp_{42}(u)-\frac{1}{4} \lambda_2 \lambda_9 ,
\label{4e59}\\
{\rm \bf 15)}\quad & \wp_{4321}(u)-2\wp_{43}(u)\wp_{21}(u)-2\wp_{42}(u)\wp_{31}(u)
-2\wp_{41}(u)\wp_{32}(u)=\frac{1}{2}\lambda_5 \wp_{41}(u)
\nonumber\\
&-\frac{1}{8} \lambda_1 \lambda_9 ,
\label{4e60}\\
{\rm \bf 16)}\quad & \wp_{4311}(u)-4\wp_{41}(u)\wp_{31}(u)-2\wp_{43}(u)\wp_{11}(u)
=-2 \Delta_5-\frac{1}{4} \lambda_0 \lambda_9 ,
\label{4e61}\\
{\rm \bf 17)}\quad &\wp_{4222}(u)-6\wp_{42}(u)\wp_{22}(u)
=6 \Delta_4 +\lambda_5 \wp_{41}(u)+\lambda_4 \wp_{42}(u)-\frac{1}{2}\lambda_3 \wp_{43}(u)
\nonumber\\
&+\lambda_2 \wp_{44}(u) 
-\frac{1}{4} \lambda_1 \lambda_9 ,
\label{4e62}\\
{\rm \bf 18)}\quad & \wp_{4221}(u)-4\wp_{42}(u)\wp_{21}(u)-2\wp_{41}(u)\wp_{22}(u)
=2 \Delta_5 +\lambda_4 \wp_{41}(u)+\frac{1}{2}\lambda_1 \wp_{44}(u)
\nonumber\\
&-\frac{1}{4} \lambda_0 \lambda_9, 
\label{4e63}\\
{\rm \bf 19)}\quad & \wp_{4211}(u)-4\wp_{41}(u)\wp_{21}(u)-2\wp_{42}(u)\wp_{11}(u)
= \frac{1}{2}\lambda_3 \wp_{41}(u)+\lambda_0 \wp_{44}(u) , 
\label{4e64}\\
{\rm \bf 20)}\quad &\wp_{4111}(u)-6\wp_{41}(u)\wp_{11}(u)
=\lambda_2 \wp_{41}(u)-\frac{1}{2}\lambda_1 \wp_{42}(u)+\lambda_0 \wp_{43}(u) , 
\label{4e65}\\
{\rm \bf 21)}\quad &\wp_{3333}(u)-6\wp_{33}(u)^2=12 \Delta_2
-3\lambda_9 \wp_{21}(u)+4 \lambda_8 \wp_{31}(u)-3\lambda_8 \wp_{22}(u)  
+\lambda_7 \wp_{32}(u)
\nonumber\\
&+\lambda_6 \wp_{33}(u)+\lambda_5 \wp_{43}(u)-3\lambda_4 \wp_{44}(u)
+\frac{1}{8} \lambda_5 \lambda_7
-\frac{1}{2} \lambda_4 \lambda_8-\frac{3}{8} \lambda_3 \lambda_9 ,
\label{4e66}\\
{\rm \bf 22)}\quad &\wp_{3332}(u)-6\wp_{33}(u)\wp_{32}(u)=-6 \Delta_3
-\frac{3}{2}\lambda_9 \wp_{11}(u)-2 \lambda_8 \wp_{21}(u)+\frac{3}{2}\lambda_7 \wp_{31}(u)
\nonumber\\
&-\frac{1}{2}\lambda_7 \wp_{22}(u)+\lambda_6 \wp_{32}(u)+\lambda_5 \wp_{42}(u)
-\frac{3}{2}\lambda_3 \wp_{44}(u)
-\frac{1}{2} \lambda_2 \lambda_9-\frac{1}{4} \lambda_3 \lambda_8 ,
\label{4e67}\\
{\rm \bf 23)}\quad & \wp_{3331}(u)-6\wp_{33}(u)\wp_{31}(u)=6 \Delta_4
+\lambda_8 \wp_{11}(u)-\frac{1}{2} \lambda_7 \wp_{21}(u)+\lambda_6 \wp_{31}(u) 
\nonumber\\
&+\lambda_5 \wp_{41}(u) -\frac{1}{4} \lambda_1 \lambda_9 ,
\label{4e68}\\
{\rm \bf 24)}\quad & \wp_{3322}(u)-4\wp_{32}(u)^2-2  \wp_{33}(u) \wp_{22}(u)
=-6 \Delta_4-2\lambda_8 \wp_{11}(u)-\frac{1}{2} \lambda_7 \wp_{21}(u)
\nonumber\\
&+\lambda_6 \wp_{31}(u)
+\frac{1}{2}\lambda_5 \wp_{41}(u)+\frac{1}{2}\lambda_5 \wp_{32}(u)
+\lambda_4 \wp_{42}(u)-\frac{1}{2} \lambda_3 \wp_{43}(u)-2\lambda_2 \wp_{44}(u)
\nonumber\\
&-\frac{1}{2} \lambda_1 \lambda_9-\frac{1}{2} \lambda_2 \lambda_8 ,
\label{4e69}\\
{\rm \bf 25)}\quad & \wp_{3321}(u)-4\wp_{32}(u)\wp_{31}(u)-2\wp_{33}(u)\wp_{21}(u)
=2 \Delta_5
+\frac{1}{2}\lambda_5 \wp_{31}(u)+\lambda_4 \wp_{41}(u)
\nonumber\\
&-\lambda_1 \wp_{44}(u)
-\frac{1}{4} \lambda_1 \lambda_8-\frac{1}{2} \lambda_0 \lambda_9 , 
\label{4e70}\\
{\rm \bf 26)}\quad &\wp_{3311}(u)-4\wp_{31}(u)^2-2\wp_{33}(u)\wp_{11}(u)
=2 \Delta_7
+\frac{1}{2}\lambda_3 \wp_{41}(u)-2\lambda_0 \wp_{44}(u)-\frac{1}{2} \lambda_0 \lambda_8 ,
\label{4e71}\\
{\rm \bf 27)}\quad &\wp_{3222}(u)-6\wp_{32}(u)\wp_{22}(u)
=-6 \Delta_5
-\frac{3}{2}\lambda_7 \wp_{11}(u)+ \lambda_5 \wp_{31}(u)+\lambda_4 \wp_{32}(u) 
+\frac{3}{2}\lambda_3 \wp_{42}(u)
\nonumber\\
&-\frac{1}{2}\lambda_3 \wp_{33}(u)-2\lambda_2 \wp_{43}(u)
-\frac{3}{2} \lambda_1 \wp_{44}(u)
-\frac{3}{4} \lambda_0 \lambda_9-\frac{1}{2} \lambda_1 \lambda_8-\frac{1}{4} \lambda_2 \lambda_7 , 
\label{4e72}\\
{\rm \bf 28)}\quad & \wp_{3221}(u)-4\wp_{32}(u)\wp_{21}(u)-2\wp_{31}(u)\wp_{22}(u)
=-2 \Delta_7
+\lambda_4 \wp_{31}(u)+ \lambda_3 \wp_{41}(u)
\nonumber\\
&-\lambda_1 \wp_{43}(u)  -\lambda_0 \wp_{44}(u)
-\frac{1}{2} \lambda_0 \lambda_8-\frac{1}{8} \lambda_1 \lambda_7,
\label{4e73}\\
{\rm \bf 29)}\quad & \wp_{3211}(u)-4\wp_{31}(u)\wp_{21}(u)-2\wp_{32}(u)\wp_{11}(u)
=\frac{1}{2}\lambda_3 \wp_{31}(u)+ \lambda_2 \wp_{41}(u)
\nonumber\\
&-\frac{1}{2}\lambda_1 \wp_{42}(u)
-\lambda_0 \wp_{43}(u)-\frac{1}{4} \lambda_0 \lambda_7 , 
\label{4e74}\\
{\rm \bf 30)}\quad &\wp_{3111}(u)-6\wp_{31}(u)\wp_{11}(u)
=\lambda_2 \wp_{31}(u)+\frac{3}{2} \lambda_1 \wp_{41}(u)-\frac{1}{2}\lambda_1 \wp_{32}(u) \nonumber\\
&-3\lambda_0 \wp_{42}(u)+\lambda_0 \wp_{33}(u)  , 
\label{4e75}\\
{\rm \bf 31)}\quad & \wp_{2222}(u)-6\wp_{22}(u)^2
=12 \Delta_7
-3\lambda_6 \wp_{11}(u)+ \lambda_5 \wp_{21}(u)+\lambda_4 \wp_{22}(u)  
+\lambda_3 \wp_{32}(u)
\nonumber\\
&+4\lambda_2 \wp_{42}(u)-3\lambda_2 \wp_{33}(u)-3 \lambda_1 \wp_{43}(u)
-3 \lambda_0 \wp_{44}(u)-\lambda_0 \lambda_8-\frac{3}{8} \lambda_1 \lambda_7
\nonumber\\
&-\frac{1}{2} \lambda_2 \lambda_6+\frac{1}{8} \lambda_3 \lambda_5 , 
\label{4e76}\\
{\rm \bf 32)}\quad & \wp_{2221}(u)-6\wp_{22}(u)\wp_{21}(u)
=-\frac{1}{2}\lambda_5 \wp_{11}(u)+ \lambda_4 \wp_{21}(u)+\lambda_3 \wp_{31}(u)  
+\lambda_2 \wp_{41}(u)
\nonumber\\
&+\frac{3}{2}\lambda_1 \wp_{42}(u)-\frac{3}{2}\lambda_1 \wp_{33}(u)
-3\lambda_0 \wp_{43}(u)
-\frac{1}{2} \lambda_0 \lambda_7-\frac{1}{4} \lambda_1 \lambda_6 ,
\label{4e77}\\
{\rm \bf 33)}\quad &\wp_{2211}(u)-4\wp_{21}(u)^2-2\wp_{22}(u)\wp_{11}(u)
=\frac{1}{2}\lambda_3 \wp_{21}(u)+ \lambda_2 \wp_{31}(u)
+\frac{3}{2}\lambda_1 \wp_{41}(u) 
\nonumber\\
&-\frac{1}{2}\lambda_1 \wp_{32}(u)+\lambda_0 \wp_{42}(u)-2\lambda_0 \wp_{33}(u)
-\frac{1}{2} \lambda_0 \lambda_6 ,
\label{4e78}
\\
{\rm \bf 34)}\quad & \wp_{2111}(u)-6\wp_{21}(u)\wp_{11}(u)
=\lambda_2 \wp_{21}(u)+\frac{3}{2}\lambda_1 \wp_{31}(u)-\frac{1}{2}\lambda_1 \wp_{22}(u) 
\nonumber\\
&+3\lambda_0 \wp_{41}(u)-2\lambda_0 \wp_{32}(u)
-\frac{1}{4} \lambda_0 \lambda_5,
\label{4e79}\\
{\rm \bf 35)}\quad &\wp_{1111}(u)-6\wp_{11}(u)^2
=\lambda_2 \wp_{11}(u)+\lambda_1 \wp_{21}(u)+4\lambda_0 \wp_{31}(u)  
-3\lambda_0 \wp_{22}(u)
\nonumber\\
&-\frac{1}{2} \lambda_0 \lambda_4+\frac{1}{8} \lambda_1 \lambda_3  ,
\label{4e80}
\end{align}
where 
\begin{align}
\Delta_2   &=(-\wp_{44}(u)\wp_{22}(u)+\wp_{42}(u)^2)+(\wp_{44}(u)\wp_{31}(u)-\wp_{43}(u)\wp_{41}(u))
\nonumber\\
&+(\wp_{43}(u)\wp_{32}(u)-\wp_{42}\wp_{33}(u)),  
\label{4e81}\\
\Delta_3   &=(\wp_{44(u)}\wp_{21}(u)-\wp_{42}(u)\wp_{41}(u))+(-\wp_{43}(u)\wp_{31}(u)+\wp_{41}(u)\wp_{33}),
\label{4e82}\\
\Delta_4   &=(\wp_{44}(u)\wp_{11}(u)-\wp_{41}(u)^2)+(-\wp_{42}(u)\wp_{31}(u)+\wp_{41}(u)\wp_{32}(u)), 
\label{4e83}\\
\Delta_5   &=(\wp_{43}(u)\wp_{11}(u)-\wp_{41}(u)\wp_{31}(u))+(-\wp_{42}(u)\wp_{21}(u)+\wp_{41}(u)\wp_{22}(u)),  
\label{4e84}\\
\Delta_7   &=(-\wp_{33}(u)\wp_{11}(u)+\wp_{31}(u)^2)+(\wp_{42}(u)\wp_{11}(u)-\wp_{41}(u)\wp_{21}(u))
\nonumber\\
&+(\wp_{32}(u)\wp_{21}(u)-\wp_{31}(u)\wp_{22}(u)) .
\label{4e85}
\end{align}
$\Delta_i$ are the sum of the special form of 
$\zeta_i(u)\Big(\overleftarrow{\partial_k}\overrightarrow{\partial_{\ell}} 
-\overrightarrow{\partial_k}\overleftarrow{\partial_{\ell}}\Big)\zeta_j(u)$
by using $\zeta_i(u)$ functions.

These differential equations agree with those of our previous paper~\cite{Hayashi}. 
In order to compare with our previous paper, we must replace 
$\wp_{ij}(u) \rightarrow 4\wp_{ij}(u)$ and $\wp_{ijk \ell}(u) \rightarrow 4\wp_{ijk \ell}(u)$ 
in our previous paper.

%%%%%%%%%%%%%%%%%%%%%%%%%%%%%%%%%%%%%%%%%%%%%%%%%%%%%%%%
\section{Summary} 
\setcounter{equation}{0}

The differential equation of the Weierstrass $\wp$ function is the differential 
equation of the static KdV equation. Then,  a set of coupled differential equations for 
the genus $g$ hyperelliptic functions can be considered 
as the $g$-dimensional generalization of the static KdV equation.
Thus,  a set of such differential equations is integrable, and the $\sigma$ function 
play the potential for such integrable differential equations.
In this way, a set of differential equations for genus $g$ hyperelliptic $\wp_{ij}$ functions 
are the special but quite important ones.

Though there are various approach to obtain a set of differential equation for genus 
$g$ hyperelliptic $\wp_{ij}$ functions, Baker's approach seems to be the most systematic one.
Thus, we review the Baker's method to obtain differential equations of the 
general genus hyperelliptic $\wp_{ij}$ functions. As Baker demonstrated his method in genus 
three case, we demonstrate in genus four case, which agree with our previous result.

%\newpage
%%%%%%%%%%%%%%%%%%%%%%%%%%%%%%%%%%%%%%%%%%%%%%%%%%%%%%%%%%%
%%%%%%%%%%%%%%%%%%%%%%%%%%%%%%%%%%%%%%%%%%%%%%%%%%%%%%%%%%%  

\end{document}